\documentclass[11pt,a4paper]{article}
\usepackage{epsfig,latexsym,amsfonts,amssymb,amsmath,amscd,
epic,graphics,theorem,epstopdf, placeins}%,pxfonts}%,pxfonts}
\theoremstyle{plain}
\newtheorem{lemma}{Lemma}[section]
\newtheorem{proposition}[lemma]{Proposition}
\newtheorem{remark}[lemma]{Remark}
\newtheorem{example}[lemma]{Example}
\newtheorem{theorem}[lemma]{Theorem}
\newtheorem{definition}[lemma]{Definition}

{\theorembodyfont{\rmfamily}  \font\ncsc=cmcsc10
 \font\ntt=cmtt12

\begin{document}
\newcommand{\pperp}{\hbox{$\perp\hskip-6pt\perp$}}
\newcommand{\ssim}{\hbox{$\hskip-2pt\sim$}}
\newcommand{\aleq}{{\ \stackrel{3}{\le}\ }}
\newcommand{\ageq}{{\ \stackrel{3}{\ge}\ }}
\newcommand{\aeq}{{\ \stackrel{3}{=}\ }}
\newcommand{\bleq}{{\ \stackrel{n}{\le}\ }}
\newcommand{\bgeq}{{\ \stackrel{n}{\ge}\ }}
\newcommand{\beq}{{\ \stackrel{n}{=}\ }}
\newcommand{\cleq}{{\ \stackrel{2}{\le}\ }}
\newcommand{\cgeq}{{\ \stackrel{2}{\ge}\ }}
\newcommand{\ceq}{{\ \stackrel{2}{=}\ }}
\newcommand{\fm}{\mathfrak{m}}
\newcommand{\N}{{\mathbb N}}
\newcommand{\A}{{\mathbb A}}
\newcommand{\K}{{\mathbb K}}
\newcommand{\Z}{{\mathbb Z}}\newcommand{\F}{{\mathbf F}}
\newcommand{\R}{{\mathbb R}}
\newcommand{\C}{{\mathbb C}}
\newcommand{\Q}{{\mathbb Q}}
\newcommand{\PP}{{\mathbb P}}
\newcommand{\mnote}{\marginpar}\newcommand{\red}{{\operatorname{red}}}
\newcommand{\Pic}{{\operatorname{Pic}}}\newcommand{\Sym}{{\operatorname{Sym}}}
\newcommand{\oeps}{{\overline\eps}}\newcommand{\Div}{{\operatorname{Div}}}
\newcommand{\oDel}{{\widetilde\Del}}
\newcommand{\real}{{\operatorname{Re}}}\newcommand{\Aut}{{\operatorname{Aut}}}
\newcommand{\conv}{{\operatorname{conv}}}\newcommand{\Ima}{{\operatorname{Im}}}
\newcommand{\Span}{{\operatorname{Span}}}
\newcommand{\Ker}{{\operatorname{Ker}}}
\newcommand{\Ann}{{\operatorname{Ann}}}
\newcommand{\Fix}{{\operatorname{Fix}}}
\newcommand{\sign}{{\operatorname{sign}}}
\newcommand{\Tors}{{\operatorname{Tors}}}
\newcommand{\Card}{{\operatorname{Card}}}
\newcommand{\alg}{{\operatorname{alg}}}\newcommand{\ord}{{\operatorname{ord}}}
\newcommand{\oi}{{\overline i}}
\newcommand{\oj}{{\overline j}}
\newcommand{\ob}{{\overline b}}
\newcommand{\os}{{\overline s}}
\newcommand{\oa}{{\overline a}}
\newcommand{\oy}{{\overline y}}
\newcommand{\ow}{{\overline w}}
\newcommand{\ot}{{\overline t}}
\newcommand{\oz}{{\overline z}}
\newcommand{\eps}{{\varepsilon}}
\newcommand{\proofend}{\hfill$\Box$\bigskip}
\newcommand{\Int}{{\operatorname{Int}}}
\newcommand{\pr}{{\operatorname{pr}}}
\newcommand{\Hom}{{\operatorname{Hom}}}
\newcommand{\Ev}{{\operatorname{Ev}}}
\newcommand{\im}{{\operatorname{Im}}}\newcommand{\br}{{\operatorname{br}}}
\newcommand{\sk}{{\operatorname{sk}}}\newcommand{\DP}{{\operatorname{DP}}}
\newcommand{\const}{{\operatorname{const}}}
\newcommand{\Sing}{{\operatorname{Sing}}\hskip0.06cm}
\newcommand{\conj}{{\operatorname{Conj}}}
\newcommand{\Cl}{{\operatorname{Cl}}}
\newcommand{\Crit}{{\operatorname{Crit}}}
\newcommand{\Ch}{{\operatorname{Ch}}}
\newcommand{\discr}{{\operatorname{discr}}}
\newcommand{\Tor}{{\operatorname{Tor}}}
\newcommand{\Conj}{{\operatorname{Conj}}}
\newcommand{\vol}{{\operatorname{vol}}}
\newcommand{\defect}{{\operatorname{def}}}
\newcommand{\codim}{{\operatorname{codim}}}
\newcommand{\tmu}{{\C\mu}}
\newcommand{\ov}{{\overline v}}
\newcommand{\ox}{{\overline{x}}}
\newcommand{\bw}{{\boldsymbol w}}
\newcommand{\bv}{{\boldsymbol v}}
\newcommand{\bn}{{\boldsymbol n}}
\newcommand{\bx}{{\boldsymbol x}}
\newcommand{\bd}{{\boldsymbol d}}
\newcommand{\bz}{{\boldsymbol z}}
\newcommand{\bp}{{\boldsymbol p}}
\newcommand{\tet}{{\theta}}
\newcommand{\Del}{{\Delta}}
\newcommand{\bet}{{\beta}}
\newcommand{\kap}{{\kappa}}
\newcommand{\del}{{\delta}}
\newcommand{\sig}{{\sigma}}
\newcommand{\alp}{{\alpha}}
\newcommand{\Sig}{{\Sigma}}
\newcommand{\Gam}{{\Gamma}}
\newcommand{\gam}{{\gamma}}\newcommand{\idim}{{\operatorname{idim}}}
\newcommand{\Lam}{{\Lambda}}
\newcommand{\lam}{{\lambda}}
\newcommand{\SC}{{SC}}
\newcommand{\MC}{{MC}}
\newcommand{\nek}{{,...,}}
\newcommand{\cim}{{c_{\mbox{\rm im}}}}
\newcommand{\clM}{\tilde{M}}
\newcommand{\clV}{\bar{V}}

\title{Relative enumerative invariants of real nodal del Pezzo surfaces}
\author{Ilia Itenberg,
\and Viatcheslav Kharlamov, \and Eugenii Shustin}
\date{}
\maketitle
\begin{abstract}
The surfaces considered are real, rational and have
a unique smooth real $(-2)$-curve.
Their canonical class $K$ is
strictly negative on any other irreducible curve
in the surface and
$K^2>0$. For
surfaces satisfying these assumptions,
we suggest
a certain signed count of real rational curves
that belong to a given divisor class
and are
simply
tangent to the $(-2)$-curve at each intersection point. We prove that this count provides
a number which
depends neither
on the point constraints
nor
on deformation of the surface preserving the real structure and
the $(-2)$-curve.
\end{abstract}

\medskip

{\bf MSC-2010 classification:} Primary 14N10, Secondary 14J26, 14P05

\medskip

{\it \hskip1.5in  Les \'el\'ephants ne sont pas autoris\'es \`a traverser la rivi\`ere,

 \hskip1.5in ils doivent rester dans leur moiti\'e de plateau.
}

{\hskip3.5in - R\`egles des Echecs Chinois
}

\tableofcontents

\section*{Introduction}

Welschinger invariants of real rational symplectic four-folds \cite{W1} provide an invariant count of real rational pseudo-holomorphic curves
in a given homology class. In particular, for real del Pezzo surfaces, Welschinger invariants count real rational algebraic curves
in a given divisor
class that pass through
a
generic conjugation-invariant
collection
of points with a given number of real points among them, and this count depends neither on the
variation
of
the point constraints, nor on the variation of the surface in its deformation class (see
\cite{IKS3}). In \cite{W2}, J.-Y.Welschinger suggested
a ``relative" version that works
 in the presence of
a smooth surface with boundary
selected
inside
the real part of
a real rational symplectic  four-fold.
His
relative
invariant
deals with
real rational
curves tangent once to the boundary of this surface, and
to achieve the invariance with respect to the
variation of
point constraints
and deformations, one has to count not only
curves
tangent to
this boundary, but also cuspidal curves, reducible curves, and curves with imposed
tangency
directions
at the fixed points.

We
introduce
relative
enumerative invariants of different nature that take place
for the real
rational surfaces $\Sig$
containing
a unique smooth real $(-2)$-curve $E\subset \Sig$ and satisfying $K_\Sig C<0$, for
any other irreducible curve $C\subset\Sigma$, and $K_\Sig^2>0$.
In Introduction, we call them for brevity real nodal Pezzo surfaces.
Our
invariants count (with signs) real rational curves
that are tangent to $E$ at {\it all their intersection points}. First, such a count and its invariance were observed
in \cite[Corollary 4.1]{IKS2} in
the
enumeration
of real rational curves disjoint from $E$, in which case the invariant can be expressed via absolute Welschinger invariants.
In the symplectic setting,
a development of this observation to
a count of
real pseudo-holomorphic curves
that intersect a number of disjoint
$(-2)$-curves transversaly (at imaginary points)
was achieved by E. Brugall\'e \cite[Theorem 3.9]{Br}, also by showing that his invariant is
a combination of absolute Welschinger invariants.

The development we suggest in this paper is two-fold:
we extend the above cited result from \cite{IKS2} to counting curves  through real collections
of points containing any amount of pairs of complex conjugate points and, overall, to divisor classes that have nonzero intersection with $E$.
In the case of real nodal del Pezzo surfaces
of degree $\ge2$, the invariance
of the count that we introduce
takes place
unconditionally
(Theorem \ref{t2}, Section
\ref{invariance_statements}).
For real nodal del Pezzo surfaces of degree $1$, the situation is unequal: on one hand, as we show in Section \ref{non-invariance}
for some divisor classes the invariance fails; on the other hand, as we show in Section \ref{invariance_statements},
it takes place under certain restrictions on the divisor class (see Theorems \ref{t4p} and \ref{cor1}). All our relative invariants
remain constant in generic (global) variations of nodal del Pezzo surfaces (Theorem \ref{t1}).

Similar invariants can be defined relative to other types of divisors than $(-2)$-curves. Namely, if we start
from a nodal del Pezzo surface, but restrict our attention to divisors that avoid a pair of disjoint complex conjugate $(-1)$-curves crossing $E$ with multiplicity 1 and then blow down these $(-1)$-curves, we come to invariants on the constructed del Pezzo surface relative to the image of $E$.
For example, iterating this procedure we obtain pure relative enumerative invariants for the projective plane relative to a conic
and for projective spherical quadrics relative to hyperplane section (Section \ref{sec-se}).

\begin{figure}
\setlength{\unitlength}{1cm}

\begin{picture}(14.5,6)(0,-0.7)
\epsfxsize 145mm \epsfbox{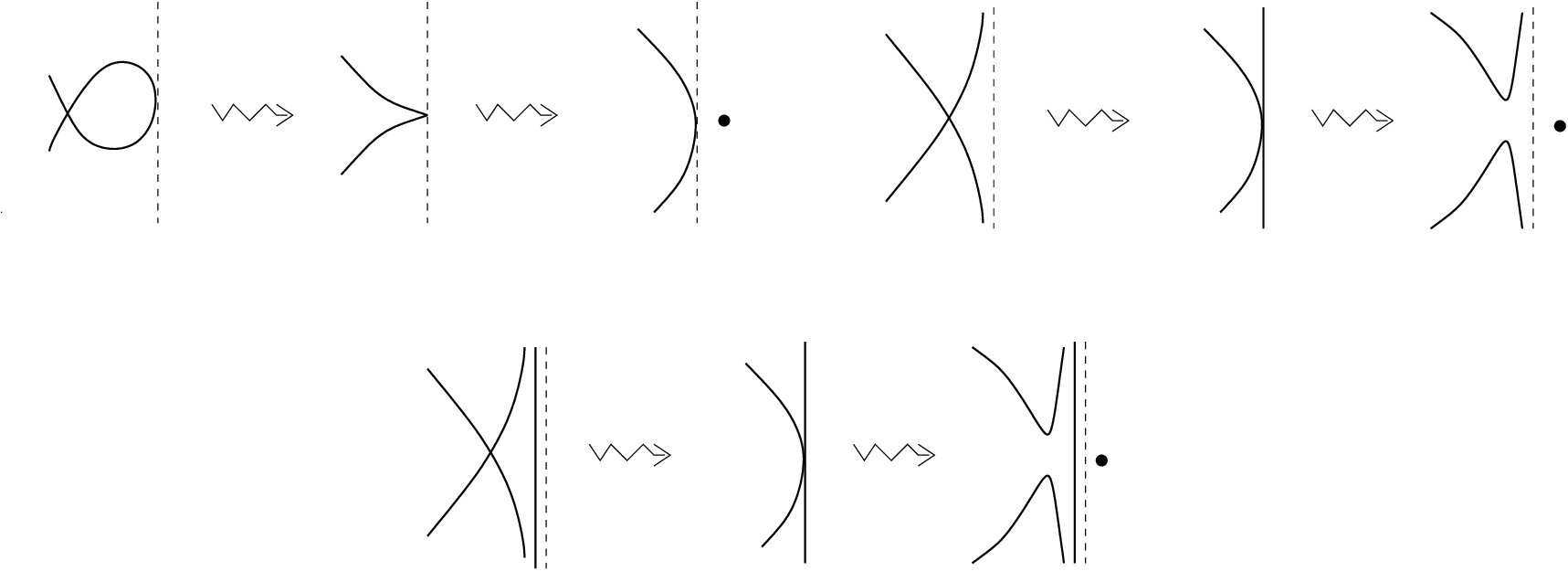}
\put(-10.8,2.5){(a)}\put(-3,2.5){(b)}\put(-7.3,-0.6){(c)}
\put(-12.9,5.1){$E$}\put(-10.4,5.1){$E$}\put(-7.9,5.1){$E$}
\put(-5.2,5.1){$E$}\put(-2.7,5.1){$E$}\put(-0.2,5.1){$E$}
\put(-9.4,1.9){$E_0$}\put(-6.9,1.9){$2E_0$}\put(-4.4,1.9){$E_0$}
\end{picture}
\caption{Bifurcations of real rational curves} \label{fig-bif}
\end{figure}
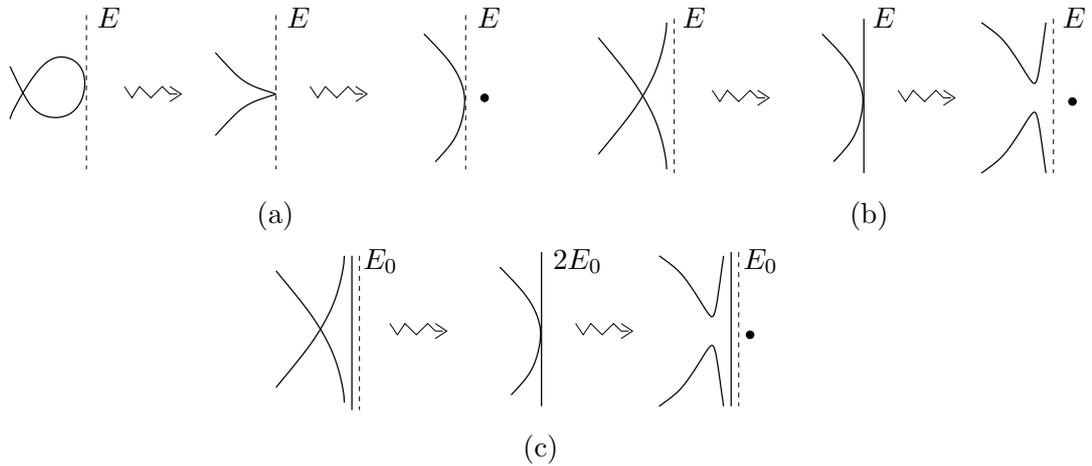

The proofs are based on the following ideas.
Similarly
to  \cite[Section 4]{IKS2}, and \cite{Br,BP}
in the symplectic setting, we use deformations of real nodal del Pezzo surfaces into genuine
real del Pezzo surfaces in order to directly relate our count
in the case of divisors disjoint from the $(-2)$-curve $E$
to absolute Welschinger invariants.
In other cases we classify all degenerations of the set of counted complex rational
curves, occurring in generic one-dimensional families of point constraints (Section \ref{sec2.3}) and then analyze
bifurcations of real curves in crossing walls associated with above degenerations (Section \ref{s3}).
The
list of bifurcations contains
a cuspidal one, where
the degenerated curve
has a cusp lying outside $E$.
As in the absolute Welschinger theory \cite{W1},
the invariance
in this bifurcation follows from the fact that two real rational curves
having opposite
weight appear or disappear. New bifurcations
include
the bifurcation via a curve having a cusp on $E$, see Figure \ref{fig-bif}(a),
and the bifurcation via a curve that splits off $E$, see Figure  \ref{fig-bif}(b).
In these bifurcations, a curve having a non-solitary real node turns into a
curve with a solitary real node (or {\it vice versa}), both curves having the same
weight
(see definitions in Section \ref{relation_Welschinger}).

Each real nodal del Pezzo surface of degree $1$ contains a unique $(-1)$-curve $E_0$ such that $EE_0=2$ (see Lemma \ref{l1}(2iv)).
For these surfaces one encounters
additional phenomena that correspond to splitting off the curve $E_0$ (see an example in Figure \ref{fig-bif}(c))
with loss of invariance.
These
phenomena impose the limits mentioned in Theorems \ref{t4p}, \ref{cor1}, and \ref{t1}.

The paper is organized as follows. In Section \ref{sec-mr}, we introduce the class of real rational surfaces
under consideration,
the sets of real rational curves on these surfaces and the counting rules,
and formulate the invariance statements. Section \ref{s1} is devoted to
the study of families of complex rational curves on nodal del Pezzo surfaces.
In
Sections \ref{s2}\;-\ref{sec-t1}, we prove the invariance
statements.
Section \ref{sec-cr}
contains simplest cases illustrating the new
invariants
and
examples
of the lack of invariance.

\section{Main results}\label{sec-mr}

Throughout the paper, we indicate the field of definition of a variety only if the variety is defined
over $\R$, otherwise it is assumed to be defined over $\C$.

\subsection{Counting rules}\label{sec1.1} A pair $(\Sig, E)$, where $\Sig$ is a smooth
rational algebraic
surface of degree $K_\Sig^2>0$ and $E\subset\Sig$
is a $(-2)$-curve,
is called a {\it nodal del Pezzo pair}\footnote{Nodality refers here exclusively to $E\subset \Sig$ and does not forbid any other, even not nodal, singularity
of $\Sig$.}
(briefly,
{\it nDP-pair}).
Recall that by a $(-2)$-curve one means
a reduced
irreducible smooth rational curve with self-intersection
$-2$, which is equivalent, by
adjunction formula,
to
an
assumption to be
an irreducible curve $E$ with
$E^2=-2$ and $EK_\Sig
=
0$.

For an nDP-pair $(\Sig,E)$, we denote by
$\Pic_+(\Sig,E)$ the semigroup in $\Pic(\Sig)$ generated by
irreducible curves $C\ne E$.
A {\it real nDP-pair} is a complex nDP-pair $(\Sig,E)$ equipped with an antiholomorphic involution $\conj:(\Sig,E)\to(\Sig,E)$.
For a real nDP-pair $(\Sig, E)$,
we put
$\R\Sig=\Fix(\conj)$ and $\R E=E\cap\R\Sig$.

Let $(\Sig,E)$ be a real nDP-pair with $\R\Sig\ne\emptyset$.
A connected component $F$ of $\R\Sigma$ is called
{\it admissible} if either $F\setminus\R E$ consists of
two connected components, $F^+$ and $F^-$, or $F\cap\R E=\emptyset$. In the latter case, we set $F^+=F$, $F^-=\emptyset$.

Fix an admissible component $F$ and consider
 a real divisor class  $D\in\Pic_+(\Sig,E)$
 such that \begin{equation}
DE=0\mod2,\quad
-DK_\Sig-1-DE/2\ge0,\quad\text{and}\ D\ \text{is primitive
if}\ D^2\le0 \ .
\label{e1}\end{equation}
Put
$
r=-DK_\Sig-1-DE/2,
$
choose an integer $m$ such that $0\le 2m\le r$, and
introduce a real structure $c_{r,m}$ on $\Sig^r$
that maps $(w_1, \ldots , w_r)\in\Sig^r$ to $(w'_1, \ldots , w'_r)\in\Sig^r$ with
$w'_i = \conj(w_i)$ if $i>2m$, and
$(w'_{2j-1}, w'_{2j}) = (\conj(w_{2j}), \conj(w_{2j-1}))$ if $j\le m$.
An $r$-tuple
$\bw=(w_1,...,w_r)$
is $c_{r,m}$-invariant
if and only if $w_i$ belongs to the real part $\R\Sig$ of $\Sig$ for
$i>2m$ and $w_{2j-1}, w_{2j}$ are conjugate to each other
for $j \le m$.
Denote by ${\mathcal P}_{r,m}(\Sig,F^+)$
the
subset of $\Sig^r$ formed by
the $c_{r,m}$-invariant $r$-tuples $\bw=(w_1,...,w_r)$ with pairwise
distinct $w_i\in\Sig$ such that $w_i \in F^+$ for all $i > 2m$.

Consider the moduli space ${\mathcal M}_{0, r}(\Sig,D)$
parametrizing the isomorphism classes $[\bn:\PP^1\to\Sig, \bp]$ of pairs $(\bn:\PP^1\to\Sig, \bp)$, where
$\bn:\PP^1\to\Sig$ is a regular map
such that
$\bn_*\PP^1= |D|$,
and $\bp$ is a sequence of $r$ pairwise distinct points in $\PP^1$.
Define  ${\mathcal V}_r(\Sig,E,D)\subset{\mathcal M}_{0, r}(\Sig,D)$
to be the
subset consisting of elements $[\bn:\PP^1\to\Sig, \bp]$
subject to the following
restriction: $\bn^*(E)=2\bd$,
where
$\bd\in\Div(\PP^1)$ is an effective divisor.
Given $\bw\in{\mathcal P}_{r,m}(\Sig,F^+)$, denote by
${\mathcal V}^\R_r(\Sig,E,D,\bw)$ the set of real elements $[\bn:\PP^1\to\Sig, \bp]\in{\mathcal V}_r(\Sig,E,D)$
such that
$\bn(\bp) = \bw$.

An nDP-pair $(\Sig, E)$ is called a {\it uninodal DP-pair} if
$-CK_\Sig>0$
for any
reduced
irreducible curve $C\ne E$ on $\Sig$.

Let $(\Sig, E)$ be a uninodal DP-pair,
and let $\bw \in {\mathcal P}_{r,m}(\Sig,F^+)$ be a generic $r$-tuple.
By Lemma \ref{lem-fin},
the set ${\mathcal V}^\R_r(\Sig,E,D,\bw)$ is finite.
Denote by ${\mathcal V}^{im,\R}_r(\Sig,E,D,\bw)\subset{\mathcal V}^\R_r(\Sig,E,D,\bw)$
the set of elements $[\bn:\PP^1\to\Sig, \bp]$ such that $\bn$ is birational
onto its image $C=\bn(\PP^1)$, the divisor $\bn^*(E)=2\bd_0$ consists of $l=DE/2$ distinct double points,
and the curve $C$ is immersed outside
$E$.
Observe that the complex conjugation on the source $\PP^1$
of $\bn$ is isomorphic to the standard one and
its fixed point set $\R P^1$
is mapped onto the one-dimensional real component of $C$:
if $r-2m>0$ it is evident,
if
$2m=r$, then either $DE/2$ or $-DK_\Sig $ is odd, whereas the former
quantity is the degree of $\bd=\bn^*(E)/2$
and the latter quantity has the same parity as
$C^2$.
Denote
by ${\mathcal V}^{im,\R}_r(\Sig,E,F,D,\bw)\subset{\mathcal V}^{im,\R}_r(\Sig,E,D,\bw)$
the set consisting of the elements $[\bn:\PP^1\to\Sig, \bp]$
such that $\bn(\R P^1)\subset F$.

Given a real
curve $C\in{\mathcal V}^{im,\R}_r(\Sig,E,D,\bw)$,
we associate two special multiplicities, $s(C,z)$  and  $ns(C,z)$, with each real
singular point $z\in C\setminus E$.
Namely, we define
$s(C,z)$ to be the sum of the intersection multiplicities $C_k\bar C_k$ of
pairs $C_k, \bar C_k$ of complex conjugate irreducible components of the germ $(C,z)$, and $ns(C,z)$ to be the sum of the intersection multiplicities $C_iC_j$  of
pairs $C_i,C_j$, $i\ne j$,  of real irreducible components of $(C,z)$.

For example, if $z$ is a real nodal point of $C$, then $s(C,z)$ equals $1$ if $z$
is solitary (given by $x^2+y^2=0$ in suitable real local
coordinates), and equals $0$ if $z$ is non-solitary (given by
$x^2-y^2=0$);
respectively
$ns(C,z)$ is equal here to $1-s(C,z)$.
The number $(-1)^{s(C, z)}$ is called the {\it Welschinger sign} of the real node $z$.
If $z\in C\setminus E$ is singular but not a node,
we can perform a real nodal equigeneric deformation of the germ $(C,z)$
({\it i.e.}, with the maximal possible number of nodes equal to
the so-called $\delta$-invariant $\delta(C,z)$). Then,
the
number of solitary
(respectively, non-solitary)
real nodes on the resulting nodal curve
is congruent to $s(C,z)$
(respectively, $ns(C,z)$) modulo $2$
({\it cf.} Lemma \ref{l4} below).

We define the following quantities for each admissible component $F\subset \R\Sigma$,
each non-empty half $F^+\subset F$,
each $\conj$-invariant class $\varphi\in H_2(\Sig\setminus F; \Z/2)$,
and each generic $\bw \in {\mathcal P}_{r,m}(\Sig,F^+)$:
\begin{equation}
RW_{m}(\Sig,E,F^+,\varphi,D,\bw)=\sum_{\xi\in{\mathcal V}^{im,\R}_r(\Sig,E,F,D,\bw)}\mu(F^+, \varphi, \xi)\ ,
\label{e22n}\end{equation}
where
\begin{equation}\mu(F^+, \varphi, \xi)=(-1)^{
\varphi\cdot C_{1/2}}
\prod_{z\in\Sing(C)\cap
F^+}(-1)^{s(C,z)}\cdot\prod_{z\in\Sing(C)\cap
F^-}(-1)^{ns(C,z)}\;,\label{esign}\end{equation}
and, for each $\xi= [\bn:\PP^1\to \Sigma, \bp]$,
the symbols
$C$ and $C_{1/2}$ stand for the images, under $\bn$, of $\PP^1$ and  of one of the connected components of $ \PP^1\setminus\R P^1$, respectively,
and $\varphi\cdot C_{1/2}$ means the intersection number.

The structure of the formula (\ref{esign}) reflects the invariance of the count of
real rational curves when the point constraints vary in generic one-parameter families. The second factor
in the right-hand side depends on the parity of the number of real solitary nodes of $C$
(or of its nodal equigeneric deformation)
as the original Welschinger sign \cite{W1}. The reason
for the solitary nodes
to be counted only in the domain $F^+$ is as follows. In almost all cases considered
in Section \ref{invariance_statements},
the one-dimensional part of $\R C$ lies entirely in $\overline{F^+}$, and hence the solitary nodes in $F^-$ cannot degenerate into cusps belonging to $F^-$ when the considered rational curves vary in generic one-parameter families, thus these solitary nodes do not matter in the invariance problem.
In these cases,
the third factor in formula (\ref{esign}) is always $1$. However,
there are special cases
where the one-dimensional part of $\R C$ jumps from $\overline{F^+}$
to $\overline{F^-}$ (and vice versa). Then, the third factor becomes non-trivial and it serves to balance the invariance in local bifurcations as shown in Figure \ref{fig-bif} when a non-solitary node in $F^-$ crosses $E$
and turns into a solitary node in $F^+$.
The $\varphi$-twisting factor $(-1)^{\varphi\cdot C_{1/2}}$
was introduced
in \cite{IKS1}. A typical
class $\varphi$
is represented by
a
(possibly empty) union of
components of $\R\Sigma$
which are all different from $F$. Then, the exponent
$\varphi\cdot C_{1/2}$ is
the number of solitary nodes of
$C$ in the chosen components
of $\R\Sigma$.

\begin{example}\label{exan}
{\rm Consider the following enumerative problem.
Let $C_2\subset\PP^2$ be a smooth real conic with $\R C_2\simeq S^1$.
Take for $\Sig = \PP^2_{(0,3)}$ the plane $\PP^2$ blown up
at three pairs of complex conjugate points on $C_2$,
and take for $E$ the strict transform of $C_2$.
The real part $\R\Sigma$ consists of one connected component,
and the complement $\R\Sigma \setminus \R E \simeq \R P^2\setminus\R C_2$ consists
of a disc $F^o$ and a M\"obius band $F^{no}$.
Put $D = 4L - E_1 - \ldots - E_6$, where $L \subset \Sigma$ is a pull-back of a generic line in $\PP^2$
and $E_1$, $\ldots$, $E_6$ are the exceptional divisors of the blow ups.
In this situation, $DE = 2$ and $D$ satisfies conditions (\ref{e1}).
Put in addition
$\varphi = 0 \in H_2(\Sigma \setminus \R\Sigma; \Z/2)$.

Given a generic collection $\bw$ of $-DK_\Sigma - 1 - DE/2 = 4$ real points in $F^+ = F^o$,
we are interested in the real rational curves $C \subset \Sigma$ which belong
to the linear system $|D|$, are tangent to $E$ at the only intersection point with $E$,
and pass through all the 4 points of $\bw$.
The signed enumeration of these real rational curves (with the sign described in (\ref{esign}))
gives rise to the number $RW_{0}(\Sig,E,F^o,\varphi,D,\bw)$.
It can be shown that each such curve $C$ is nodal;
furthermore, the one-dimensional part of $\R C$ is entirely contained in $\overline{F^o}$.
So, the sign $(-1)^s$ of $C$ is determined
by the parity of the number $s$ of solitary nodes of $C$ which belong to $F^o$.

Another option is to choose a generic collection $\bw'$ of $4$ real points in $F^{no}$
and consider the number $RW_{0}(\Sig,E,F^{no},\varphi,D,\bw')$.
As we will see, the number $RW_{0}(\Sig,E,F^o,\varphi,D,\bw)$ (respectively, $RW_{0}(\Sig,E,F^{no},\varphi,D,\bw')$)
does not depend on the choice of a generic collection $\bw$ (respectively, $\bw'$)
provided that the points of the collection are in $F^o$ (respectively, $F^{no}$).
However, these numbers $RW_{0}(\Sig,E,F^o,\varphi,D,\bw)$ and $RW_{0}(\Sig,E,F^{no},\varphi,D,\bw')$
are not equal (for the precise values and a further discussion of this example, see Section \ref{sec-se}).}
\end{example}

\subsection{Relation to Welschinger invariants}\label{relation_Welschinger}

Theorems of this section relate, in the special case $DE=0$, the numbers $RW_m(\Sig,E,F^+,\varphi,D, \bw)$
to certain modified Welschinger invariants.
In particular, this proves the invariance of the numbers $RW_m(\Sig,E,F^+,\varphi,D, \bw)$ in the case $DE = 0$.

An
instance of such a special case is provided by a slight modification of
Example \ref{exan}:
the surface $\Sig$ and the curve $E$ remain the same, but we put
$D = -K_\Sigma = 3L - E_1 - \ldots - E_6$.

A {\it perturbation} of
a uninodal DP-pair $(\Sig,E)$ is a proper submersion $f$
of a smooth variety $\mathfrak X$  to
$\Delta_a=\{z\in\C : \vert z\vert < a\}$,
$a>0$, with
$f^{-1}(0)=\Sig$ and such that $f^{-1}(z)=\Sig_z$ is a del Pezzo surface for each $z \ne 0$.
A perturbation $f : \mathfrak X \to \Delta_a$ is called {\it real} if
$\mathfrak X$ is equipped with a real structure $c : \mathfrak X\to \mathfrak X$
such that $f \circ c = \conj\;\circ f$.

Let $f: \mathfrak X \to \Delta_a$ be
a real perturbation
of a real uninodal DP-pair $(\Sig,E)$.
The real part $\Pic^\R({\mathfrak X}_t)$ of $\Pic({\mathfrak X}_t)$, $t\in (-a,a)$, is naturally
identified with $\Pic^\R(\Sig)$. Given a divisor class $D\in\Pic(\Sig)$, a connected component $F$ of $\R\Sig$, and a $\Conj$-invariant class $\varphi\in H_2(\Sig\setminus
F; \Z/2)$, we obtain a continuous family of tuples $(D, F_t,\varphi_t)$, $t\in (-a,a)$.

Thus, in particular, the modified Welschinger invariants $W_m({\mathfrak
X}_t,D-2sE,F_t,\varphi_t)$ are well defined for each $t \in (-a, a)$,
$t\ne 0$, and each $s\ge 0$. They are given by taking the sum
$\sum
(-1)^{C_+\cdot C_- +
\varphi_t \cdot C_+}$
over all immersions $\nu: \PP^1 \to
{\mathfrak X}_t$ representing
the given divisor class
$D'=D-2sE$ on ${\mathfrak X}_t$ and interpolating a given generic collection
$\bw \in {\mathcal P}_{r,m}(
{\mathfrak X}_t,F_t)$,
where $r = -K_{{\mathfrak X}_t}D' - 1$
and $C_{\pm}=\nu(\PP^1_\pm)$ with $\PP^1_+,\PP^1_-$ being the two connected components
of $\PP^1\setminus\R\PP^1$ (see \cite[Section 1]{IKS3}).

\begin{theorem}\label{t3n}
Let $(\Sig,E)$ be a real uninodal DP-pair,
$F\subset\R\Sig$ an admissible connected component, $\varphi\in H_2(\Sig\setminus F; \Z/2)$
a $\Conj$-invariant class,
and $D\in\Pic_+(\Sig,E)$ a real
divisor class matching conditions (\ref{e1}) and such that $DE=0$.
Assume that $[\R E]\ne0\in H_1(\R\Sig; \Z/2)$ or $\R E=\emptyset$
{\rm (}in particular, $F = F^+${\rm )}.
Then, for any $0\le m\le r/2$,
where $r=-DK_\Sig-1$, any generic configuration $\bw\in {\mathcal P}_{r,m}(\Sig,F^+)$, and
any real perturbation ${\mathfrak X} \to
\Delta_a$ of $(\Sig,E)$, we have
\begin{eqnarray}RW_m(\Sig,E,F^+,\varphi,D,\bw)&=&W_m({\mathfrak
X}_t,D,F_t,\varphi_t)\nonumber\\ & &+2\sum_{s\ge1}(-1)^s\; W_m({\mathfrak
X}_t,D-2sE,F_t,\varphi_t)\label{e2504c}\end{eqnarray}
for all $t \in (-a, a)$,
$t\ne 0$.
In particular, $RW_m(\Sig,E,F^+,\varphi,D,\bw)$ does not depend on the choice of a generic configuration
$\bw\in {\mathcal P}_{r,m}(\Sig,F^+)$.
\end{theorem}

Suppose that $\R E\ne\emptyset$ and $[\R E]=0\in H_1(\R\Sig; \Z/2)$. Consider a real proper map
$f$ of a smooth real variety $\mathfrak X$  to
$\Delta_a$
such that:
\begin{itemize}
\item $f$ is a submersion at all but one point; the latter point
is a simple

critical point and belongs to $f^{-1}(0)$;
\item $f^{-1}(z)={\mathfrak X}_z$ is a del Pezzo surface for each $z \ne 0$;
\item  $f^{-1}(0)={\mathfrak X}_0$, where  ${\mathfrak X}_0$
is obtained by a regular map $\varsigma : \Sig \to {\mathfrak X}_0$
that contracts $E\subset \Sig$ to a (non-solitary real nodal) point.
\end{itemize}
Denote by $G$ the connected component of $\R \Sig$ such that $G\cap \R E\ne\emptyset$.
The map $f$ is called a {\it dividing
surgery}
of $(\Sig, E)$, if for each
$t \in (0, a)$
the real part $\R {\mathfrak X}_t$ of ${\mathfrak X}_t$ has two connected components, $G_t^+$ and $G_t^-$,  that converge, as $t\to 0$, to the closures of
two connected components of
$\varsigma(G)\setminus
\varsigma(E)$.
If $G$ coincides with $F$, we assume that $G^+_t$ converges to $F^+$.
For each connected component $H\ne G$ of $\R \Sig$, the image $\varsigma(H)\subset {\mathfrak X}_0$ is non-singular and is included in a topologically trivial family of connected components $H_t$ of
$\R {\mathfrak X}_t$, $t\in (-a,a)$.

For any dividing
surgery
$f :  \mathfrak X \to \Delta_a$,
the real part $\Pic^\R({\mathfrak X}_t)$ of $\Pic({\mathfrak X}_t)$,
$t \in (0, a),$
is
identified with $\{D\in\Pic^\R(\Sig)\ :\ DE=0\}$, see \cite[Proposition 4.2]{IKS2}.
In
addition,
any $\conj$-invariant class $\varphi\in H_2(\Sig\setminus G; \Z/2)$
(respectively, $\varphi\in H_2(\Sig\setminus (G\cup H); \Z/2)$, $H\ne G$)
continuously deforms into
classes
$\varphi_t\in H_2({\mathfrak X}_t\setminus(G_t^+\cup G_t^-); \Z/2)$
(respectively, $\varphi_t\in H_2({\mathfrak X}_t\setminus(G_t^+\cup G_t^-\cup H); \Z/2)$),
$t \in (0, a)$,
that are invariant under the real structure.

\begin{theorem}\label{t3}
Let $(\Sig,E)$ be a real uninodal DP-pair,
$F\subset\R\Sig$ an admissible connected component,
and $D\in\Pic_+(\Sig,E)$ a real
divisor class matching conditions (\ref{e1}) and such that $DE=0$.
Pick any $0\le m\le r/2$, where $r=-DK_\Sig-1$, and any generic configuration
$\bw\in {\mathcal P}_{r,m}(\Sig,F^+)$.
Assume that $\R E\ne\emptyset$ and $[\R E]=0\in H_1(\R\Sig; \Z/2)$.
Let ${\mathfrak X} \to
\Delta_a$
be a dividing
surgery of $(\Sig,E)$.

(1) If $\R E\subset F$ and $\varphi\in H_2(\Sig\setminus F; \Z/2)$ is
a $\Conj$-invariant class, then
we have
$$RW_m(\Sig,E,F^+,\varphi , D,\bw)=W_m({\mathfrak X}_t,D,F_t^+,\varphi_t+[F_t^-])
$$
for any
$t \in (0, a)$.

(2) If $\R E$ lies in a component $G\ne F$ of $\R\Sig$ {\rm (}in particular, $F = F^+${\rm )}
and $\varphi\in H_2(\Sig\setminus (F\cup G); \Z/2)$ is
a $\Conj$-invariant class, then
we have
$$RW_m(\Sig,E,F^+,\varphi , D,\bw)=W_m({\mathfrak X}_t,D,F_t,\varphi_t+[G_t^-])
=W_m({\mathfrak X}_t,D,F_t,\varphi_t+[G_t^+])$$
for any
$t \in (0, a)$.

(3) The numbers $RW_m(\Sig,E,F^+,\varphi , D,\bw)$ as defined in paragraphs (1) and (2) do not depend on the choice of a generic
configuration $\bw\in {\mathcal P}_{r,m}(\Sig,F^+)$.
\end{theorem}

\subsection{Invariance statements}\label{invariance_statements}
\begin{theorem}\label{t2}
Let $(\Sig,E)$ be a real uninodal DP-pair  with $\deg\Sig=K_\Sig^2\ge2$. If $F\subset \R\Sigma$ is an admissible component,
$\varphi\in H_2(\Sig\setminus F; \Z/2)$ is a $\conj$-invariant class, and
$D\in\Pic_+(\Sig,E)$ is a real
divisor class matching conditions (\ref{e1}) and satisfying $r=-DK_\Sig-DE/2-1>0$,
then, for any integer $0\le m \le r/2$,
the number
$RW_m(\Sig,E,F^+,\varphi,D,\bw)$ does not depend on the choice of
a generic
$\bw\in{\mathcal P}_{r,m}(\Sig,F^+)$.
\end{theorem}

Theorem \ref{t2} implies,
in particular, that
in Example \ref{exan}
the difference between the number of real rational curves under consideration
in $\Sig=\PP^2_{(0,3)}$
which pass through $\bw$
and
have an even (respectively, odd) number of solitary nodes in $F^o$
does not depend on the choice
of a generic collection $\bw \subset F^o$.
Similarly, the difference between the number of real rational curves under consideration
passing through $\bw'$ and
having an even (respectively, odd) number of solitary nodes in $F^{no}$
does not depend on the choice
of a generic collection $\bw' \subset F^{no}$.
The resulting numbers $RW_{0}(\Sig,E,F^o,\varphi,D,\bw)$ and  $RW_{0}(\Sig,E,F^{no},\varphi,D,\bw')$
are equal to $48$ and $80$, respectively
(see Section \ref{sec-se}).
In particular, this implies that through any generic collection $\bw$
(respectively, $\bw'$) of $4$ points in $F^o$ (respectively, $F^{no}$)
one can always trace at least $48$ (respectively, $80$) real rational curves
tangent to $E$ and belonging to the linear system $|4L - E_1 - \ldots - E_6|$.

\smallskip

In order to extend the statement of Theorem \ref{t2} to uninodal DP-pairs of degree $1$, we have to introduce
additional restrictions, and these restrictions are essential as we explain in Section \ref{sec-cr}.

For any uninodal DP-pair $(\Sig,E)$ of degree $1$,
there exists a unique $(-1)$-curve
in $|-(K_\Sig + E)|$
(see Lemma \ref{l1} in Section \ref{sec2}).
We denote this curve by $E_0$.

A uninodal DP-pair $(\Sig, E)$
is said to be {\it tangential}
if it is of degree $1$ and the curve $E_0$ is
tangent to $E$.
The tangential DP-pairs are not generic among uninodal DP-pairs of degree $1$ (see Proposition \ref{prop1}).

\begin{theorem}\label{t4p}
Let
$(\Sig,E)$ be a real uninodal DP-pair of degree $1$.
Let $F\subset\R\Sigma$ be an admissible component,
$\varphi\in H_2(\Sig\setminus F; \Z/2)$ a $\conj$-invariant class, and
$D\in\Pic_+(\Sig,E)$ a real
divisor class matching conditions (\ref{e1}) and satisfying $r=-DK_\Sig-DE/2-1>0$.
If
either $DE = 0$, or $DE = 2$ and the uninodal DP-pair $(\Sig, E)$ is not
tangential, or $\R E \cap F = \emptyset$ and $\R E_0\cap F = \emptyset$,
then,
for any integer $0\le m \le r/2$,
the number
$RW_m(\Sig,E,F^+,\varphi,D,\bw)$ does not depend on the choice of
a generic
$\bw\in{\mathcal P}_{r,m}(\Sig,F^+)$.
\end{theorem}

For
any $\bw\in{\mathcal P}_{r,m}(\Sig,F^+)$ and any $\xi = [\bn:\PP^1\to\Sig, \bp] \in {\mathcal V}^{im,\R}_r(\Sig,E,F,D,\bw)$,
the one-dimensional real component
$\bn(\R\PP^1)$ of
$\bn(\PP^1)$ does not traverse $\R E$, and thus
$\bn(\R\PP^1)$ is entirely contained either in the
closure $\overline {F^{\;+}}$ of $F^+$, or in the closure $\overline {F^{\;-}}$ of $F^-$. If $r>2m$, then
$\bn(\R\PP^1)\subset\overline {F^{\;+}}$, since the real point constraints lie in $F^+$. In the case $r=2m$,
under certain additional conditions, the statements of Theorems \ref{t2} and \ref{t4p} can be refined
in order to
obtain invariants that count separately curves
with
$\bn(\R\PP^1)\subset\overline {F^{\;+}}$ and with
$\bn(\R\PP^1)\subset\overline {F^{\;-}}$.

\begin{theorem}\label{cor1}
Let
$(\Sig,E)$ be a real uninodal DP-pair, $F\subset\R\Sigma$ an admissible component,
$\varphi\in H_2(\Sig\setminus F; \Z/2)$ a $\conj$-invariant class, and
$D\in\Pic_+(\Sig,E)$ a real
divisor class matching conditions (\ref{e1}).
In the case $\deg\Sig=1$, suppose that
either $DE = 0$, or $DE = 2$ and the uninodal DP-pair $(\Sig, E)$ is not
tangential.
If $r=2m$, $F\cap\R E\ne\emptyset$, and $D$ is not
representable in the form $D=D'+\conj_*D'$ with $|D'|$ containing an irreducible rational curve
and satisfying $D'E\equiv1\mod2$,  then the numbers
$$RW^+_{m}(\Sig,E,F^+,\varphi,D,\bw)=\sum_{\renewcommand{\arraystretch}{0.6}
\begin{array}{c}
\scriptstyle{\xi=[\bn]\in{\mathcal V}^{im,\R}_r(\Sig,E,F,D,\bw)}\; :\\
\scriptstyle{\bn(\R P^1)\subset\overline {F^{\;+}}}
\end{array}}
\mu(F^+, \varphi, \xi),$$
$$RW^-_{m}(\Sig,E,F^+,\varphi,D,\bw)=\sum_{\renewcommand{\arraystretch}{0.6}
\begin{array}{c}
\scriptstyle{\xi=[\bn]\in{\mathcal V}^{im,\R}_r(\Sig,E,F,D,\bw)}\; :\\
\scriptstyle{\bn(\R P^1)\subset\overline {F^{\;-}}}
\end{array}}
\mu(F^+, \varphi, \xi)
$$
do not depend on the choice of
a generic
$\bw\in{\mathcal P}_{r,m}(\Sig,F^+)$, and
$$RW_{m}(\Sig,E,F^+,\varphi, D)=RW^+_{m}(\Sig,E,F^+,\varphi, D)+RW^-_{m}(\Sig,E,F^+,\varphi, D).$$
\end{theorem}

\smallskip

By
an {\it elementary deformation} of an nDP-pair $(\Sig,E)$ we understand a
proper submersion $f$
of a pair of smooth varieties $(\mathfrak X, \mathfrak E)$  to
$\Delta_a$
such that $f^{-1}(z)=(\Sig_z, E_z)$ is an nDP-pair for each $z\in \Delta_a$ and
$f^{-1}(0)=(\Sig, E)$.
If $f: {\mathfrak X} \to \Delta$ is such an elementary deformation
and $z_1, z_2 \in \Delta_a$, we also say that $f$ is an elementary deformation
between $f^{-1}(z_1)$ and $f^{-1}(z_2)$.
An elementary deformation $f : (\mathfrak X, \mathfrak E)\to \Delta_a$ is called {\it real} if the pair
$(\mathfrak X, \mathfrak E)$ is equipped with a real structure $c : (\mathfrak X, \mathfrak E)\to (\mathfrak X, \mathfrak E)$
such that $f \circ c = \conj\;\circ f$.

Given a divisor class $D\in\Pic_+(\Sig, E)$, an admissible connected component $F$ of $\R\Sig$, a connected component $F^+$
of $F\setminus\R E$,
and a $\Conj$-invariant class $\varphi\in H_2(\Sig\setminus F;\Z/2)$,
for any real elementary deformation  $f : (\mathfrak X, \mathfrak E)\to \Delta_a$
we obtain a continuous family of tuples $T_t=(\Sig_t, E_t, F_t, F^+_t, \varphi_t, D_t), t\in (-a,a)$,
where
\begin{itemize}
\item $(\Sig_t, E_t)=f^{-1}(t)$ and $D_t\in\Pic_+(\Sig_t, E_t)$,
\item
$\varphi_t\in H_2(\Sig_t \setminus F_t; \Z/2)$ is a class invariant under the restriction of $c$ to $\Sig_t$,
\item
$F_t$ is a connected component of $\R \Sig_t$ and
$F^+_t$ is a connected component
of $F_t\setminus \R E_t$.
\end{itemize}
The tuples $T_t=(\Sig_t, E_t, F_t, F^+_t, \varphi_t, D_t), t\in (-a,a)$, are said to be {\it elementary deformation equivalent}.
Tuples $ T=(\Sig, E, F, F^+, \varphi, D)$ and
$\widetilde T=(\widetilde\Sig, \widetilde E, \widetilde F, \widetilde F^+, \widetilde \varphi, \widetilde D)$
whose underlying nDP-pairs $(\Sig, E)$ and $(\widetilde\Sig, \widetilde E)$ are uninodal
are called {\it deformation equivalent} if they can be connected by a chain
$T=T^{(0)}$, \dots, $T^{(k)}=\widetilde T$ so that any two neighboring tuples in the chain are isomorphic
to elementary deformation equivalent tuples with uninodal underlying nDP-pairs.

\begin{theorem}\label{t1}
Let
$(\Sig,E)$ be a real uninodal DP-pair, $F\subset\R\Sigma$ an admissible component,
$\varphi\in H_2(\Sig\setminus F; \Z/2)$ a $\conj$-invariant class, and
$D\in\Pic_+(\Sig,E)$ a real
divisor class matching conditions (\ref{e1}).
In the case $\deg\Sig=1$, suppose that
either $DE = 0$, or $DE = 2$ and the uninodal DP-pair $(\Sig, E)$ is not tangential,
or $\R E \cap F = \emptyset$ and $\R E_0\cap F = \emptyset$.
If a tuple $\widetilde T=(\widetilde\Sig, \widetilde E, \widetilde F, \widetilde F^+, \widetilde \varphi, \widetilde D)$,
where $(\widetilde\Sig, \widetilde E)$ is a real uninodal DP-pair
which is supposed to be not tangential in the case $\widetilde D \widetilde E = 2$, is deformation equivalent
to  $T=(\Sig, E, F, F^+, \varphi, D)$, then
$$RW_m(\Sig,E,F^+,\varphi,D)=RW_m(\widetilde\Sig, \widetilde E, \widetilde F^+, \widetilde \varphi, \widetilde D)
$$
for any $0\le m \le r/2$, where $r=-DK_\Sig-1-DE/2$.
\end{theorem}

\section{Families of rational curves on nodal del Pezzo
surfaces}\label{s1}

\subsection{Auxiliary statements}
For
curve germs $(C_1,z)$, $(C_2,z)$ on a smooth algebraic surface,
denote by $(C_1\cdot C_2)_z$ the
intersection multiplicity at $z$, and by $\ord(C_1,z)$ the order of
$C_1$ at $z$ ({\it i.e.,} the intersection multiplicity with a
generic smooth curve through $z$).

\begin{lemma}\label{gus}
Let $(B,b_0)$ be a germ of a reduced
analytic space of dimension $d\ge1$, and let
$\{{\mathcal C}_b,\ b\in(B,b_0)\}$ be a flat equisingular family of
reduced irreducible curves
of geometric genus $g$
on a smooth algebraic surface $S$.
Let $Z,W\subset {\mathcal C}_{b_0}$ be disjoint finite sets {\rm (}possibly empty{\rm )}.
For each point $z\in Z$, fix a finite collection
of pairwise transversal smooth curve germs $(L_{z,i}, z)\subset (S,z)$,
$1\le i\le m_z$, $m_z\ge 1$, and for each point $w\in W$, fix
a smooth curve germ
$(M_w, w)\subset (S,w)$.
Suppose that:
\begin{itemize}
\item
${\mathcal C}_b \cap {\mathcal C}_{b_0}$ is finite if
$b\ne b_0$,
\item$({\mathcal C}_b\cdot L_{i,z})_z=({\mathcal C}_{b_0}
\cdot L_{i,z})_{z}$, for all $b\in B$, $z\in Z$, $1\le i\le m_z$,
\item there are sections
$\sigma_w:B\to{\mathcal C}$, $w\in W$, such that
$\sigma_w(b_0)=w$, $\sigma_w(b)\in M_w$, and
$({\mathcal C}_b\cdot M_w)_{\sigma_w(b)}=({\mathcal C}_{b_0}\cdot M_w)_w$
for all $b\in B$, $w\in W$.
\end{itemize}
Then,
\begin{eqnarray}-{\mathcal C}_{b_0}K_S&\ge&2-2g+\sum_{z\in Z}
\left(\ord({\mathcal C}_{b_0},z)+
\sum_{i=1}^{m_z}
\left(({\mathcal C}_{b_0}\cdot
L_{i,z})_z-\ord({\mathcal C}_{b_0},z)\right)\right)\nonumber\\
& &+\sum_{w\in W}(({\mathcal C}_{b_0}\cdot M_w)_w-\ord({\mathcal
C}_{b_0},w))+\sum_{k=1}^s(\ord P_k-1)+(d-1),
\label{egus}\end{eqnarray}
where $P_1, \dots, P_s$ are all the singular local branches of ${\mathcal C}_{b_0}$.
\end{lemma}

{\bf Proof.} First, we reduce the consideration to a one-dimensional
family by fixing a generic set
$Q\subset{\mathcal C}_{b_0}$ of $d-1$ points.
Then, we apply lower bounds for local intersection
multiplicities obtained in \cite[Theorem 2]{GuS}.
Namely, for $b\in B\setminus\{b_0\}$ close enough to $b_0$, we have
\begin{eqnarray}({\mathcal C}_{b_0}\cdot{\mathcal C}_b)_{U(z)}&
\ge&2\delta({\mathcal C}_{b_0},z)
+2\,\ord({\mathcal C}_{b_0},z)-\br({\mathcal C}_{b_0},z)\nonumber\\ & &+
\sum_{i=1}^{m_z}(({\mathcal C}_{b_0}\cdot
L_{i,z})_{z}-\ord({\mathcal C}_{b_0},z)),\quad z\in Z,\nonumber\\
({\mathcal C}_{b_0}\cdot{\mathcal C}_b)_{U(w)}&\ge&2\delta
({\mathcal C}_{b_0},w)
+({\mathcal C}_{b_0}\cdot M_w)_{w}-\br({\mathcal C}_{b_0},w),
\quad w\in W,\nonumber\\
({\mathcal C}_{b_0}\cdot{\mathcal C}_b)_{U(p)}&\ge&2\delta({\mathcal C}_{b_0},p)
+\ord({\mathcal C}_{b_0},p)-\br({\mathcal C}_{b_0},p),\quad p\in
\Sing({\mathcal C}_{b_0})
\setminus(Z\cup W),\nonumber\\
({\mathcal C}_{b_0}\cdot{\mathcal C}_b)_Q&\ge&d-1\ ,\nonumber
\end{eqnarray} where
$({\mathcal C}_{b_0}\cdot{\mathcal C}_b)_{Y}$ denotes the sum of the intersection multiplicities taken at the points of a set $Y$,
the symbol
$\br$ stands for the number of irreducible components of
a given curve germ, and $U$
for a small neighborhood of a given point. Taking into account
that the sum of the right-hand sides of the above inequalities
does not exceed $({\mathcal C}_{b_0})^2$ and combining this with the genus formula
$$({\mathcal C}_{b_0})^2+{\mathcal C}_{b_0}K_S+2=2\sum_{q\in
\Sing({\mathcal C}_{b_0})}\delta({\mathcal C}_{b_0},q)+2g\ ,$$ we derive (\ref{egus}).
\proofend

For each reduced curve germ $(C,z)$ on a smooth surface $\Sig$,
there exists a smooth miniversal
embedded deformation (see
\cite{KasSchlessinger}
for an explicit construction). Denote the base of such a deformation by
$B(C,z)$, the local curves corresponding to elements $t\in B(C,z)$ by $C_t$,
and
by
$B^{eg}(C,z)\subset B(C,z)$ the base of equigeneric deformations of
$(C,z)$ ({\it i.e.}, $\delta$-constant, or equinormalizable).

\begin{lemma}\label{l4}
$B^{eg}(C,z)$ is an irreducible analytic germ
and the local curves $C_t$ are nodal
for generic elements $t\in B^{eg}(C,z)$.
If the
irreducible
components of $(C,z)$ are smooth, then $B^{eg}(C,z)$ is smooth.
If in addition
$\Sig$
and $(C,z)$ are real,
 then,
for any real nodal element $C_t$ of an equigeneric deformation of $(C,z)$, one has:
\begin{itemize}\item the parity of the number $s(C_t)$ of solitary nodes
of $C_t$ coincides with the
parity of $s(C,z)$;
\item the parity of
the number $ns(C_t)$ of non-solitary real nodes of $C_t$ coincides
with the parity of $ns(C,z)$.
\end{itemize}
\end{lemma}

{\bf Proof.}
The first two
statements are
proven in \cite[Proposition
4.17]{DH}. The
third statement is immediate as soon as the numbers $s(C_t)$, $s(C,z)$,  $ns(C_t)$,
and $ns(C,z)$ are
interpreted as linking numbers.
\proofend

\begin{remark}\label{r637}
For an arbitrary isolated curve singularity $(C,z)$ with irreducible components
$P_1,...,P_s$, $s\ge1$, the
stratum $B^{eg}(C,z)\subset B(C,z)$ is not necessarily smooth, but possesses a tangent cone
\begin{equation}
TB^{eg}(C,z)=\{g\in B(C,z)\ :\ \ord\;g\big|_{P_i}\ge2\delta(P_i)+
\sum_{j\ne i} P_i\cdot P_j,\ i=1,...,s\}\ .
\label{e637}\end{equation}
This follows from \cite[Proposition 4.19]{DH} and the fact that the right-hand side is just
the conductor ideal $J^{cond}(C,z)\subset\C\{x,y\}$ (annulator of the module $(\nu_*{\mathcal O}_{C^\vee}/{\mathcal O}_C)_z$,
where $\nu:C^\vee\to C$ is the normalization).
\end{remark}

Assume that
$z$ belongs to a nonsingular curve $E\subset \Sig$. Consider isolated
boundary singularities with boundary $(E,z)$, that is reduced
holomorphic germs $(C,z)$ not containing $(E,z)$. Choose local coordinates
$x,y$ in a neighborhood of $z$ so that
$z=(0,0)$
and $E=\{y=0\}$.
Then the germ at zero $B(C,z,m)$ of the space $\C\{x,y\}/\fm_z^m$, $\fm_z=\langle x,y\rangle^m$ for $m\gg1$
represents both a versal deformation base of the singularity $(C,z)$ and
a versal deformation base
of the boundary singularity $(C,z)$ relative to $E$ (see, for example, \cite{SS}).
Introduce the
substratum $B_E^{eg}(C,z,m)\subset B(C,z,m)$ formed by the elements
$t\in B^{eg}(C,z,m)\subset B(C,z,m)$
with the property that each irreducible component of  the corresponding
local curve $C_t$
meets $E$ at only one point.

\begin{lemma}\label{l5}
(1) Assume that
each irreducible component of $(C,z)$ intersects $E$
with multiplicity $2$. Then, $B_E^{eg}(C,z,m)$ is
smooth and the local curves $C_t$ are nodal for generic elements
$t\in B_E^{eg}(C,z)$.
The tangent space to $B_E^{eg}(C,z,m)$ at $(C,z)$
is
\begin{equation}TB_E^{eg}(C,z,m)\simeq\{g\in B(C,z,m)\ :\ \ord\; g\big|_C\ge2\delta(C,z)+1\}
\ .\label{e2234}\end{equation}

Assume in addition that
$\Sig$ and $E$ are real, $\R E\ne\emptyset$, $z\in\R E$, and $(C,z)$ is a real curve germ.
Let $U$ be a regular neighborhood of $z$ in $\R\Sig$, $U^+,U^-$ connected
components of $U\setminus\R E$. Then, for a nodal
element
$C_t$, $t\in B^{eg}_E(C,z,m)$, the parity of the sum of the number of solitary nodes of $C_t$ in
$U^+$ and the number of non-solitary nodes of $C_t$ in $U^-$ does not depend on the choice
of $t\in B^{eg}_E(C,z,m)$.

(2) Let $(C,z)\subset(\Sig,E)$ be a boundary singularity with respect to the pair $(\Sig,E)$ such that
$C$ is smooth and $(C\cdot E)_z=4$. Then the closure in $B_E(C,z,m)$ of the stratum, parameterising
local curves with two simple tangency points with $E$, is smooth, has codimension $2$ and its tangent space is
\begin{equation}\{g\in B_E(C,z,m)\ :\ \ord\;g\big|_C\ge2\}\;.\label{e2235}\end{equation}
\end{lemma}

{\bf Proof.} (1) Since, in equigeneric deformations, the components
of $(C,z)$ deform separately and independently, to prove the smoothness of
$B^{eg}_E(C,z,m)$ and formula (\ref{e2234}) we have to consider
only the case of an irreducible germ $(C,z)$. For a smooth $(C,z)$
simply tangent to $E$ the smoothness of
$B^{eg}_E(C,z)$ is evident. In the remaining cases,
$(C,z)$ is of type
$A_{2s}$, $s\ge1$, and intersects
$E$ with multiplicity $2$, so that
we have
$z=(0,0)$, $E=\{y=0\}$, $(C,z)=\{y^{2s+1}+x^2=0\}$ in suitable local
coordinates.
A simple computation yileds that the elements of $B_E^{eg}(C,z,m)$ are given by
\begin{equation}
\left((x+\sum_{i\ge0}\alpha_ix^i)^2+y\left(y^s+\sum_{i,j\ge0}
\beta_{ij}x^iy^j\right)^2\right)\left(1+\sum_{i,j\ge0}\gamma_{ij}x^iy^j\right)
\label{egus22}\end{equation}
modulo $\langle x,y\rangle^m$, proving the smoothness of $B_E^{eg}(C,z,m)$. The terms in formula (\ref{egus22})
linear in the parameters, generate the tangent space to $B_E^{eg}(C,z,m)$
$$T_{(C,z)}B_E^{eg}(C,z,m)=\langle x,y^{s+1}\rangle/\langle x,y\rangle^m\ ,$$
which can easily be identified with the right-hand side of (\ref{e2234}).

To prove the invariance modulo 2 of $s^+(C_t)+ns^-(C_t)$, where $s^+(C_t)$ is the
number of solitary nodes of $C_t$ in
$U^+$ and $ns^-(C_t)$ is
the number of non-solitary nodes of $C_t$ in $U^-$, it is sufficient to assign to it
a topological meaning.
We do it using intersection numbers between
some auxiliary Arnol'd cycles in the Milnor ball. For simplicity, we treat
separately the following basic cases
of a nontrivial input:
a pair of complex conjugate branches, a real branch, and a pair of real branches.

The input into $s^+(C_t)$ of a pair of conjugate irreducible components $B_t,
\bar B_t$ of $C_t$ originated by a pair of complex conjugate branches $(B,z)$,
$(\bar B,z)$ of $(C,z)$
is equal to the $\Z/2$-intersection number of $(B,z)$ with the cycle formed by
$U^+$ and a half
of $E$. This follows directly from the definition of $s^+(C_t)$ as soon as we
pick that half of $E$
which does not contain the point of tangency between $B_t$ and $E$. On the other
hand, the parity of the intersection number
introduced does not depend on the choice of a half. Indeed, the sum of
the two cycles corresponding to the two choices is formed by $E$ plus twice $U^+$,
and
by assumption the intersection number of $B$ with $E$ is even. Furthermore,
the independence on the choice of a branch in a pair of complex conjugate brunches
follows now from
$\conj$-invariance
of intersection numbers. The input of the pair $(B_t, \bar B_t)$ into $ns^-(C_t)$
is zero, since there no any cross point in
$B_t\cup \bar B_t$.

Similarly, if a  real branch $(B,z)$ is contained in $U^-$, then the  input
into $(s^+(C_t)+ns^-(C_t)) \mod 2$
of local curves $B_t$ originated by $B$
is equal to the $\Z/2$-intersection number of the cycle
formed by a half of $B$ and the part bounded by $\R B$ in $U^-$ with the cycle
formed by $U^+$ and that half of $E$ which induces on $\R E$
the orientation with a direction at $z$ opposite of that induced by the
chosen half of $B$.
To check this equality, it is sufficient to consider perturbations given
by  formula (\ref{egus22}) and to note that the number of cross points
that is appearing
in $U^-$ is even if and only if the transferred direction of $\R B_t$
remains opposite to the direction $\R E$, and that the input of the tangency point is $0$
under the same assumption.

If a  real branch $(B,z)$ is contained in $U^+$, then its input into
$s^+(C_t)+ns^-(C_t)$ is equal to zero, since then, by genus argument,
the real part
has no crossing points in  $U^-$, and, by B\'ezout theorem, has no solitary
points in $U^+$.

Finally, the input of a pair of real branches $(B_1,z), (B_2,z)$ is equal
to the $\Z/2$-linking number in $\partial U$ between the boundary points
$\partial (\R B_1)\cap \partial U$ and $\partial (\R B_2)\cap \partial U$.

\smallskip

(2) In suitable local
coordinates $x,y$, we have $z=(0,0)$, $E=\{y=0\}$, $C=\{y+x^4=0\}$.
In the space $\{y+x^4\}+B_E(C,z,m)$, the local curves twice simply tangent to $E$
form the family
$$B^{eg}=\{y(1+\text{h.o.t.})+(x+\alpha)^2(x+\beta)^2(1+\text{h.o.t.})\ :\ \alpha,\beta\in(\C,0)\}\ ,$$
which is smooth of codimension $2$, and
has the tangent space given by formula (\ref{e2235}).
\proofend

\begin{lemma}\label{l637}
Let $C$ be a reduced
irreducible curve in a smooth algebraic surface
$\Sig$ such that $H^1(\Sig,{\mathcal O}_\Sig)=0$,
let $\nu:C^\vee\to C$ be the normalization map, and let ${\mathcal J}\subset
{\mathcal O}_C$ be the ideal sheaf
such that
at each point $z\in C$
with the
local branches $P_1,....,P_s$,
$s\ge1$,
it holds
$${\mathcal J}_z=\{\varphi\in
{\mathcal O}_C\ :\ \ord\varphi\big|_{P_i}\ge2\delta(P_i)+
\sum_{j\ne i} P_i\cdot P_j
+k_{z,i},\ i=1,...,s\},$$
where $k_{z,i}\ge0$ and
$\sum_{z,i}k_{z,i}<\infty$. Then,
\begin{equation}{\mathcal J}\otimes{\mathcal O}_\Sig(C)=\nu_*{\mathcal O}_{C^\vee}(\bd)\ ,\label{edi}\end{equation}
where $\deg\bd=C^2-2\sum_{z\in C}\delta(C,z)-\sum_{z,i}k_{z,i}$.
\end{lemma}

{\bf Proof.} Straightforward from \cite[Section 2.4]{CH} or
\cite[Section 4.2.4]{DI}).\proofend

\subsection{Complex nDP-pairs and their deformations}\label{sec2}

\begin{lemma}\label{l1}
Let $(\Sig,E)$ be a complex uninodal DP-pair of degree $k \geq 1$. Then:
\begin{enumerate}\item[(1)] $-K_\Sig$ is an effective divisor class
represented by a smooth elliptic curve,
\mbox{$\dim|-K_\Sig|=K_\Sig^2=k$}; \item[(2)] $-(K_\Sig+E)$ is an
effective divisor class represented by a smooth rational curve
different from $E$, and the following holds:
\begin{enumerate}\item[(2i)] $\dim|-(K_\Sig+E)|=k-1$, $(K_\Sig+E)^2=k-2$;
\item[(2ii)] if $k\ge3$, then
$-(K_\Sig+E)C\ge0$ for each irreducible curve $C$, and
$(K_\Sig+E)C=0$ only if $C$ is a $(-1)$-curve crossing $E$
transversally at one point;
\item[(2iii)] if $k=2$, then
$-(K_\Sig+E)C\ge0$ for each irreducible curve $C$;
furthermore,
$(K_\Sig+E)C=0$ if and only if $C$ is either a $(-1)$-curve crossing $E$
transversally at one point, or a
smooth rational curve
representing divisor class $-(K_\Sig+E)$; \item[(2iv)] if $k=1$,
then there exists a unique smooth rational curve
$E_0\in|-(K_\Sig+E)|$, and we have $E_0^2=-1$, $EE_0=2$,
$-(K_\Sig+E)E_0=-1$, and $-(K_\Sig+E)C\ge0$ for each irreducible
curve $C\ne E_0$. In addition, for each $(-1)$-curve $E'_0\ne E_0$,
either $E'_0E=1$, $E'_0E_0=0$, or $E'_0E=0$, $E'_0E_0=1$.
\end{enumerate}
\item[(3)] if $k=1$, then $\Sigma$ can be represented as the blow-up
of the plane at $6$ distinct points on a smooth conic and at two more points
outside the conic, while $E$ is the strict transform of that conic.\end{enumerate}
\end{lemma}

{\bf Proof.}
Claims (1) and (3) are well known and can be found, for example,  in \cite{DemazurePart3}.
All the statements in claim (2) are straightforward consequences of
the possibility to represent each uninodal DP-pair
$(\Sig, E)$
of degree $k$ such that $\Sig$ is not the quadratic cone
as the blowing up of a del Pezzo surface $\Sig'$ of degree $k+1$ at a
point belonging to exactly one
$(-1)$-curve of $\Sig'$.
In particular, the last statement in (2iv) follows from the fact
that $E'_0K_\Sig=1$ and $-(K_\Sig+E)E'_0\ge0$.\proofend

For each uninodal DP-pair $(\Sig,E)$ of degree $1$, we use the presentation
described in Lemma \ref{l1}(3). Then,
the pull-back $L$ of a generic line in the plane, the exceptional divisors
$E_1,...,E_6$ of the blown up points
of a conic, and the exceptional divisors $E_7,E_8$ of the blow-ups at the
points outside the conic form a so called {\it geometric basis}
in $\Pic(\Sig)$. In such a basis we have
\begin{equation}
\begin{cases}
&-K_\Sig=3L-E_1-...-E_8,\\
&E=2L-E_1-...-E_6, \\
&-K_\Sig-E=E_0=L-E_7-E_8.
\end{cases}
\label{e5n}
\end{equation}

An nDP-pair $(\Sig, E)$ is called {\it ridged}
if either the linear system $|-K_\Sig|$ contains
a cuspidal curve, or
the linear system
$|-2K_\Sig-E|$ contains a curve with a cusp on $E$ or a curve with a cusp in
$\Sig \setminus E$ and tangent to $E$.

We say
that a uninodal DP-pair $(\Sig,E)$ of degree $1$ possesses {\it property $T(1)$}, if:
\begin{enumerate}
\item[(i)]
$(\Sig, E)$ is not ridged,
\item[(ii)] no two $(-1)$-curves intersect $E$ at the same point;
\item[(iii)] for any divisor $D'\in\Pic_+(\Sig,E)$ such that $\dim|D'|=1$, no
irreducible rational
curve $C\in|D'|$ hits two points in the intersection of $E$ with the
union of all $(-1)$-curves of $\Sig$.
\end{enumerate}

\begin{lemma}\label{deg2}
Any uninodal DP-pair $(\Sig,E)$ of degree $k \ge2$, blown up at $k - 1$
generic points in $\Sig\setminus E$, becomes a non-tangential uninodal DP-pair of degree $1$
possessing property $T(1)$.
\end{lemma}

{\bf Proof.}
Each of the conditions (i) and (ii) in the definition of property $T(1)$
imposes non-trivial
algebraic
conditions on
the position of extra blown up points.
The same is true for the condition to be non-tangential.
Thus,
it is sufficient to show that any uninodal DP-pair of degree $2$ satisfies the condition (iii).
Let
$(\Sig,E)$ be a uninodal DP-pair of degree $2$.
Assume that there exist $D'\in\Pic_+(\Sig,E)$, $\dim|D'|=1$, and
an irreducible rational curve
$C\in|D'|$ hiting two points in the intersection of $E$ with the
union of all $(-1)$-curves of $\Sig$. Hence,
$$\frac{(D')^2-D'K_\Sig}{2}=1\quad\text{and}\quad\frac{(D')^2+D'K_\Sig}{2}+1\ge0,$$
which implies the inequalities
\begin{equation}
1\le-D'K_\Sig\le2\ .\label{el2.5}\end{equation}
Since $D'E\ge2$, we have $-(K_\Sig + E)D' \leq 0$.
Thus,
$D' = -(K_\Sig + E)$, see Lemma \ref{l1}(2iii).
However, no irreducible curve
$C \in |-(K_\Sig + E)|$
can hit any of $(-1)$-curves intersecting $E$,
see Lemma \ref{l1}(2iii).
\proofend

%In view of Lemma \ref{deg2}, we extend the class of uninodal DP-pairs possessing
%the property $T(1)$ and include in this class all uninodal DP-pairs of degree $\geq 2$.

An nDP-pair $(\Sig, E)$
is
called
\begin{itemize}
\item {\it binodal DP-pair}, if
$\Sig$ contains a
$(-2)$-curve $E'$ disjoint from $E$ and such that
$-CK_\Sig > 0$
for any reduced irreducible curve $C \ne E,E'$;
\item {\it cuspidal DP-pair},
if $\Sig$ contains a
$(-2)$-curve $E'$ such that
$EE' = 1$,
and
$-CK_\Sig > 0$
for any reduced irreducible curve
$C\ne E,E'$.
\end{itemize}

\begin{lemma}\label{negtype}
(1) Each binodal DP-pair $(\Sigma, E)$ can be viewed, after
blowing up $K_\Sigma^2-1$ generic points, as
a blow up of the plane at $8$ distinct points: $6$
points on a smooth conic $C_2$ and two points outside $C_2$ which lie
on a smooth conic $C'_2$ intersecting $C_2$ at $4$ blown up points,
so that the curves $E$ and $E'$ become the strict transforms of
$C_2$ and $C'_2$, respectively.

(2) Each cuspidal DP-pair $(\Sigma, E)$ can be viewed, after
blowing up $K_\Sigma^2-1$ generic points, as
a blow up of the plane at $8$ distinct points: $6$
points on a smooth conic $C_2$ and two points outside $C_2$ which lie
on a straight line $C_1$ intersecting $C_2$ at one blown up point,
so that the curves $E$ and $E'$
become the strict transforms of $C_2$ and $C_1$, respectively.

(3) Each tangential DP-pair $(\Sig, E)$ can be viewed as
a blow up of
the plane at $8$ distinct points: $6$
points on a smooth conic $C_2$ and two points outside $C_2$ which lie
on a straight line $C_1$
tangent to $C_2$, so that the curves $E$ and $E_0$
become the strict transforms of $C_2$ and $C_1$, respectively.
\end{lemma}

{\bf Proof.}
These statements are well known and can be easily deduced from the fact
that all del Pezzo surfaces of degree $1$ under consideration
admit $\PP^2$ as a minimal model.
The first two statements can be found,
for example, in \cite{DemazurePart3}.
\proofend

\begin{lemma}\label{l10} (1) Let $(\Sig,E)$ be a non-ridged binodal
DP-pair of degree $1$. Then,
\mbox{$\dim|-K_\Sig|=1$}, a generic curve in $|-K_\Sig|$
is smooth elliptic, and the other curves in $|-K_\Sig|$ are either
uninodal rational, or $E\cup E_0$, or $E'\cup E'_0$, where $E_0$
(resp. $E'_0$) is the unique curve in $|-(K_\Sig+E)|$ (resp. in
$|-(K_\Sig+E')|$); the curves $E_0,E'_0$ are smooth rational
$(-1)$-curves and satisfy $$EE_0=E'E'_0=2,\quad EE'_0=E'E_0=0,\quad E_0E'_0=1\ .$$

(2) Let $(\Sig,E)$ be a non-ridged cuspidal
DP-pair of degree $1$. Then,
$\dim|-K_\Sig|=1$, a generic curve in $|-K_\Sig|$
is smooth elliptic, and the other curves in $|-K_\Sig|$ are either
uninodal rational, or $E\cup E'\cup E_{-1}$, where $E_{-1}$ is a
smooth rational $(-1)$-curve intersecting each of $E,E'$ at one
point different from $E\cap E'$.

(3) Let $(\Sig,E)$ be a non-ridged tangential DP-pair. Then,
$\dim|-K_\Sig|=1$, a generic curve in $|-K_\Sig|$
is smooth elliptic, and the other curves in $|-K_\Sig|$ are either
uninodal rational, or $E\cup E_0$.
\end{lemma}

{\bf Proof.} Straightforward consequence of Lemma \ref{negtype}
and
elementary properties of plane cubics.\proofend

\begin{proposition}\label{prop1}
Given two elementary deformation equivalent uninodal DP-pairs of degree $1$,
any generic elementary deformation between these pairs
can be represented
as the blow up of the plane at $8$ points which vary in such a way that at least $6$
of them remain distinct and lie on a smooth conic.
All but
finitely many nDP-pairs in such a generic
elementary deformaton
are uninodal, non-tangential,
and have  property $T(1)$,
while the exceptional members of the deformation
are either non-tangential uninodal DP-pairs lacking property $T(1)$,
or binodal, cuspidal, or tangential DP-pairs which are non-ridged.
\end{proposition}

{\bf Proof.}
We start by proving a partial result: all but
finitely many nDP-pairs in a generic elementary deformaton
between two uninodal DP-pairs of degree $1$ are uninodal DP-pairs,
while the exceptional members of the deformation
are either cuspidal or binodal DP-pairs.

Note
that  if a
rational surface $X$ with $K_X^2>0$
contains two distinct smooth irreducible rational curves $E, D$ with
$E^2=-2$ and $D^2\le -2$, then $ED\le 1$. Indeed, if $D^2=-2$ it follows
from negative definiteness of the orthogonal complement to $K$; if
$D^2\le -3$ and $ED\ge 2$,
then, since $\dim \vert -K_X\vert\ge 1$ (as it follows from Serre duality and Riemann-Roch theorem), the divisor
$-K$ splits into
$E+D+C$, where $C$ is a nonzero effective divisor, and, since the
divisor $E+D+C$ as
any effective representative of the anticanonical divisor on a
rational surface
is connected,
we come to a contradiction due to
$C(E+D)=0\mod 2$ (by adjunction applied to $C$) and $1=p_a(-K)\ge (DE-1) +(C(E+D)-1)$.

For
any effective divisor $D$ on $X$, denote by
${\mathcal T}_{X\Vert D}$ the subsheaf of the tangent sheaf ${\mathcal T}_{X}$ generated
by vectors fields tangent to $D$, and by ${\mathcal N}'_{D/X}$
their quotient, so that we obtain the following short exact sequence of
sheaves:
$$
0\to{\mathcal T}_{X\Vert D}\to{\mathcal T}_X\to {\mathcal N}'_{D/ X}\to 0.
$$
According to the well known theory of deformations of pairs (see
\cite[Section 3.4.4]{Sernesi}),
 the long exact cohomology sequence associated to this
short sequence, and the above remark, to prove the partial result
stated at the beginning of the proof,
it is sufficient to show that $H^2({\mathcal T}_{X\Vert D})=0$ and $h^1(
{\mathcal N}'_{(D+E)/X})\ge 3$
if $D$ is either a rational irreducible curve with
$D^2\le - 3$ or $D=D_1
\cup D_2$
where
$D_i$ are $(-2)$-curves.

The equality $H^2({\mathcal T}_{X\Vert D})=0$
follows from Serre duality,
\mbox{$H^2({\mathcal T}_{X\Vert D})=
(H^0(\Omega^1_X(\log D)\otimes K))^*$} (see \cite{Deligne}),
and Bogomolov-Sommese vanishing $H^0(\Omega^1_X(\log D)\otimes K)=0$
(see \cite{Viehweg});
the latter holds in our case since $X$ is a rational
surface
with $K^2\ge 1$, and thus its anticanonical
Iitaka-Kodaira dimension is equal to $2$ (see \cite{Sakai}).

If $D^2\le -3$,
the inequality $h^1(
{\mathcal N}'_{(D+E)/X})\ge 3$
follows
from Serre-Riemann-Roch duality
and
from the exactness of
the fragment
$H^0({\mathcal N}_{D/X})\to H^1({\mathcal N}_{E/X})\to H^1({\mathcal N}'_{(D+E)/X})\to H^1(
{\mathcal N}_{D/X})$
of the long cohomology sequence
associated with the exact sequence of sheaves $0\to {\mathcal N}_{E/X}\to {\mathcal N}'_{(D+E)/X}\to
{\mathcal N}_{D/X}\to 0$.
In the second case, $D=D_1
\cup D_2$
where $D_i$ are $(-2)$-curves,
the argument is similar, but the splitting principal is to be used twice.

This proves that all but finitely many members in the family are uninodal DP-pairs,
while the exceptional ones are cuspidal or binodal DP-pairs.
Being combined with  Lemma \ref{l1} and Fujiki-Nakano-Horikawa deformation
stability of blow-ups (see \cite[Theorem 4.1]{H3})
it implies the first statement of the proposition.

It remains to check that in a generic
elementary deformation between two uninodal DP-pairs of degree $1$, first, all binodal and cuspidal members
are non-ridged, and, second,
all uninodal members are non-tangential and have property $T(1)$ except for a finite number of members that are
either non-ridged tangential DP-pairs, or uninodal DP-pairs lacking property $T(1)$.
Clearly, it is sufficient to show that each of the conditions in the definitions of tangential
DP-pairs and of property $T(1)$
imposes a finite number of proper algebraic conditions on the coordinates of the blown up points $p_1$, $\dots$, $p_8$
(we assume that the points are numbered in such a way that the points $p_1$, $\ldots$, $p_6$
lie on a conic).

For the condition to be tangential,
it is immediate, since this restriction means that the straight line through $p_7, p_8$ is tangent to the conic.

For
the condition to be non-ridged
in a generic family of uninodal, binodal, cuspidal, or
tangential DP-pairs,
we argue
by contradiction as follows.
Consider a generic member $(\Sig , E)$ of the family,
the geometric bases $L, E_1, \dots, E_8$ in $\Pic (\Sig)$ as in
 (\ref{e5n}) under the numbering of points $p_i$ in a way that $p_1, p_2$
do not belong neither to $E'$ or $E_0$ (in notation of Lemmas \ref{l1} and \ref{negtype}), and the surface $\Sig'$ obtained by contracting
$\sigma:\Sig\to\Sig'$ of  $E_1,E_2$.
If $|-K_\Sig|$ contains a cuspidal curve $C$
for any position of $p_1, p_2$,
we obtain at least a two-dimensional family of cuspidal curves $C'=\sigma(C)$
in $|-K_{\Sig'}|$, which turns the inequality (\ref{egus}) of Lemma \ref{gus} in a contradiction
$-C'K_{\Sig'}=3\ge4$.
Similarly, if $|-2K_{\Sig'}-E|=|4L-E_3-...-E_6-2E_7-2E_8|$
contains either a curve with a cusp in $\Sig -E$ and tangent to $E$, or a curve with a cusp on $E$, for any position of $p_1, p_2$,
the inequality (\ref{egus}) applied to the family of curves $C'=\sigma(C)$ obtained by variation of $p_1, p_2$
leads to a condraction: $-C'K_{\Sig'}=4\ge5$.

The complete list of $(-1)$-curves on $\Sig$,
intersecting $E$, is as follows:
\begin{equation}\begin{cases}
&E_i,\quad 1\le i\le 6,\quad L-E_i-E_j,\quad 1\le i<j,\ 7\le j\le8,\\
&2L-E_{i_1}-E_{i_2}-E_{i_3}-E_7-E_8, \quad 1\le i_1<i_2<i_2\le6,\\
&3L-\sum_{1\le i\le8,\ i\ne j,k}
E_i-2E_k,\quad 1\le j\le 6,\ 7\le k\le 8,\\
&4L-\sum_{1\le i\le6,\ i\ne j}
E_i-2E_j-2E_7-2E_8,\quad 1\le j\le 6.\end{cases}\label{e2232a}\end{equation}
This implies the finiteness for part (ii) of property $T(1)$.
To show the properness of this restriction, observe
that if two curves
$C', C''$ belonging to the list
intersect, then there is $E_i$ such that
$C'E_i=1$ and $C''E_i=0$ (up to permutation of $C'$ and $C''$). Then, we blow down $E_i$
by $\sigma_i:\Sig\to\Sig_i$ and obtain that $(\sigma_i(C'))^2=0$ and $\dim|\sigma_i(C')|=1$.
Hence, shifting the point $\sigma_i(E_i)$ in $\Sig_i$ (along $\sigma_i(E)$ is $1\le i\le 6$)
and blowing up the shifted point, we make the sets $C'\cap E$ and $C''\cap E$ disjoint.

In
part (iii) of
property $T(1)$ we have to consider divisor
classes presented by irreducible rational curves with $D'E\ge2$ and $\dim|D'|=1$.
Together with adjunction inequaity, $2\ge (D')^2-D'K_\Sig$, it yields that $(D')^2=0$
and $-D'K_\Sig=D'E=2$. Writing
$D'=dL-d_1E_1-...-d_8E_8$ in the basis (\ref{e5n}), one can easily extract from the
preceding relations that
$d<8$,
and hence we have to deal with only finitely many divisor
classes. Furthermore, the Cremona base
changes in $\Pic(\Sig)$
\begin{equation}\begin{cases}&L=2L'-E'_i-E'_j-E'_7,\\
&E_i=L'-E'_j-E'_7,\\ &E_j=L'-E'_i-E'_7,\\
&E_7=L'-
E'_i-E'_j,\end{cases}\quad1\le i<j\le 6,\label{e2232}\end{equation}
combined with permutations of $E_7$ and $E_8$ bring any of the considered divisor classes
to the form $D'=L-E_7$. Indeed, assuming that $d_1\ge...\ge d_6$, $d_7\ge d_8$, and $d_1+d_2+d_7>d$,
we obtain via the transformation (\ref{e2232}) specified to $i=1$, $j=2$, that $D'=d'L'-d'_1E'_1-...$ with
$d'=2d-d_1-d_2-d_7<d$. Using the relation $D'E=2$ and the above reduction, we end up with
two minimal expressions $D'=L-E_7$ or $D'=4L-E_1-...-E_6-2E_7-2E_8$, where the latter one does not meet the condition
$(D')^2=0$.
Next, it is easy to check that, for any two $(-1)$-curves
$C', C''$ in the list
(\ref{e2232a}) such that $C'D'>0$, $C''D'>0$, there exists $E_i$, $1\le i\le 6$, satisfying
$E_iC'=1$ and $E_iC''=E_iD'=0$ (up to permutation of $C',C''$). Then, varying $E_i$ and
$C'$ as in the preceding paragraph and keeping $D'$ and $C''$ fixed, we can make
$C'\cap E$ disjoint from the curve $C\in|D'|$ passing through $C''\cap E$.
\proofend

\subsection{Rational curves on uninodal, binodal, cuspidal, and tangential DP-pairs}\label{sec2.3}

In this section, we consider
uninodal, binodal, cuspidal, and
tangential DP-pairs $(\Sig,E)$ of degree $1$.
For any $D\in\Pic_+(\Sig,E)$, the moduli space $\overline{\mathcal M}_{0, r}(\Sig,D)$
of stable pointed maps $(\bn:\widehat C\to\Sig, \bp)$ of connected curves $\widehat C$ of genus $0$ such that
$\bn_*\widehat C\in|D|$ is a projective variety (see \cite[Theorem 1]{FP}). We deal with its subvariety
$\overline{\mathcal M}^*_{0, r}(\Sig,D)$  that is, by definition, the closure of ${\mathcal M}_{0, r}(\Sig,D)$.
We use the notation $\rho$ for the natural morphism
 \begin{equation}
 \rho:\overline{\mathcal M}^*_{0, r}(\Sig,D)\to|D|,\quad
[\bn:\widehat C\to\Sig, \bp]\mapsto
\bn_*\widehat
C\;.\label{e2233}
\end{equation}
and the notation $\pi_i$ for the forgetful morphism
\begin{equation}
 \pi_i:\overline{\mathcal M}^*_{0, r}(\Sig, D) \to \overline{\mathcal M}^*_{0, r - i}(\Sig,D),\quad
\label{e2233-new}
\end{equation}
provided by removing the last $i$ marked points.

Furthermore, for any $D\in\Pic_+(\Sig,E)$
and any nonnegative integer $l\le DE/2$,
denote by ${\mathcal V}^l_r(\Sig,E,D)$
the subset of ${\mathcal M}_{0, r}(\Sig,D)$
consisting of elements \mbox{$[\bn:\PP^1\to\Sig, \bp]$} subject to the following
restriction: $\bn^*(E)=2\bd_0+\bd'$,
where
$\bd_0,\bd'\in\Div(\PP^1)$ are effective divisors,
$\deg\bd_0=l$. Respectively, by $\overline{\mathcal V}^{\;l}_r(\Sig,E,D)$
we denote the closure of ${\mathcal V}^l_r(\Sig,E,D)$ in
$\overline{\mathcal M}^*_{0, r}(\Sig,D)$. To simplify notations, we
write ${\mathcal V}_r(\Sig,E,D)$ in the case of $\bd'=0$ ({\it i.e.}, $l=DE/2$),
and we abrreviate  ${\mathcal V}_0(\Sig,E,D)$, ${\mathcal V}^l_0(\Sig,E,D)$ and $\overline{\mathcal V}^{\;l}_0(\Sig,E,D)$)
to ${\mathcal V}(\Sig,E,D)$, ${\mathcal V}^l(\Sig,E,D)$ and $\overline{\mathcal V}^{\;l}(\Sig,E,D)$, respectively,
and write $[\bn:\widehat C\to\Sig]$ (instead of $[\bn:\widehat C\to\Sig, \emptyset]$)
for their elements.
Put
$$
r(\Sig,D,l)=-DK_\Sig-1-l.
$$

For an irreducible family
${\mathcal V}\subset\overline{\mathcal M}^*_{0, r}(\Sig,D)$, we set $\idim{\mathcal V}=\dim\rho({\mathcal V})$; the latter
numerical characteristic
can be viewed as the maximal number of generic points in $\Sig$
through which one can trace a curve $C=\bn_*\widehat C$ for some
$[\bn:\widehat C\to\Sig, \bp]\in{\mathcal V}$. We say that an irreducible family
${\mathcal V}\subset\overline{\mathcal V}^{\;l}_r(\Sig,E,D)$ is {\it equisingular}
if
all the curves $C=\bn_*\widehat C$, where
$[\bn:\widehat C\to\Sig, \bp] \in {\mathcal V}$,
split into distinct irreducible components
$C=m_1C_1\cup...\cup m_sC_s$
with the same multiplicities $m_1,...,m_s$
and the topological types of the singular points of the curves
$C_1\cup....\cup C_s$ persist in the induced family.

\subsubsection{Codimension zero: the case of uninodal DP-pairs}\label{sec-zero}
Consider a uninodal DP-pair  $(\Sig,E)$  of degree $1$ and, as above, denote by $E_0$ both
the divisor class $-(K_\Sig + E)$
and the unique curve belonging to this divisor class,
see Lemma \ref{l1}.
We fix a divisor $D\in\Pic_+(\Sig,E)$ and an integer $l$, $0\le l\le DE/2$.

\begin{lemma}\label{lp1-1}
If $D=sE_0$
then either
${\mathcal
V}^l(\Sig,E,D)=\emptyset$ or $\idim{\mathcal
V}^l(\Sig,E,D)=0$. If $D\ne sE_0$
and ${\mathcal V}^l(\Sig,E,D)\ne\emptyset$ then
$r(\Sig,D,l)\ge0$.
\end{lemma}

{\bf Proof.}
The former statement is trivial. The
latter
one follows from
the identity
\begin{equation}r(\Sig,D,l)=-\frac{DK_\Sig}{2}-\frac{D(K_\Sig+E)}
{2}+\frac{DE-2l}{2}-1\
, \label{e4}\end{equation} since the first summand is positive and
next two are non-negative (see Lemma \ref{l1}).
\proofend

\begin{lemma}\label{lp1-2}
Let  $D\ne sE_0$,
${\mathcal V}^l(\Sig,E,D)\ne\emptyset$, and
$r(\Sig,D,l)=0$. Then
\begin{enumerate}\item[(1)] $\idim{\mathcal V}^l(\Sig,E,D)=0$.
\item[(2)] If the pair $(\Sig,E)$ has property $T(1)$, the elements
$[\bn:\PP^1\to\Sig]\in {\mathcal V}^l(\Sig,E,D)$ are as follows:
\begin{enumerate}
\item[(2i)] either $-DK_\Sig=1$, $0\le DE\le 1$, $l=0$, and $\bn$
takes $\PP^1$ birationally onto a smooth or uninodal curve
$C\in|D|$;
\item[(2ii)] or $-DK_\Sig=2$, $DE=2$, $l=1$, and $\bn$ takes $\PP^1$
birationally onto a
smooth or uninodal rational curve $C\in|D|$
simply tangent to
$E$ at one point;
\item[(2iii)] or $D=2C$, $C$ is a
$(-1)$-curve crossing $E$ transversally at one point, $\bn:\PP^1\to
C$ is a double covering with two ramification points, one of which
is the intersection point $C\cap E$.
\end{enumerate} All the curves in (2i)-(2iii) are
disjoint from $E\cap E_0$.

\item[(3)]
If the pair
$(\Sig,E)$
lacks property $T(1)$, the elements
$[\bn:\PP^1\to\Sig]\in {\mathcal V}^l(\Sig,E,D)$ different from those mentioned in (2i)-(2iii),
are as follows: \begin{enumerate}\item[(3i)] $\bn$ takes $\PP^1$
birationally onto a rational curve $C\in|-K_\Sig+E_0|$ which either
has a cusp in $\Sig\setminus E$ and is
simply tangent to $E$
at one point, or has a cusp on $E$;
\item[(3ii)] $\bn$ takes $\PP^1$ birationally onto a cuspidal curve
$C\in|-K_\Sig|$.
\end{enumerate}\end{enumerate}
\end{lemma}

{\bf Proof.}
Formula (\ref{e4}) and $r(\Sig,D,l)=0$ yield
\begin{eqnarray} &\text{(a)}&\ \text{either}\ -DK_\Sig=2,\ DE=2, \ l=1,
\nonumber\\ &\text{(b)}&\
\text{or}\ -DK_\Sig=1,\ DE=1,\ l=0,\label{eR0}\\ &\text{(c)}&\ \text{or}\
-DK_\Sig=1,\ DE=0,\ l=0.\nonumber\end{eqnarray}
Choose
the standard basis (\ref{e5n}) in
$\Pic(\Sig)$.
Note
that for $D=sE_i$, $1\le i\le8$, $s\ge1$, we have $1\le s\le2$, and
then the statement (2) is immediate.
Thus, we can suppose that
$D=dL-d_1E_1-...-d_8E_8$ with $d>0$, $d_1,...,d_8\ge0$.
In the situation (a),
the Cremona base
change (\ref{e2232}) used as in the proof of Proposition \ref{prop1}
brings $D$ to the form
$D= L-E_7$, or
$D=4L-E_1-...-E_6-2E_7-2E_8$,
where $\bn:\PP^1\to\Sig$ is birational onto its image, or to the form $2E_i$,
$1\le i\le 6$, where $\bn:\PP^1\to\Sig$ is a double covering of $E_i$ considered above.
This implies the statement (2) in Case (a).

In the same manner, we can see that, in Case
(b), $\bn$ takes $\PP^1$ isomorphically onto a $(-1)$-curve crossing
$E$ transversally at one point.

At last, in case (c), $\bn$ must be
birational onto its image, and, in the above notation, we have
$$2d=d_1+...+d_6,\quad d=d_7+d_8+1\;,$$ which together with
inequality $d_1+d_2+d_7\le d$ yields $d\le3$. Thus, if $d\le2$, $\bn$
takes $\PP^1$ isomorphically onto a $(-1)$-curve disjoint from $E$.
If $d=3$, $\bn$ takes $\PP^1$ birationally onto a curve
$C\in|-K_\Sig|$ disjoint from $E$ and having a node or a cusp.

The relation $\idim{\mathcal V}^l(\Sig,E,D)=0$ is evident in all the
considered cases except for $D=4L-E_1-...-E_6-2E_7-2E_8=-K_\Sig+E_0$, $[\bn:\PP^1\to \Sig]
\in{\mathcal V}^l(\Sig,E,-K_\Sig+E_0)$ with $\bn$ birational onto its image $C$,
and $\Card(\bn^{-1}(E))=1$. In
such an exceptional
case,
the assumption $\idim V>0$ leads to a contradiction, since
then (\ref{egus}) implies
$$-DK_\Sig\ge2
+DE-
\Card\, \bn^{-1}(E)=
 -DK_\Sig+1\;.$$

The relation $C\cap E\cap E_0=\emptyset$ follows from the fact that in
all the cases (2i)-(2iii) either $CE=0$ or $CE_0=0$.

For statement (3) we repeat the above analysis and come either to the
case $D=-K_\Sig+E_0$
or to the case of $D=-K_\Sig$, which then must be as asserted in Lemma.
\proofend

\begin{lemma}\label{lp1-3}
Let  $D\ne sE_0$,
${\mathcal V}^l(\Sig,E,D)\ne\emptyset$, and
$r(\Sig,D,l)>0$. Then, for any irreducible component
${\mathcal V}\subset{\mathcal V}^l(\Sig,E,D)$ whose generic element
is represented by a map $\bn:\PP^1\to\Sig$ birational onto its image
$C=\bn(\PP^1)$,
we have $\idim{\mathcal V}=r(\Sig,D,l)$ and the following properties:
\begin{enumerate}\item[(i)] the family ${\mathcal V}$ has no base
points,
\item[(ii)]
$\bn$ is an immersion outside of $E$,
\item[(iii)] the divisor $\bn^*(E)$
is supported
on
$DE-l$ distinct points so that $l$ of them
have multiplicity $2$ and the others multiplicity $1$,
\item[(iv)]
either $C$ is disjoint from $E_0$, or it has $DE_0$ local branches
centered on $E_0\setminus E$ and intersecting $E_0$ with
multiplicity $1$,
\item[(v)] if $r\ge2$, then $C$ is immersed and smooth along $E\cup E_0$, it
intersects $E_0$ at $DE_0$ distinct points, and intersects $E$ at $DE-l$ distinct points,
transversally at $DE-2l$ of them and with
simple tangency
at the other $l$ ones.\end{enumerate} Moreover, if $C'\subset\Sig$ is any reduced, irreducible curve
different from $E$, the subset of elements $[\bn:\PP^1\to\Sig]\in{\mathcal V}$ such that
$\bn^*(C')$ consists of only simple points, is open and dense.
\end{lemma}

{\bf Proof.}
Suppose that $\idim{\mathcal V}>r(\Sig,D,l)$. Abbreviate
$r=r(\Sig,D,l)$.
By (\ref{egus}), applied to a generic part of ${\mathcal V}$ (which
is equisingular, since being equigeneric),
$$-DK_\Sig\ge2+
DE - \Card\,\bn^{-1}(E)+r
\ge -DK_\Sig+1\;,$$
a
contradiction. To prove that $\idim{\mathcal V}=r$, we need only to
show that $\idim{\mathcal V}\ge r$, but the latter immediately follows from
\cite[Theorem II.1.2]{K} and from the fact that the condition
$\bn^*(E)=2\bd_0+\bd'$,
where
$\bd_0,\bd'\in\Div(\PP^1)$ are effective divisors,
$\deg\bd_0=l$, reduces
the intersection dimension by at most $l$.

To establish required geometric properties of the curve
$C$, we again apply (\ref{egus}) which takes now the form
\begin{equation}
-DK_\Sig\ge 2+DE-\Card\, \bn^{-1}(E)+r-1+\sum(\ord P-1)\ge -DK_\Sig+\sum(\ord P-1)\;,\label{egen}\end{equation}
where $P$ runs over all singular local branches of $C$ in
$\Sig\setminus E$. This yields that
$DE-\Card\,\bn^{-1}(E)=l$, which, in particular, means that $C$ has $l$
local branches
centered along $E$ and intersecting $E$ with multiplicity $2$ and each
of the other $DE-2l$
local branches intersects
$E$ with multiplicity $1$. Furthermore,
$\bn$
must be an
immersion outside of $E$,
the curve $C$ avoids the set $E\cap E_0$ and cannot have
local branches
centered on $E_0$ and intersecting $E_0$ with multiplicity $>1$,
since otherwise one would have an extra positive contribution to the
right-hand side of (\ref{egen}), which is a contradiction.
A similar reasoning proves claim (v).

Assume now that $r\ge2$ and suppose that $C$ has a
singular point on $E\cup E_0$, {\it i.e.} has multiplicity $s\ge2$
at some point on $E\cup E_0$. Fixing the position of this point we
obtain a subfamily ${\mathcal V}'\subset{\mathcal V}$ of dimension
$\idim{\mathcal V}'\ge\idim
{\mathcal V}-1=r-1>0$. Applying inequality (\ref{egus}) to the family
${\mathcal V}'$, we
get a contradiction:
$$-DK_\Sig\ge2+s+DE-\Card\,\bn^{-1}(E)+r-2
\ge
-DK_\Sig+1\;,$$ which proves the last statement.
\proofend

\begin{definition}\label{d181015}
Denote by
${\mathcal V}^{l, im}(\Sig,E,D)$
the union of the sets ${\mathcal V}^{im}$ over all irreducible components ${\mathcal V}
\subset{\mathcal V}^l(\Sig,E,D)$ of intersection dimension $r(\Sig,D,l)$
whose generic element is represented by a map $\bn:\PP^1\to\Sig$ birational
onto its image, where ${\mathcal V}^{im}\subset{\mathcal V}$ is the
{\rm (}open{\rm )} subset formed by the elements $[\bn:\PP^1\to\Sig]$ satisfying
properties (ii) and (iii) of Lemma \ref{lp1-3}.
Denote by $\overline{{\mathcal V}}^{\;l, im}(\Sig,E,D)$ the closure of
${\mathcal V}^{l, im}(\Sig,E,D) \subset \overline{\mathcal M}^*_{0, 0}(\Sig, D)$.
For $l = DE/2$, we abrreviate ${\mathcal V}^{l, im}(\Sig,E,D)$
and $\overline{{\mathcal V}}^{\;l, im}(\Sig,E,D)$ to ${\mathcal V}^{im}(\Sig,E,D)$
and $\overline{{\mathcal V}}^{\;im}(\Sig,E,D)$, respectivement.
\end{definition}

\begin{lemma}\label{lp1-4}
Let $D\ne sE_0$,
${\mathcal V}^l(\Sig,E,D)\ne\emptyset$, and
$r(\Sig,D,l)>0$. Let
${\mathcal V}\subset{\mathcal V}^l(\Sig,E,D)$ be an irreducible component with $\idim{\mathcal V}\ge
r(\Sig,D,l)$,
whose generic element
is represented by a map $\bn:\PP^1\to\Sig$
that
multiply covers its image
$C=\bn(\PP^1)$.
Then
$\idim{\mathcal V} = r(\Sig,D,l)$, and
$\bn:\PP^1\to C$ is a ramified double covering, $D=2D'$, and
\begin{enumerate}\item[(i)] either $C$ is an embedded smooth curve with $CE=2$
and $C^2=0$ intersecting $E$ at two distinct points,
\item[(ii)] or $D'=-K_\Sig+E_0$ and $C
\in|D'|$ is a uninodal
rational curve having along $E$ two smooth local branches that cross $E$
transversally.\end{enumerate}
In both cases, $CE_0=0$, the ramification set is $C\cap E$, and $r(\Sig,D,l)=1$.
\end{lemma}

{\bf Proof.}
Let a generic member $[\bn:\PP^1\to\Sig]\in{\mathcal V}$ be
an $s$-multiple cover, where $s\ge2$. Then $D=sD'$.
The family of reduced irreducible rational curves in $|D'|$ has
dimension $-D'K_\Sig-1$ or is empty (see Lemmas \ref{lp1-1} - \ref{lp1-3}).
Hence,
$$r(\Sig,D,l)=-DK_\Sig-1-l \leq \idim{\mathcal V} \leq -D'K_\Sig-1\ ,$$
and thus,
$$
-(s-1)D'K_\Sig\le l\le\frac{s}{2}D'E\ .$$
In view of $-D'K_\Sig>0$ and $-D'(K_\Sig+E)\ge0$, we get
$s=2$, $D'(K_\Sig+E)=0$, and $l=D'E$. Furthermore,
$\idim{\mathcal V} = r(\Sig, D, l)$ and
by Lemma \ref{lp1-3} the image $C=\bn(\PP^1)$ is an
immersed rational curve intersecting $E$ at $CE=D'E$ distinct points
transversally.

Since $l=D'E$, the map $\bn$ must have ramification at $C \cap E$,
and therefore, the Riemann-Hurwitz formula yields $1\le CE\le 2$. Next, the
assumption $r(\Sig, D, l) >0$ implies $CE=2$ and $r(\Sig, D, l) = 1$.
We have $D'$
satisfying $D'E=2$ and $D'E_0=0$.
The Cremona transformations
(\ref{e2232}) lead to the two options:
either $D'= L-E_7$,
or $D'=4L-E_1-...-E_6-2E_7-2E_8$ ({\it cf.} the proof of Proposition
\ref{prop1} and Lemma
\ref{lp1-2}).
In the former
case, the curve $C$ is smooth;
in the
latter case,
$C$ has a node in $\Sig\setminus E$.
\proofend

The following statement is an immediate
consequence of Lemmas \ref{lp1-1}, \ref{lp1-3}, and \ref{lp1-4}.

\begin{lemma}\label{lem-fin}
Let $(\Sig,E)$ be a uninodal DP-pair, and let $D\in\Pic_+(\Sig,E)$ be a
divisor class matching conditions
(\ref{e1}).
Put $l=DE/2$.
If $r(\Sig,D,l)=0$, then the set ${\mathcal V}(\Sig,E,D)$ is finite. If $r=r(\Sig,D,l)>0$
and $\bw$ is a generic $r$-tuple
of points of $\Sig$, then the set
$$\{[\bn:\PP^1\to\Sig]\in{\mathcal V}(\Sig,E,D)\ :\ \bw\subset\bn_*(\PP^1)\}$$ is finite.
\proofend
\end{lemma}

\subsubsection{Codimension zero: the case of binodal, cuspidal, and
tangential DP-pairs}\label{sec-zero1}

The following properties of general members of families of curves on
binodal, cuspidal,
and tangential DP-pairs are used in the proof of Theorem \ref{t1} in
Section \ref{sec-t1}. They are very similar to those for families of curves on
uninodal DP-pairs.

\begin{lemma}\label{p3} Let $(\Sig,E)$ be a binodal {\rm (}resp. cuspidal, resp.
tangential {\rm )}
DP-pair of degree $1$.
Assume that $(\Sig, E)$ is ridged.
Let
$D\in\Pic_+(\Sig,E)$ and $0\le l\le DE/2$.

If $D=sE_0$ or $sE'$ {\rm (}resp. $D=sE_{-1}$
or $sE'$,
resp. $D=sE_0${\rm )} with $s\ge1$,
then either ${\mathcal
V}^l(\Sig,E,D)=\emptyset$
or $\idim{\mathcal
V}^l(\Sig,E,D)=0$.

Let $D\ne sE_0$, $sE'$ {\rm (}resp. $D\ne E_{-1}$, $sE'$, resp. $D\ne sE_0${\rm )}
with $s\ge1$. If ${\mathcal V}^l(\Sig,E,D)\ne\emptyset$ and a irreducible component ${\mathcal V}\subset
{\mathcal V}^l(\Sig,E,D)$ satisfies $\idim{\mathcal V}\ge r(\Sig,D,l)$, then
$\idim V=r(\Sig,D,l)$, and
a generic element
$[\bn:\PP^1\to\Sig]\in{\mathcal V}$ is as follows:
\begin{enumerate}\item[(i)] either $\bn$ birationally takes $\PP^1$ onto
its image so that

- it is an immersion outside $E$,

- the divisor $\bn^*(E)$ consists of $DE-l$ points so that $l$ of them have
multiplicity $2$ and the other ones are simple,

- the divisor $\bn^*(E')$ consists of $DE'$ simple points,

- the divisor $\bn^*(E_0)$ {\rm (}resp. $\bn^*(E_{-1})${\rm )} consists of $DE_0$
{\rm (}resp. $DE_{-1}${\rm )} simple points;
\item[(ii)] or $D=2C$, $l=1$, $C$ is a
$(-1)$-curve crossing $E$ transversally at one point, $\bn:\PP^1\to
C$ is a double covering with two ramification points, one of which
is the intersection point of $C$ and $E$, \item[(iii)] or $D=2C$, $l=2$, $C$ is a
smooth or uninodal curve satisfying $-(K_\Sig+E)C=0$ and
transversally intersecting $E$ at two distinct points, and
$\bn:\PP^1\to C$ is a double covering ramified at
$C\cap E$.
\end{enumerate} Furthermore, the curve $\bn(\PP^1)$ does not hit the points
of $E\cap(E_0\cup E')$ {\rm (}resp. $E\cap(E_{-1}\cup E')$, resp. $E\cap E_0${\rm )}.
\end{lemma}

{\bf Proof.} All the statements can be established using the same
argumentation as
in the proof of Lemmas \ref{lp1-1} - \ref{lp1-3} and \ref{lp1-4}.
It literally applies to the tangential case.
The only statement whose proof requires modification is the case
(i) for binodal or cuspidal DP-pairs.

\smallskip

Suppose that $r=r(\Sig,D,l)>0$ and $\bn:\PP^1\to
C=\bn(\PP^1)\subset\Sig$ birational. Lemma \ref{gus} yields:
\begin{itemize}
\item for $(\Sig, E)$ binodal,
$$-DK_\Sig\ge2+DE-
\Card\, \bn^{-1}(E)+DE'-\Card\, \bn^{-1}(E')+\sum_P(\ord P-1)+r-1$$
$$=-DK_\Sig+(DE-
\Card\, \bn^{-1}(E)-l)+(DE'-
\Card\, \bn^{-1}(E'))+\sum_P(\ord P-1)
\ ,$$ where $P$ runs over all singular branches of $C$ in
$\Sig\setminus(E\cup E')$.
\item for $(\Sig, E)$ cuspidal,
$z=E\cap E'$,
\begin{eqnarray}& & -DK_\Sig
\ge 2
+r-1+\sum_P(\ord P-1)+(DE-(C\cdot
E)_z-\nonumber\\& &
\Card\, \bn^{-1}(E\setminus\{z\}))
+(DE'-(C\cdot
E')_z-
\Card\, \bn^{-1}(E'\setminus\{z\}))+
\nonumber\\& &
2\cdot\ord(C,z)-
\Card\, \bn^{-1}(z)
+((C\cdot e)_z-\ord(C,z))+((C\cdot
E')_z- \nonumber\\& &
\ord(C,z))=
-DK_\Sig
+\sum_P(\ord P-1)+
\nonumber\\
& &(DE-
\Card\, \bn^{-1}(E)-l)+(DE'-
\Card\, \bn^{-1}(E'))+
\Card\, \bn^{-1}(z)\;,\nonumber\end{eqnarray} where
$P$ runs over all singular branches of $C$ in $\Sig\setminus(E \cup E')$.
\end{itemize}
These inequalities imply the required properties of $\bn:\PP^1\to\Sig$.
\proofend

\subsubsection{Codimension one}\label{sec-one}
Throughout this section we
consider non-tangential uninodal DP-pairs $(\Sig,E)$ having degree $1$ and possessing property $T(1)$.

\begin{lemma}\label{lp1-5} Let  $D\ne sE_0$,
let
$r(\Sig,D,l)>0$,
and let ${\mathcal V}\subset{\mathcal V}^l(\Sig,E,D)$
be a non-empty irreducible component of the family defined by the
condition that, for
all
$[\bn:\PP^1\to\Sig]\in{\mathcal V}$, the divisor $\bn^*(E)$ has
$p\le DE$
components, and the images of $m$ of them
are fixed on $E$, $1\le m\le p$. Suppose that $\idim{\mathcal V}\ge-D(K_\Sig+E)-1+p-m$.
Then: \begin{enumerate}
\item[(i)] if
 $D(K_\Sig+E)=0$ and $m<DE-l$,
 or if $-D(K_\Sig+E)>0$,
we have $\idim{\mathcal V}=-D(K_\Sig+E)-1+p-m$,
\item[(ii)] if $D(K_\Sig+E)=0$ and $m=DE-l$, we have
$\idim{\mathcal V}=0$.\end{enumerate}
\end{lemma}

{\bf Proof.}
Let a generic element ${\mathcal V}$ be presented by
a map $\bn:\PP^1\to\Sig$ birational onto its image $C$.
Assume that $-D(K_\Sig+E)>0$.
Then
$-D(K_\Sig+E)-1+p-m\ge0$.
Suppose
that $\idim{\mathcal V}>-D(K_\Sig+E)-1+p-m$. Then inequality (\ref{egus}) applied to
the family ${\mathcal V}$ results in a contradiction:
$$-DK_\Sig\ge2+(DE-p
+m)+(-D(K_\Sig+E)-1+p-m)=
-DK_\Sig+1\;.$$
The same argument settles the case of $D(K_\Sig+E)=0$. If a generic
element of ${\mathcal V}$ is a multiple covering, then the required
statement follows from Lemma \ref{lp1-4}.
\proofend

\begin{lemma}\label{lp1-6}
Let $D=-2K_\Sig+E_0$.
Then,
for $DE=2$, the space $\overline{\mathcal V}(\Sig,E,D)$
contains
no element $\bn:\widehat C\to\Sig$ such that $\bn_*(\widehat C)$ is supported at $E\cup E_0$.
\end{lemma}

{\bf Proof.}
We argue by contradiction.
Suppose that $[\bn:\widehat C_1\cup\widehat C_2\to\Sig]
\in\overline{\mathcal V}(\Sig, E, D)$
is mapping $\widehat C_1$
onto $E$ and $\widehat C_2$
onto $E_0$. Denote by $c_i$ the number of irreducible
components of $\widehat C_i$, $i=1,2$,
and denote by $z_1,z_2$ the two intersection points of $E$ and $E_0$.

Note that $-2K_\Sig+E_0=2E+3E_0$, and hence $c_1\le2$, $c_2\le 3$.

Consider, now, $[\bn:\widehat C\to\Sig]$ as a limit of a family $[\bn_t: \PP^1\to\Sig]$, $t\ne 0$,
belonging to ${\mathcal V}(\Sig,E,D)$.

Suppose that $c_1=1$. Then $\bn:\widehat C_1\to E$ is a double covering
ramified at two points. Since
$DE_0=1$, the only possible structure of $\bn:\widehat C\to E\cup E_0$ is as follows:
$\bn:\widehat C_1\to E$ is ramified at $z_1$ and $z_2$, and $2\le c_2\le 3$, while
one of the components of $\hat C_2$ is attached to $\widehat C_1$ at
$(\bn\big|_{\widehat C_1})^{-1}(z_1)$ and another one
at $(\bn\big|_{\widehat C_1})^{-1}(z_2)$. However, then a generic curve $\bn_t(\PP^1)$, $t\ne0$,
must intersect $E$ with multiplicity $\ge3$, contrary to $DE=2$.

Thus, $c_1=2$, and $\widehat C_1$ consists of two disjoint
components isomorphically mapped onto $E$. If $c_2=3$, then one of the three
irreducible components
$\widehat C'_2$, $\widehat C''_2$, $\widehat C'''_2$ of $\widehat C_1$, say
$\widehat C'_2$, meets both the components $\widehat C'_1$, $\widehat C''_1$, while
the other two meet each only one (otherwise it would contradict the condition
$DE=2$); this implies that $\widehat C'_1$, $\widehat C''_1$ are disjoint,
which in its turn implies that the restriction of $\bn$ to the germ of each of
$\widehat C'_1$, $\widehat C''_1$ in $\widehat C$ (such a restriction being considered
as a relative to boundary cycle in a small tubular neighborhood of $E$) has the
total intersection number with $E$ equal to $0$; the remaining
two intersection points of $\widehat C''_2\cup \widehat C'''_2$ with $E$ are
transversal, hence $n_t$ with small $t\ne 0$ can not be tangent to $E$; contradiction.

Thus, $c_2=2$ and $\widehat C_2$ consists of two irreducible components, one
of them, $\widehat C'_2$, is mapped with degree $2$
onto $E_0$, and the other one, $\widehat C''_2$, is mapped onto $E_0$ isomorphically.
To avoid transversal intersections with $E$,
the component $\widehat C''_2$ should be
joined by
a node both with $\widehat C'_1$ and $\widehat C''_1$, the two ramification points
of $\bn$ restricted to $\widehat C'_2$  should be the points $(\bn\big|_{\widehat
C'_2})^{-1}(z_1)$ and
$(\bn\big|_{\widehat C'_2})^{-1}(z_2)$,
and  $\widehat C'_2$  should be joined
at one of the branching points, say, $(\bn\big|_{\widehat C''_2})^{-1}(z_1)$,
by a node with one component of
$\widehat C_1$, say, with $\widehat C'_1$.
However, then the relative to boundary cycle realized in a small tubular
neighborhood of $E$
by the
restriction of $\bn$ to a neighborhood of $\widehat C''_1$ in $\widehat C$
has a negative total intersection number with $E$
(since $\widehat C''_1$ is adjacent in $\widehat C$ only to $\widehat C''_2$,
which is mapped onto $E_0$ isomorphically, this intersection number is equal
to $-2+1=-1$), which contradicts the existence of a smoothing family $\bn_t$.\proofend

Starting from here, and up to end of section, we introduce the following additional
assumptions:
\begin{itemize}\item $D\in\Pic_+(\Sig,E)$;
\item
$DE>0$ is even,
and $r=r(\Sig,D,l)\ge1$;
\item ${\mathcal V}(\Sig,E,D)\ne\emptyset$.
\end{itemize}
Notice here that the assumptions $DE>0$ and $r=-DK_\Sig-\frac{DE}{2}-1\ge1$ together with the genus inequality
\mbox{$(D^2+DK_\Sig)/2+1\ge0$} yield
\begin{equation}D^2\ge-2-DK_\Sig\ge\frac{DE}{2}>0\ .
\label{ed2g0}\end{equation} Furthermore,
in view of Lemma \ref{lp1-3} and property $T(1)$, ${\mathcal V}^{im}(\Sig,E,D)$ is an open dense subset of
${\mathcal V}(\Sig,E,D)$.

\begin{lemma}\label{lp2-1}
Let ${\mathcal V}\subset
\overline{\mathcal V}(\Sig,E,D)\setminus{\mathcal V}^{im}(\Sig,E,D)$ be an
irreducible equisingular family such that $\idim{\mathcal V}=r-1$, and let $l=DE/2$. If
$[\bn:\PP^1\to\Sig]$
is a generic element of ${\mathcal V}$, then the map $\bn$ and its image
$C=\bn(\PP^1)$ are as follows:
\begin{enumerate}\item[(i)] the map $\bn:\PP^1\to C$ is birational
onto $C$, but is not an immersion outside $\bn^*(E)$, and the divisor
$\bn^*(E)$ consists of $l$ distinct double points;
\item[(ii)] or $l\ge2$, the map $\bn:\PP^1\to C$ is birational
onto $C$,
and the divisor $\bn^*(E)$ consists of $l-2$ double points
and one more point of multiplicity $4$; furthermore, in this case, if $r\ge3$,  then $\bn$ is an immersion;
\item[(iii)] or $l=2$, $\bn:\PP^1\to
C$ is a double covering, $D=2C$, $-CK_\Sig=3$, $CE=2$, and $C$ is an
immersed rational curve transversally intersecting $E\setminus E_0$
at two points, which are ramification points of the covering.
\end{enumerate}
\end{lemma}

{\bf Proof.}
The statement is straightforward if $-K_\Sig-D$ is effective. Thus, we suppose in
the sequel that $-K_\Sig-D$ is not effective.

Let
$\bn:\PP^1\to C$ be birational and $r\ge2$. Inequality (\ref{egus})
yields
$$-DK_\Sig\ge2+(r-2)+(DE-\Card\, \bn^{-1}(E))
+\sum_Q(\ord Q-1)$$
\begin{equation}=-DK_\Sig+\left(\sum_Q(\ord Q-1)\right)+\left(\frac{DE}{2}-
\Card\, \bn^{-1}(E)
\right)
-1\;,\label{e18}\end{equation} where
the
sum
runs over the singular local branches of $\bn$ centered in
$\Sig\setminus E$.
Statement (i) and statement
(ii), except for the condition on $\bn$ to be an immersion everywhere, follow
immediately.
Suppose that $r\ge3$. Fixing the position of $z=\bn(p)\in E$, where $4p\le\bn^*(E)$, we obtain a subfamily of dimension $r-2\ge1$, and again
apply inequality (\ref{egus}) and obtain the following analogue of (\ref{e18}):
$$-DK_\Sig\ge2+(r-3)+(DE-\Card\, \bn^{-1}(E))
+\sum_Q(\ord Q-1)+\ord(C,z)$$
$$=-DK_\Sig+\left(\sum_Q(\ord Q-1)\right)+\left(\frac{DE}{2}-
\Card\, \bn^{-1}(E)
\right)
+(\ord(C,z)-1)-1\;,$$ which yields that $\bn$ is an immersion.

Let
$\bn:\PP^1\to C$ be birational, and $r=1$. The latter equality reads
$$-(K_\Sig+E)D+\frac{DE}{2}=2\;.$$ In view of $-(K_\Sig+E)D\ge0$, we
have the following possibilities:
\begin{equation}\begin{cases}\text{either}\quad & DE=0,-DK_\Sig=2,\\
\text{or}\quad & DE=2,
-DK_\Sig=3,\\ \text{or}\quad & DE=-DK_\Sig=4.\end{cases}\label{1i}\end{equation}
Claim (i) automatically holds in the first two cases. In the third case, one always has either (i), or (ii).

Finally, consider the case of
$\bn:\PP^1\to C$ being an $s$-multiple covering, $s\ge2$. Then $D=sD'$
and, according to our assumptions (\ref{ed2g0}),
$D'^2> 0$,
so that, in particular,
$D'$ is different from $E$ and $E_0$. We have
$$r-1=-DK_\Sig-2-\frac{DE}{2}=DE_0+\frac{DE}{2}-2\le D'E_0+\frac{D'E}{2}-1=
-D'K_\Sig-1\;,$$ which
implies
$$-\frac{s-2}{2}D'K_\Sig-\frac{s}{2}(K_\Sig+E)D'\le 1\;..$$ Hence,
\begin{enumerate}\item[(a)] either $s=2$, $-(K_\Sig+E)D'=1$, \item[(b)] or
$s=2$, $-(K_\Sig+E)D'=0$, \item[(c)] or $s=4$, $-D'K_\Sig=D'E=1$.
\end{enumerate}
The inequality $(D')^2>0$ excludes
case (c). In case (b) we have $D'E_0=0$ and $D'E\le 2$ (since $s=2$
and the intersection points with
$E$ must be ramification points of the covering), that is, $D'K_\Sig\ge-2$,
which together with $(D')^2>0$ leads to $p_a(D')\ge1$. Hence, a rational curve $C\in|D'|$
has singular points with the total $\delta$-invariant $((D')^2+D'K_\Sig)/2+1$,
which yields that,
in a deformation into a family of rational curves
$[\bn_t: \PP^1\to\Sig]$, $t\ne 0$,
belonging to ${\mathcal V}^{im}(\Sig,E,D)$,
we will get curves with
the total $\delta$-invariant
at least $4[((D')^2+D'K_\Sig)/2+1]$. However, the genus bound yields
$$4\left(\frac{(D')^2+D'K_\Sig}{2}+1\right)\le\frac{(2D')^2+2D'K_\Sig}{2}+1,
$$
which implies $D'K_\Sig\le-3$
contrary to $D'K_\Sig\ge-2$ pointed above.

In case (a), the same ramification argument gives
$D'E\le2$, and then the same genus argument rules out the options
$D'E=0$ and $D'E=1$.
At last,
if $D'E=2$, the curve $C$ appears to be the image of a generic
element of the family ${\mathcal V}^0(\Sig,E,D')$ with
$r(\Sig,D',0)=2$, and hence claim (iii) follows by Lemma \ref{lp1-3}.
\proofend

\begin{lemma}\label{lp2-2}
Assume
$D\in\Pic_+(\Sig,E)$. Consider an irreducible
equisingular family
${\mathcal V}\subset\overline{\mathcal V}
(\Sig,E,D)
\setminus{\mathcal V}
(\Sig,E,D)$
with $\idim{\mathcal V}=r-1$. If
 $\bn(\widehat C)\not\supset E$ for
a generic element
$[\bn:\widehat C\to\Sig]\in{\mathcal V}$, then
the following holds:
\begin{enumerate}\item[(i)] either, $\widehat C=\widehat
C_1\cup\widehat C_2$, where $\widehat C_1\simeq\widehat C_2\simeq
\PP^1$, $|\widehat C_1\cap\widehat C_2|=1$, and, for $i=1,2$,
$[\bn:\widehat C_i\to\Sig]$ represents a generic element in a family
$V_i\subset{\mathcal V}^{l_i}
(\Sig,E,D_i)$ with $\idim
V_i=r(\Sig,D_i,l_i)$, $l_i=D_iE/2$, with some
$D_1,D_2\in\Pic_+(\Sig,E)$ such that $D_1+D_2=D$, $D_1D_2>0$, the
numbers $D_1E$, $D_2E$ are non-negative and even, and, in addition,
\begin{enumerate} \item[(a)] either $\bn:{\widehat C_i}\to C_i=\bn(\widehat
C_i)$, $i=1,2$, are
immersions, and the curves $C_1,C_2$ intersect transversally,
\item[(b)] or
$\bn:\widehat C_1\to C_1$ is an immersion, $\bn:\widehat C_2\to C_2$
is a double covering, where $C_2$ is as in
Lemma \ref{lp1-4}, and the curves $C_1,C_2$ intersect
transversally,
\item[(c)] or $\bn:\widehat C_2\to C_2$ is a
double covering of $C_2=E_0$ with ramification at $E\cap E_0$, and $\bn:\widehat
C_1\to\Sig$ is birational onto its image, which is disjoint from $E\cap E_0$,
and $(\bn\big|_{\widehat C_1})^*(E_0)$ is supported at $DE_0-1$ distinct points,
\end{enumerate} furthermore, in each
of the
cases (a), (b), and (c), one has $C_1\cap C_2\cap E=\emptyset$;
\item[(ii)] or $\widehat C=\widehat
C_1\cup\widehat C_2$, where $\widehat C_1\simeq\widehat C_2\simeq
\PP^1$, $|\widehat C_1\cap\widehat C_2|=1$, and $[\bn:\widehat
C_i\to\Sig]$ is a generic element in a
family
${\mathcal V}_i\subset{\mathcal V}^{l_i}(\Sig,E,D_i)$ such that $\idim V_i=r(\Sig,D_i,l_i)$, $D_iE=2l_i+1$, $D_iE_0\ge0$
for $i=1,2$, and
$D_1+D_2=D$; furthermore, $C_1=\bn(\widehat C_1)$ and
$C_2=\bn(\widehat C_2)$ intersect along $E$ at one point $z$, where
each of these curves has a smooth local branch transversal to $E$;
\item[(iii)] or $\widehat C=\widehat
C_1\cup\widehat C_2\cup\widehat C_3$, where \begin{itemize}\item
$\widehat C_1\simeq\widehat C_2\simeq\widehat C_3\simeq \PP^1$,
$\widehat C_1,\widehat C_2$ are disjoint and each of them is joined
by one node with $\widehat C_3$,\item $\bn:\widehat C_i\to\Sig$,
$i=1,2$, are immersions and represent distinct generic elements in
some ${\mathcal V}(\Sig,E,D_i)$, where
$D_iE_0\ge1$, $i=1,2$,
\item $\bn:\widehat C_3\to\Sig$ is a double covering of $E_0$ ramified
at $E\cap E_0$.\end{itemize}
\end{enumerate}
\end{lemma}

{\bf Proof.} Under the hypotheses of the lemma, the source curve $\widehat C$ is reducible.
In every
deformation of $[\bn:\widehat C\to\Sig]$ into a generic
element $[\widetilde\bn:\PP^1\to\Sig]\in {\mathcal V}
(\Sig,E,D)$, the local branches of
$\bn:\widehat C\to\Sig$ centered along $E$ and crossing $E$ with
odd multiplicity (briefly {\it odd branches}) must
glue up
pairwise,
and
hence they are among the
nodal branches
of $\widehat C$. Write $\widehat
C$ in the form $\widehat
C=\widehat C^{(1)}\cup\widehat C^{(2)}$
so that
$\bn(\widehat
C^{(1)})\not\supset E_0$
while that each irreducible
component of $\widehat C^{(2)}$
is mapped onto $E_0$.

Now form a graph $\Gamma$ whose vertices
represent
the irreducible components of $\widehat C$ and arcs
represent the nodes
that
are intersections
of odd branches. This graph is a forest. Its isolated points
correspond to the irreducible components of $\widehat C$ which have no odd
branches; denote by $\widehat C^{(i)}_{even}$, $i=1,2$, the union of such
components inside $\widehat C^{(i)}$, and by $\widehat
C^{(i)}_{odd}$ the union of the remaining components.

The irreducible
components $\widehat C'$ of $\widehat C^{(1)}_{odd}$ such that
$E_0
\cdot \bn_*\widehat C'=0$
and all
the
local branches
of $\bn\vert_{\widehat C'}$
centered on $E$ are
odd, will be called irregular, while the other ones - regular.
Observe that only regular components of $\widehat C^{(1)}_{odd}$ can
be joined with
an irreducible component of $\widehat C^{(2)}_{odd}$ by
an arc
in
$\Gamma$.

For
each
nontrivial connected component of $\Gamma$, we
choose some "initial" vertex and orient all its arcs in the outward
direction with respect to the initial vertex.

We restrict our choice of
the
initial vertex
by the following conditions:
if there are regular vertices,
we choose one of them;
and if there is no regular vertex, but a
vertex from $\widehat C^{(2)}$, we choose one of these latter ones.

Denote
by $\widehat C_{odd}^{(1)
,irr}$ the union of those irregular
components of $\widehat C_{odd}^{(1)}$, which are terminal vertices in the
oriented graph
$\Gamma$, and by $\widehat
C_{odd}^{(1),reg}$ the union of the other components of $\widehat
C_{odd}^{(1)}$.
Denote by
$c_{even}^{(1)}$ and $c_{odd}^{(1),reg}$
the number of components on
$\widehat C_{even}^{(1)}$ and $\widehat C_{odd}^{(1),reg}$,
respectively.

Now we estimate $\idim{\mathcal V}$ from above by performing the following
procedure:
\begin{itemize}\item extend the partial order on the vertices of $\Gamma$
up to a
linear one and replace ${\mathcal V}$ by its finite cover ${\mathcal V}'$
parameterizing elements
$[\bn:\widehat C\to\Sig]$ with ordered components of $\widehat C$,
\item project ${\mathcal V}'$ to the family of elements
$[\bn:\widehat C_1\to\Sig]$, where $\widehat C_1$ is the first component of
$\widehat C$,
then take a generic fiber of the projection and project it to the family of
elements
$[\bn:\widehat C_2\to\Sig]$, where $\widehat C_2$ is the second component of
$\widehat C$, and so on;
\item notice that $[\bn:\widehat C_k\to\Sig]$, $k\ge1$, varies in the family
restricted by the condition that
all odd branches belonging to the nodes that join $\widehat C_k$ with the
preceding components of $\widehat C$
have fixed centers on $E$;
\item then apply Lemma \ref{lp1-5}.
\end{itemize}
Hence,
\begin{equation}\idim{\mathcal V}\le-D^{(1)}K_\Sig-c_{even}^{(1)}-
c_{odd}^{(1),reg}-\frac{D_{even}^{(1)}E}
{2}-l_{odd}^{(1),reg}-
\frac{D_{odd}^{(1)}E-2l_{odd}^{(1),reg}+a}{2}\;,\label{e-expr}\end{equation}
where $D^{(1)}=D_{even}^{(1)}+D_{odd}^{(1)}$, $D_{even}^{(1)}=\bn_*\widehat
C_{even}^{(1)}$,
$D_{odd}^{(1)}=\bn_*\widehat C_{odd}^{(1)}$,
$2l_{odd}^{(1),reg}$ is the total intersection multiplicity of even
local branches of $\bn:\widehat C_{odd}^{(1),reg}\to\Sig$ with $E$,
and
$a$ is the number of arcs joining in $\Gamma$ components of
$\widehat C_{odd}^{(1)}$ with $\widehat C^{(2)}$.
Observe that in case of equality, for a generic $[\bn:\widehat C\to\Sig]\in
{\mathcal V}$,
the element $[\bn:\widehat C_1\to\Sig]$ satisfies the conclusions of Lemmas
\ref{lp1-2}(2),
\ref{lp1-3}, or \ref{lp1-4}; in particular, $(\bn\big|_{\widehat C_1})^*(E)$ is the sum of
distinct simple and double points.

On the other hand,
the value (\ref{e-expr}) must
not be less than
$$r-1=-DK_\Sig-2-\frac{DE}{2}\;.
$$
Plugging
$D=D_{even}^{(1)}+D_{odd}^{(1)}+D^{(2)}$, where
$D^{(2)}=\bn_*\widehat C^{(2)}$, into the latter expression and
comparing it with (\ref{e-expr}), we arrive
to
\begin{equation}c_{even}^{(1)}+c_{odd}^{(1),reg}+\frac{a}{2}\le2\;.
\label{en1}\end{equation} We
now analyze possible values of the parameters in this relation:
\begin{itemize}\item If $a=4$, then $\widehat
C_{odd}^{(1),reg}=\emptyset$, which is impossible, since
each of the
irreducible components of $\widehat C^{(2)}$ has an even number of
odd branches, and thus neither of them
can be terminal in $\Gamma$.
\item If $a=2$, then
$c_{odd}^{(1),reg}=1$, which again
is not possible, since otherwise the only component of $\widehat
C_{odd}^{(1),reg}$
is joined by two nodes with a component of
$\widehat C^{(2)}$, thus forming a curve of positive arithmetic
genus. In particular, we conclude that $\widehat
C_{odd}^{(2)}=\emptyset$ and
$a=0$.
\item If $c_{even}^{(1)}=2$,
then $\widehat C_{odd}^{(1),reg}=\emptyset$, which also yields
$\widehat C_{odd}^{(1),irr}=\emptyset$. This implies the properties
enumerated at the beginning
of item (i), and, by Lemmas \ref{lp1-2} - \ref{lp1-3}, under
assumption that
$\widehat C^{(2)}=\emptyset$ leads to assertions (i-a,b).
Assume that $\widehat C_{even}^{(2)}\ne\emptyset$. Then,
$\bn_*\widehat
C_{even}^{(2)}=2sE_0$, $s\ge1$. By Lemma \ref{lp1-3} the
components of $\widehat C_{even}^{(1)}$
have in total
$(D-2sE_0)E_0=DE_0+2s$ local branches centered on $E_0$. So, at most two
of these branches glue up with $\widehat C_{even}^{(2)}$, and hence
at least $DE_0+2s-2$ branches centered along $E_0$ persist in a
deformation into a generic element of ${\mathcal V}
(\Sig,E,D)$,
whence $s=1$, and we fit the
assertions of
item (iii). Notice
that the two components $D_1, D_2$ of $\widehat C^{(1)}_{even}$ are mapped to
distinct curves. Indeed, the dimension count $D_iE_0+{D_iE}=1$ leaves the only
option
of $D_1=D_2=E+E_0=-K_\Sig$ for possibly coinciding images $C_1,C_2$ of the
components of $\widehat C^{(1)}_{even}$ ({\it cf.} the proof of Lemma \ref{lp1-2}
and note that $D_i^2>0$ since $D^2>0$):
however in such a case, the curve
$C_1+C_2+2E_0$ would deform into a rational curve with four nodes in
a neighborhood of the node of $C_1=C_2$,
but this contradicts the
arithmetic genus formula $p_a(-2K_\Sig+2E_0)=3$.
\item If
$c_{even}^{(1)}=
c_{odd}^{(1),reg}=1$,
then $\widehat C^{(2)}_{odd}=\emptyset$,
$\widehat
C_{odd}^{(1),irr}\ne\emptyset$, and each component of $\widehat
C_{odd}^{(1),irr}$ is joined to $\widehat C^{(1),reg}_{odd}$ by a node and has
a unique odd branch.
Each component $C'$ of $C^{(1),irr}_{odd}=\bn_*\widehat C^{(1),irr}_{odd}$
has $C'E_0=0$ and, because of equality in (\ref{e-expr}),
$C'E=1$.
It is easy to see
that then $C'$ is a $(-1)$-curve intersecting $E$ at one point.

If the curve
$C^{(1),reg}_{odd}=\bn_*\widehat C^{(1),reg}_{odd}$ either positively
intersects $E_0$, or
has local branches intersecting $E$ with even multiplicity,
Lemma \ref{lp1-5}(i) applies, and we get the following upper bound to
$\idim{\mathcal V}$:
$$\left(-D^{(1)}_{even}K_\Sig-\frac{D^{(1)}_{even}E}{2}-1\right)+
\left(-D^{(1),reg}_{odd}K_\Sig-l'-(D^{(1),reg}_{odd}E-2l')-1\right)\ ,$$
where $l'$ is defined by
$\big[\bn\vert_{{C^{(1),reg}_{odd}}}:
C^{(1),reg}_{odd}\to\Sig]
\in{\mathcal V}^{l'}(\Sig,E,
D^{(1),reg}_{odd})$. Taking into account that $D=D^{(1)}_{even}+D^{(1),reg}_{odd}+
D^{(1),irr}_{odd}+2sE_0$, $s\ge0$, we
obtain the estimate
$$r-1=\idim{\mathcal V}\le \left(-DK_\Sig-\frac{DE}{2}-2\right)-\left(-D^{(1),irr}_{odd}K_\Sig-\frac{D^{(1),
irr}_{odd}E}{2}\right)$$
$$-\left(\frac{D^{(1),reg}_{odd}E}{2}-l'\right)-s<-DK_\Sig-\frac{DE}{2}-2=r-1\ ,$$
which
is a contradiction. Thus, $C^{(1),reg}_{odd}E_0=0$ and $C^{(1),reg}_{odd}E=1$,
which implies, similar to the above, that
$C^{(1),reg}_{odd}=C'$
and, hence, $\bn_*(\widehat C^{(1),reg}_{odd}\cup\widehat C^{(1),irr}_{odd})
=sC'$, $s\ge2$. Note, that $\bn(C^{(1)}_{even})\ne C'$, since otherwise we
would get a non-primitive
divisor $D$ with $D^2\le0$ contrary to the hypotheses of Lemma. By Lemma
\ref{lp1-3}, we can assume that the curve
$\bn\big|_{\widehat
C^{(1)}_{even}}$ intersects $C'$ transversally. So, the relative to boundary
cycle realized in a small tubular neighborhood of $C'$
by the
restriction of $\bn$ to a neighborhood of $\widehat C^{(1),reg}_{odd}\cup
\widehat C^{(1),irr}_{odd}$ in $\widehat C$
has a negative total intersection number with $C'$
(indeed, this intersection number is equal to $-s+1\le -1$), which contradicts
the existence of a deformation family $\bn_t$.
\item If
$c_{even}^{(1)}=0$, $c_{odd}^{(1),reg}=1$, then similarly we get $\widehat
C^{(2)}_{odd}=\emptyset$ and
$\widehat C^{(1), irr}_{odd}\ne\emptyset$ so that each component $C'$ of
$\bn_*\widehat C^{(1),irr}_{odd}$ satisfies $C'E_0=0$ and has a unique
branch centered on $E$,
and this branch intersects $E$ with odd multiplicity. The preceding item
argument
leads to the bound
$$r-1=\idim{\mathcal V}\le(r-1)+1-\left(-D^{(1),irr}_{odd}K_\Sig-
\frac{D^{(1),irr}_{odd}E}{2}\right)-
\left(\frac{D^{(1),reg}_{odd}E}{2}-l'\right)-s$$
that leaves the only possibility $D^{(1),reg}_{odd}E=2l'+1$, $\widehat C^{(1),
irr}_{odd}$ is irreducible, and $s=0$, which
fits the assertions of item (ii), or leaves the only case when $c^{(1),reg}_{odd}=
c^{(1),irr}_{odd}=1$, and both $\widehat C^{(1),reg}_{odd}$ and $\widehat C^{(1),
irr}_{odd}$ are mapped onto the same $(-1)$-curve
intersecting $E$ at one point. However, the latter option corresponds to the
forbidden, by the hypotheses of Lemma, case of
a non-primitive $D$ with $D^2\le0$.
\item If
$c_{odd}^{(1),reg}=2$, then $\widehat C_{even}^{(1)}=\emptyset$
and $\widehat C^{(2)}_{odd}=\emptyset$.
If $\widehat C_{odd}^{(1),irr}=\emptyset$, then the preceding item
arguments bring us to
conclusions of item (ii). If $\widehat
C_{odd}^{(1),irr}\ne\emptyset$, then we have $\widehat C=\widehat
C_1\cup \widehat C_2\cup\widehat C_3$,
and, once more by the similar arguments and using additionally the property $T(1)$,
we come to the upper bound
$r-1=\idim{\mathcal V}\le \sum(-KC_i - l_i -1) - 2$, which contradicts to
$r-1=-KD -l -2$, since $l=\sum l_i +2$.
\item If $c_{even}^{(1)}=1$, $c_{odd}^{(1),reg}=0$, then $\widehat
C_{odd}^{(1)}=\emptyset$. Since $\widehat C$ is reducible, we derive
that $\widehat C_{even}^{(2)}\ne\emptyset$, $\bn_*\widehat
C_{even}^{(2)}=2sE_0$, $s\ge1$, and $C^{(1)}E_0=DE_0+2s$, where
$C^{(1)}=\bn_*\widehat C_{even}^{(1)}$. In particular,
\begin{equation}r=r(\Sig,D,l)=r(\Sig,C^{(1)},C^{(1)}E/2)=\frac{C^{(1)}E}{2}+
DE_0+2s-1\ .\label{r12}\end{equation}
Show, first, that
$[\bn:\widehat C^{(1)}\to\Sig]$ cannot have $C^{(1)}E_0=DE_0+2s$ smooth transversal
branches along $E_0$. Indeed, otherwise, in a deformation of
\mbox{$[\bn:\widehat C\to\Sig]$} into a generic element of ${\mathcal V}
(\Sig,E,D)$, at least $C^{(1)}E_0-s=DE_0+s>DE_0$ local branches
of $C^{(1)}$ centered on $E_0$ would persist. Hence,
$[\bn:\widehat C^{(1)}\to\Sig]$ is generic in an $(r-1)$-dimensional subfamily
of ${\mathcal
V}^{C^{(1)}E/2}(\Sig,E,C^{(1)})$ subject to condition $\Card(\bn\big|_{\widehat
C^{(1)}})^{-1}(E_0)\le
C^{(1)}E_0-s$.

Let $r-1>0$. Then by Lemma \ref{lp2-1}(iii) and
the smoothing argument, \mbox{$\bn:\widehat C^{(1)}\to C^{(1)}$} cannot be a multiple
covering. Thus, $\bn:\widehat C^{(1)}\to C^{(1)}$ is birational.
Inequality (\ref{egus}) immediately yields that
$\bn:C^{(1)}\to\Sig$ has $C^{(1)}E_0-1$ local branches
centered on $E_0$, and, furthermore, $s=1$.
That is, we
fit the assertions of item (i-c).

If $r=1$, then (\ref{r12}) implies that $C^{(1)}E=DE_0=0$ and $s=1$. Hence, we again
fit the assertions of item
(i-c),
provided we show that $\bn:\widehat C^{(1)}\to C^{(1)}$ is not a double cover,
{\it i.e.} $C^{(1)}\ne 2C'$. Indeed, if it were so, from $C'E=0$ and $C'E_0=1$,
we could derive that
$C'$ is a $(-1)$-curve, disjoint from $E$ and crossing $E_0$ at one point, which
would yield
that $D=2D_1$, $D_1^2=0$ contrary to the hypotheses of Lemma.\proofend
\end{itemize}

\begin{lemma}\label{lp2-3}
Assume that $D\in\Pic_+(\Sig,E)$.
Consider an irreducible
equisingular family
${\mathcal V}\subset\overline{\mathcal V}
(\Sig,E,D)
\setminus{\mathcal V}
(\Sig,E,D)$
with $\idim{\mathcal V}=r-1$.
If
 $\bn(\widehat C)\supset E$ for
a generic element
$[\bn:\widehat C\to\Sig]\in{\mathcal V}$, then
the following holds:
\begin{enumerate} \item[(i)] either $\bn_*(\widehat C)=2E+k_0E_0$ with $k_0=2$ or $4$;
\item[(ii)] or $\widehat C=\widehat
C_1\cup\widehat C_2$, where $\widehat C_1\simeq\widehat C_2\simeq
\PP^1$, $|\widehat C_1\cap\widehat C_2|=1$, $\bn\big|_{\widehat
C_1}$ is an immersion onto $C_1\subset\Sig$, $\bn\big|_{\widehat
C_2}$ is an isomorphism onto $C_2=E$, and \mbox{$[\bn:\widehat
C_1\to\Sig]$} is a generic element in a component of ${\mathcal
V}(\Sig,E,D-E)$ of intersection dimension $r-1$;
\item[(iii)] or $\widehat C=\widehat C_1\cup\widehat C_2\cup\widehat C_3$, where
$\widehat C_1\simeq\widehat C_2\simeq\widehat C_3\simeq \PP^1$,
$\widehat C_3$ is isomorphically taken onto $E$, and the following
holds: \begin{enumerate}\item[(a)] either $\widehat C_1\cap\widehat
C_2=\emptyset$, $|\widehat C_1\cap\widehat C_3|=|\widehat
C_2\cap\widehat C_3|=1$, and $\bn\big|_{\widehat C_i}$, $i=1,2$, are
immersions, and the following holds: $C_1\cap C_2\cap E=\emptyset$,
$[\bn:\widehat C_i\to\Sig]\in{\mathcal V}^{l_i}(\Sig,E,D_i)$,
$i=1,2$, are generic elements in the corresponding families for some
$D_1,D_2\in\Pic_+(\Sig,E)$ such that $D_1+D_2=D-E$, $D_iE=2l_i+1$,
$i=1,2$; \item[(b)] or $|\widehat C_1\cap\widehat C_2|=|\widehat
C_2\cap\widehat C_3|=1$, $\widehat C_1\cap\widehat C_3=\emptyset$,
$[\bn:\widehat C_1\to\Sig]\in{\mathcal V}(\Sig,E,D - E - 2E_0)$ is a
generic element,
$\bn\big|_{\widehat C_1}$ is an immersion, and
$\bn\big|_{\widehat C_2}$ is the double covering of $E_0$
ramified at $E\cap E_0$.
\end{enumerate}
\end{enumerate}
\end{lemma}

{\bf Proof.}
Write
$\widehat C=\widehat C^{(1)}\cup\widehat C^{(2)}\cup\widehat
C^{(3)}$, where $\bn(\widehat C^{(1)})$ does not contain neither
$E$ nor $E_0$, $\bn_*\widehat C^{(2)}=k_0 E_0$ with $k_0\ge 0$,
and $\bn_*\widehat C^{(3)}=kE$
with $k\ge1$. We follow the lines of the proof of Lemma \ref{lp2-2},
and
keep the notations and definitions introduced there, except for the following
modification:
by $\widehat C^{(1),irr}_{odd}$ we denote the union of those irregular components of $\widehat C^{(1)}$,
which are terminal in $\Gamma$ and disjoint from $\widehat C^{(3)}$.

Assume first
that $\widehat C^{(1)}_{odd}=\emptyset$.
Then
$$r-1=-DK_\Sig-2-\frac{DE}{2}=\idim
\mathcal{V}
\le-D^{(1)}K_\Sig-
c^{(1)}_{even}-\frac{D^{(1)}E}{2}\;.$$
Substituting
$D=D^{(1)}+kE+k_0E_0$, we obtain
$$k+
c^{(1)}_{even}\le2\;.
$$
\begin{itemize}\item If $k=2$, then $\widehat C^{(1)}=\widehat C^{(1)}_{even}=\emptyset$. In
view of $DE\ge0$, $DE_0\ge0$, we get $2\le k_0\le4$.
The case $k_0=3$
is excluded by Lemma \ref{lp1-6}.
The left out options $k_0=2$ or $4$ fit the
the assertions of item (i).

\item
If $k=1$, then  $\widehat C^{(3)}\simeq\PP^1$ and
 $\bn$ maps it isomorphically onto $E$.
If, in addition,
$C^{(1)}=\emptyset$, then $D=E+k_0E_0$, in which case
$r=r(\Sig,D,l)=0$ contrary to
our assumption that $r>0$. Thus,
$\widehat C^{(1)}\simeq\PP^1$ and we come to the properties enumerated in item (ii), if
$\widehat C^{(2)}=\emptyset$,  and to those of item (iii-b), otherwise.
\end{itemize}

Now
assume that $C^{(1)}_{odd}\ne\emptyset$. Observe that
$\widehat C^{(1),reg}_{odd}\ne\emptyset$,
while as always $\widehat
C_{odd}^{(1),irr}\cap\widehat C^{(3)}=\emptyset$
and
$|\widehat
C_{odd}^{(1),reg}\cap\widehat C^{(3)}|\le
c_{odd}^{(1),reg}+c^{(3)}-1$, where $c^{(3)}\le k$ is the number of
components of $\widehat C^{(3)}$. This implies the following upper
bound for $\idim{\mathcal V}$ ({\it cf.} the bound (\ref{e-expr}) in
the proof of Lemma \ref{lp2-2}):
$$-D^{(1)}K_\Sig-c^{(1)}_{even}-c^{(1),reg}_{odd}-
\frac{D^{(1)}_{even}E}{2}-l^{(1),reg}_{odd}-\frac{D_{odd}^{(1)}E-2l_{odd}^{(1),reg}+
a-c_{odd}^{(1),reg}-c^{(3)}+1}{2},$$ which together with $\idim{\mathcal V}=r-1$ yields
\begin{equation}
c^{(1)}_{even}+\frac{c^{(1),reg}_{odd}}{2}+\frac{a}{2}+\left(k-\frac{c^{(3)}}{2}\right)\le\frac{3}{2}\
.\label{elast2}\end{equation}
In veiw of
$k\ge c^{(3)}\ge 1$, it gives
$c^{(1)}_{even}=0$, $1\le
c^{(1),reg}_{odd}\le2$ and $1\le k\le 2$.
\begin{itemize}
\item If $k=2$ then $c^{(3)}=2$, $c^{(1),reg}_{odd}=1$, $a=0$, and
$\widehat C^{(1),reg}_{odd}\simeq\PP^1$.
Note that
according to Lemma\ref{lp1-5}(ii)
the map
$\bn:\widehat C^{(1),irr}_{odd}$
constrained by the position of its intersection points with $E$
is rigid. Hence,
according to Lemma \ref{lp1-5}(i)
$$-C^{(1),reg}_{odd}K_\Sig-1-C^{(1),reg}_{odd}E+
\Card(\bn\big|_{\widehat C^{(1),reg}_{odd}})^{-1}(E)
\ge-DK_\Sig-\frac{DE}{2}-2\ .$$ Substituting there $D=C^{(1),reg}_{odd}+C^{(1),irr}_{odd}+k_0E_0+2E$
and taking into account that $-C^{(1),irr}_{odd}K_\Sig=C^{(1),irr}_{odd}E$, we arrive to
$$
2\,\Card(\bn\big|_{\widehat C^{(1),reg}_{odd}})^{-1}(E)
\ge C^{(1),irr}_{odd}E+C^{(1),reg}_{odd}E+2\ .
$$
Furthermore, selecting in $(\bn\big|_{\widehat C^{(1),reg}_{odd}})^*(E)$ only components of odd multiplicity, we obtain
\begin{equation}
\Card(\bn\big|_{\widehat C^{(1),reg}_{odd}})^{-1}(E)
\ge C^{(1),irr}_{odd}E+2
\label{eodd}\end{equation} with the equality only if
each even component of $(\bn\big|_{\widehat C^{(1),reg}_{odd}})^*(E)$ has multiplicity $2$ and each odd
component has multiplicity $1$. On the other hand,
$$
\Card(\bn\big|_{\widehat C^{(1),reg}_{odd}})^{-1}(E)
\le\Card(\bn\big|_{\widehat C^{(1),irr}_{odd}})^{-1}(E)+2\ ,$$
which yields the equality in (\ref{eodd}) and the fact that $[\bn:
\widehat C^{(1)}_{odd}\to\Sig]$ represents a generic element in the family
${\mathcal V}^{l_1}(\Sig,E,C^{(1)})$,
where $l_1$ is the number of even components in
$(\bn\big|_{\widehat C^{(1),reg}_{odd}})^*(E)$.
In particular, $\bn:\widehat C^{(1),reg}_{odd}\to\Sig$ avoids $E\cap E_0$, and
the divisor $\bn\big|_{\widehat C^{(1),reg}_{odd}}^*(E_0)$
is supported at $C^{(1),reg}_{odd}E_0$ points.

Let us show that $\widehat C^{(2)}=\emptyset$. If it were not so, then
$D=C^{(1)}+2E+k_0E_0$, $k_0>0$. One can see that $\widehat
C^{(2)}=\widehat C^{(2)}_{even}$, and each component of $\widehat
C^{(2)}$ would be joined by a node either with $\widehat C^{(1)}$ or with
$\widehat C^{(3)}$. Then at least
$$C^{(1),reg}_{odd}E_0+4-\frac{k_0}{2}>DE_0=C^{(1),reg}_{odd}E_0+4-k_0$$ local branches of
$\bn:\widehat C^{(1)}\cup\widehat C^{(3)}\to\Sig$ centered on $E_0$
would persist in a deformation into a generic element of ${\mathcal V}
(\Sig,E,D)$, which is a contradiction.

Thus, $D=C^{(1)}+2E$, but then $DE=C^{(1)}E-4$ which contradicts the
fact that two local branches of $\bn:\widehat C^{(1)}\to\Sig$ centered on $E$ glue up with components of
$C^{(3)}$ in a deformation into a generic
element of ${\mathcal V}
(\Sig,E,D)$, and these branches intersect $E$ with multiplicity $1$.
\item It remains to
treat
the
case $k=1$. First, we exclude the case
$c^{(1),reg}_{odd}=1$. Indeed, since $DE$,
as well as the number of nodes joining in $\Gamma$ the components of $C^{(2)}$ with the only component
of $C^{(3)}$,
are
even, the number of odd branches of $\bn:\widehat C^{(1),reg}_{odd}\to\Sig$ on $E$ that are not
nodes joining
$\widehat C^{(1),reg}_{odd}$ with $\widehat C^{(1),irr}_{odd}$
is
also
even, and if such branches do exist, in the deformation, they glue up with $\widehat C^{(3)}\simeq E$
into a curve of a positive genus. Hence, all odd branches of $\bn:\widehat C^{(1),reg}_{odd}\to\Sig$
on $E$ are joined to odd branches of $\bn:\widehat C^{(1),irr}_{odd}\to\Sig$. In particular, it implies
$$\Card(\bn\big|_{\widehat C^{(1),reg}_{odd}})^{-1}(E)
\le C^{(1),irr}_{odd}E\ .$$

On the other hand,
writing $D=C^{(1),reg}_{odd}+C^{(1),irr}_{odd}+k_0E_0+E$, we transform,
as in the previous item,
the dimension inequality
$$-DK_\Sig-\frac{DE}{2}-2\le-C^{(1),reg}_{odd}K_\Sig-C^{(1),reg}_{odd}E+\Card(\bn\big|_{\widehat
C^{(1),reg}_{odd}})^{-1}(E)-1$$ into
\begin{equation}\Card(\bn\big|_{\widehat
C^{(1),reg}_{odd}})^{-1}(E)
\ge C^{(1),irr}_{odd}E\ ,\label{elast1}\end{equation}
and arrive to the conclusion that
all the branches of $\bn:\widehat C^{(1),reg}_{odd}\to\Sig$ on $E$
intersect $E$ with multiplicity $1$ or $2$, and
that
$[\bn:\widehat C^{(1),reg}_{odd}\to\Sig]$ is a generic member in the corresponding family
${\mathcal V}^{l_1}(\Sig,E,C^{(1),reg}_{odd})$.
Furthemore, since (\ref{elast1}) must
turn into
an equality,
it follows that:
$\widehat C^{(2)}\ne\emptyset$;
$\widehat C^{(2)}_{odd}=\emptyset$;
$\widehat C^{(1),reg}_{odd}$
is joined by one node with one component
of $\widehat C^{(2)}$ ({\it cf.} inequality (\ref{elast2}) which yields $a=1$);
and
each but one component of $\widehat C^{(2)}_{even}$ is joined by a node with $\widehat C^{(3)}$.

By
Lemma \ref{lp1-3}(iv), $\Card(\bn\big|_{\widehat C^{(1),reg}_{odd}})^{-1}(E_0\setminus E)=C^{(1),reg}_{odd}E_0$. Hence, if $\widehat C^{(2)}_{even}$ contains $s\ge 1$ components, then, in the deformation
of $[\bn:\widehat C\to\Sig]$ into a generic element of ${\mathcal V}
(\Sig,E,D)$, we encounter at least
$$C^{(1),reg}_{odd}E_0-1+\max\{0,2-(s-1)\}=(DE_0+k_0-2)-1+\max\{0,3-s\}$$
$$\ge DE_0+2s-3+\max\{0,3-s\}\ge DE_0+1$$ local branches intersecting $E_0$, which is a contradiction.

Thus, $c^{(1),reg}_{odd}=2$ and $a=0$.
This time
the dimension inequality reads
\begin{equation}
-DK_\Sig-\frac{DE}{2}-2\le-G^{(1),reg}_{odd}K_\Sig-C^{(1),reg}_{odd}E+\Card(\bn\big|_{\widehat
C^{(1),reg}_{odd}})^{-1}(E)-2\ ,
\label{elast4}\end{equation} and
then
transforms into
\begin{equation}\Card(\bn\big|_{\widehat
C^{(1),reg}_{odd}})^{-1}(E)_{odd}\ge C^{(1),irr}_{odd}E+2. \label{elast3}
\end{equation}
Arguing as above, we come to the conclusion that:
all the branches of $\bn:\widehat C^{(1),reg}_{odd}\to\Sig$ on $E$
intersect $E$ with multiplicity $1$ or $2$;
(\ref{elast3}) must be an equality; the two components
$(\widehat C^{(1),reg}_{odd})'$ and $(\widehat C^{(1),reg}_{odd})''$ of
$\widehat C^{(1),reg}_{odd}$ are disjoint, each of them is joined by one
node with $\widehat C^{(3)}$; and each component of
$C^{(1),irr}_{odd}$ intersects $E$ with multiplicity $1$. In particular, all components
of $\widehat C^{(1),irr}_{odd}$ are mapped onto $(-1)$-curves intersecting $E$ at one point.
Furthermore, $(\widehat C^{(1),reg}_{odd})'$ and $(\widehat C^{(1),reg}_{odd})''$ cannot be
mapped onto these $(-1)$-curves. Indeed, if $E'$ were such a $(-1)$-curve, the deformation of
the union of the corresponding components of $\widehat C^{(1),reg}_{odd}$ and
$\widehat C^{(1),irr}_{odd}$, mapped onto $E'$, would produce in a tubular neighborhood of $E'$ in
$\Sig$ a relative holomorphic cycle negatively intersecting $E'$. Thus, if $\widehat C^{(1),irr}_{odd}\ne\emptyset$, its components represent rigid curves,
so that,
by Lemma
\ref{lp1-5}(i), in (\ref{elast4}) we have to subtract from the right-hand side the number of
components of $\widehat C^{(1),irr}_{odd}$, which then will break the equality in (\ref{elast3}) established above. Hence $\widehat C^{(1),irr}_{odd}=\emptyset$.

At last, the emptiness of
$\widehat C^{(2)}$ can be proved by the argument used in the preceding item.
So,
in the case under consideration we fit the
the assertions of item (iii-a).\proofend
\end{itemize}

\section{Proof of Theorems \ref{t3n} and \ref{t3}}\label{s2}

{\bf Proof of Theorem \ref{t3}.}
(1)
Since we identify divisor classes
orthogonal to $E$
along the
considered family, we write $|D|_X$ to specify that the linear
system $|D|$ is considered on the surface $X$.

In the setting and notations of
\cite[Section 4.2]{IKS2},
take disjoint sections $w_i:[0,\eps)\to{\mathfrak X}$, $1\le i\le
r$, such that: $r-2m$ of them are real, targeted in $F_t^+$, the
others form $m$ pairs of complex conjugate, and $\{w_i(0)\}_{1\le
i\le r}=\bw$. Consider the limits of real rational curves in
$|D|_{{\mathfrak X}_t}$  passing through $\{w_i(t)\}_{1\le i\le r}$,
$t\ne0$. By \cite[Theorem 4.2]{Va} and \cite[Lemma 4.1]{IKS2},
the limit curves are either real rational curves
$C\in|D|_\Sig$, passing through $\bw$, or have the form
$C+\sum_{1\le i\le s}(L_i+\overline L_i)$, $s\ge1$, where
$C\in|D-sE|_\Sig$ is a real rational curve passing through $\bw$ and
crossing $E$ at $2s$ distinct imaginary points, while
$L_i,\overline L_i$,
$1\le i\le s$, are $s$ pairs of complex conjugate ruling lines of
the quadric $Z\subset{\mathfrak X}_0$
 passing through all the points of $C\cap E$. Furthermore, each
rational curve $C\in|D-sE|_\Sig$ as above gives rise to $2^s$ limit
curves: for each pair $w,\overline w\in C\cap E$ we have two choices
for a pair of conjugate ruling lines of $Q$ passing through
$w,\overline w$. Now observe that one of these pairs has a solitary
node in $Q^+$ which, in deformation, goes to $F_t^+$ and contributes
factor $(-1)$ to the weight of the corresponding rational curve
$C_t\in|D|_{{\mathfrak X}_t}$ in \cite[Formula (1)]{IKS3} for
$W_m({\mathfrak X}_t,D,F_t^+,\varphi_t+[F_t^-])$, whereas the other
pair of ruling lines has a solitary node in $Q^-$, which therefore
contributes factor $1$ to the weight of the corresponding rational
curve $C'_t\in|D|_{{\mathfrak X}_t}$.
Since the other singularities
of $C_t$ and $C'_t$ are of the same real type and are close to each
other,
the weights of these curves
in \cite[Formula (1)]{IKS3} cancel each other.
Thus, we are left out
only
with the
limit curves $C\in|D|_\Sig$,
which provide the same contribution to
$RW_m(\Sig,D,F^+,\varphi,\bw)$ and to $W_m({\mathfrak
X}_t,D,F_t^+,\varphi_t+[F_t^-])$.

(2) The second statement of Theorem \ref{t3}
is
proved in the same way as the first one.
\proofend

\vskip10pt

{\bf Proof of Theorem \ref{t3n}.}
We keep the same notation and
proceed in the same way as in the proof of Theorem
\ref{t3}.
The real rational curves in $|D|_{{\mathfrak X}_t}$
passing through $\{w_i(t)\}_{1\le i\le r}$ , $t>0$, form families
that have limits of the following kind:
either real rational curves $C\in|D|_\Sig$ passing through $\bw$,
or of the form $C+\sum_{1\le i\le s'}(L_i+\overline L_i)$,
$s'\ge1$, where $C\in|D-s'E|_\Sig$ is a real rational curve passing
through $\bw$ and crossing $E$ at $2s'$ distinct imaginary points,
while
$L_i,\overline L_i$, $1\le i\le s'$, are $s$ pairs of complex
conjugate ruling lines of
$Z$
passing through all the points of
$C\cap E$. Notice that $L_i$ and $\overline L_i$ belong to the same
ruling, and hence we must have $s'=2s$. Thus, each rational curve
$C\in|D-2sE|_\Sig$, $s\ge0$, passing through $\bw$ and
crossing $E$ at $4s$ distinct imaginary points, gives rise to $\binom{2s}{s}$ limit
curves obtained by picking a ruling of
$Z$,
choosing $s$ pairs of complex conjugate points in $C\cap E$,
and attaching a line from the
preferred ruling to each of the chosen $2s$ points, whereas to the remaining points in $C\cap E$ we attach the lines from the other ruling. Since the attached lines have in total an even intersection
with $\varphi$ ({\it i.e.}, with any representing cycle), we get
\begin{equation}W_m({\mathfrak
X}_t,D,F_t,\varphi_t)=\sum_{s\ge0}\binom{2s}{s}\widetilde
W_m(\Sig,D-2sE,F,\varphi,\bw)\,.\label{e2504a}\end{equation} where
$\widetilde W(\Sig,D-2sE,F,\varphi,\bw)$ counts real rational irreducible
curves $C\in|D-2sE|$ on $\Sig$, passing through $\bw$, disjoint from $\R E$, and equipped
with the weights $(-1)^{C_+\circ\varphi}\mu(F^+,C)$
({\it cf.} formula (\ref{e22n})).
Applying the same counting procedure to the divisor classes $D-2kE$,
$k\ge1$, we obtain
\begin{equation}W_m({\mathfrak
X}_t,D-2kE,F_t,\varphi_t)=\sum_{s\ge0}\binom{2s+2k}{s}\widetilde
W_m(\Sig,D-2(s+k)E,F,\varphi,\bw)\;.\label{e2504b}\end{equation}
Then (\ref{e2504c}) follows from (\ref{e2504a}) and (\ref{e2504b}).
\proofend

\section{Proof of Theorems \ref{t2}, \ref{t4p}, and \ref{cor1}}\label{s3}

The proof is organized as follows. First, we simultaneously establish
Theorem \ref{t2} and Theorem \ref{t4p} for uninondal DP-pairs possessing
property $T(1)$ (Sections \ref{sec1} and \ref{sec1a}), and then we prove Theorem \ref{t4p}
for uninodal DP-pairs lacking property $T(1)$
(Section \ref{sec1b}). In Section \ref{sec4.5} we prove Theorem \ref{cor1}.

\subsection{Preparation}\label{preparation}
Suppose that the
tuple $(\Sigma,E,F,D)$
satisfies the hypotheses of one of
Theorems
\ref{t2} and \ref{t4p}; in particular, we suppose that
$r=r(\Sig,D,l)>0$, where
$l=DE/2$. In view of Theorems \ref{t3n} and \ref{t3}, we can also suppose
that $DE>0$.

Furthermore, till the end of Section \ref{sec1a} we assume that the uninodal DP-pair $(\Sig,E)$
possesses property $T(1)$. Recall that, according to Lemma \ref{deg2}, if  $(\Sig,E)$ has degree $k \ge2$, then its blow up at $k - 1$
generic points is a non-tangential uninodal DP-pair of degree $1$
possessing property $T(1)$.
This allows one to treat degree $\ge 2$ tuples $(\Sigma, E,F,D)$ as degree $1$ tuples.

Throughout Section \ref{s3} we shortly write
$RW(\bw)$ for $RW_m(\Sig,E,F^+,\varphi,D,\bw)$.

Let $0 \leq i \leq r$ be an integer.
For any sequence $\bv$ of $i$ pairwise distinct points of $\Sig$,
put
\begin{eqnarray}
\overline{\mathcal V}^{\;im}_i(\Sig,E,D,\bv)&=&
\Big\{[\bn:\widehat C\to\Sig,\bp]\in\overline{\mathcal M}^*_{0,i}(\Sig,D)\ :\nonumber\\
& &[\bn:\widehat C\to\Sig]\in\overline{\mathcal V}^{\;im}(\Sig,E,D),\ \bn(\bp)=\bv\Big\}. \nonumber
\end{eqnarray}

Let $\gamma: I = [0,1]\to{\mathcal P}_{r,m}(\Sig,F^+)$ be
a smooth simple
real-analytic path that connects two generic points of ${\mathcal P}_{r,m}(\Sig,F^+)$.
The path $\gamma$ is said to be {\it qualified} if all sets $\overline{\mathcal V}^{\;im}_r(\Sig,E,D,\gamma(t))$,
$t\in I$,
are finite, and there exists a finite subset $I_1\subset I$ such that
\begin{itemize}
\item if $t\in I\setminus I_1$, then
$\pi_r(\overline{\mathcal V}^{\;im}_r(\Sig,E,D,\gamma(t)))\subset{\mathcal V}^{im}(\Sig,E,D)$;
\item
if
$t_0 \in I_1$, then
each $\xi = [\bn:\hat C\to\Sig, \widehat\bw] \in \overline{\mathcal V}_r^{\;im}(\Sig,E,D,\gamma(t_0))$
with
\mbox{$\pi_r(\xi) \not\in \mathcal V^{im}(\Sig,E,D)$}
is generic in some $(r - 1)$-dimensional stratum
listed in Lemmas \ref{lp2-1} - \ref{lp2-3}, $\widehat\bw$ is a generic $r$-tuple
on $\hat C$, and for each $i = 1$, $\ldots$, $r$, the path $t \mapsto w_i(t)$
is not tangent to $C=\bn_* \hat C$ at $w_i(t_0)$.
\end{itemize} The existence of qualified paths follows from Lemmas \ref{lem-fin}, \ref{lp2-1} - \ref{lp2-3}.

To prove Theorems
\ref{t2} and \ref{t4p}, we choose appropriate qualified paths
and apply the following {\it localization principle}.

Let $\gamma:I \to {\mathcal P}_{r,m}(\Sig,F^+)$ be a qualified path
connecting two generic points of ${\mathcal P}_{r,m}(\Sig,F^+)$.
Denote by $\Gamma(\gamma)$ (respectively, $\Gamma_1(\gamma)$) the subset of
$\bigcup_{t \in I}\overline{\mathcal V}^{\;im}_r(\Sig,E,D,\gamma(t))$
(respectively, $\bigcup_{t \in I\setminus I_1}\overline{\mathcal V}^{\;im}_r(\Sig,E,D,\gamma(t))$)
formed by the real elements.
Let $\tau: \Gamma(\gamma) \to I$ be the tautological projection.
By definition, the function $RW: I\setminus I_1 \to \Z$,
$t \mapsto RW(\gamma(t))$,
is the direct image $\tau_*\mu_{\Gamma_1}$ of the function $\mu_{\Gamma_1}: \Gamma_1 \to \Z$
defined by $\xi \mapsto \mu(F^+, \varphi, \xi)$.
Since $\tau$ is a proper map and $\Gamma_1(\gamma) \subset \Gamma(\gamma)$
is a dense open subset,
to prove that the function $RW$ is constant, it is sufficient to check
that for any $\xi \in \Gamma(\gamma)$ the direct image $\tau_*\mu_{\Gamma_1(\gamma), \xi}$
of the restriction $\mu_{\Gamma_1(\gamma), \xi}$ of $\mu_{\Gamma_1}$ to $(\Gamma(\gamma), \xi) \cap \Gamma_1(\gamma)$
is a constant function, where $(\Gamma(\gamma), \xi)$ is the germ of $\Gamma(\gamma)$ at $\xi$.

\begin{lemma}\label{les-net1}
Let $\gamma:I \to {\mathcal P}_{r,m}(\Sig,F^+)$ be a qualified path
connecting two generic points of ${\mathcal P}_{r,m}(\Sig,F^+)$.
Then, no element of $\Gamma(\gamma)$
admits a description given in
Lemma \ref{lp2-2}(i-c), Lemma \ref{lp2-2}(iii), Lemma \ref{lp2-3}(i), or
Lemma \ref{lp2-3}(iii-a).
\end{lemma}

{\bf Proof.} Each of the elements described in  Lemma \ref{lp2-2}(i-c), Lemma \ref{lp2-2}(iii),
and Lemma \ref{lp2-3}(i) involves the curve $E_0$, which is impossible in the case of a blow up
of a given degree $\ge 2$ tuple at generic points. Hence, these cases are
not relevant for Theorem \ref{t2}.

Assume that $(\Sig, E)$ is as in Theorem \ref{t4p} and possesses property T(1).

Then,  in the case of Lemma \ref{lp2-2}(i-c) we have
$DE\ge2EE_0=4$, and hence $\R E_0\cap F^+=\emptyset$. This implies that the
curve $\bn_*\widehat C_1$ intersects $E_0$ only in complex conjugate points, thus, the divisor $(\bn\big|_{\widehat C_1})^*(E_0)$
cannot consist of $\bn_*\widehat C_1\cdot E_0-1$ points.
In the case of Lemma \ref{lp2-2}(iii)
we have $DE\ge4$, and then $\R E_0\cap F^+=\emptyset$. It follows that the components $\widehat C_1,
\widehat C_2$ of $\widehat C$ intersect $\widehat C_3$ in complex conjugate points, and hence
$\widehat C_1$ and $\widehat C_2$ are complex conjugate. On the other hand, we obtain that $D=2E_0+2D'$ and
$$r-1=-DK_\Sig-\frac{DE}{2}-2\le-D'K_\Sig-\frac{D'E}{2}-1\ ,$$ which yields
$$-D'K_\Sig-\frac{D'E}{2}-1\le 0\ ,$$ thus, $r=1$ and any $\bw\in\Ima(\gamma)$ consists of one real point.
Then, both the curves
$\bn(\widehat C_1)$, $\bn(\widehat C_2)$ hit the real point $\bw$, which
contradicts the fact that $\bw$ must be a non-singular point of $\bn(\widehat C_1\cup
\widehat C_2)$.
The case of Lemma \ref{lp2-3}(i) with $k_0=2$ is excluded by the assumption $DE>0$.
In the case of
Lemma \ref{lp2-3}(i) with $k_0=4$, we have $DE=4$ and $r=1$, and hence the hypotheses
of Theorem \ref{t4p} yield
that $\R E_0\cap F^+=\emptyset$,
which, in particular, means that $\Ima(\gamma)$ is disjoint from $E\cup E_0$, which contradicts the appearence of the considered degeneration.

Finally, we exclude the case of Lemma \ref{lp2-3}(iii-a).
Indeed, each of the curves
$C_1,C_2$ has exactly one smooth local branch transversally crossing
$E$, and, since their real parts entirely (up to a finite set) lie
in $F^+$, they must be complex conjugate. It follows that $D=E+2D'$ and
$$r-1=-DK_\Sig-\frac{DE}{2}-2=-2D'K_\Sig-D'E-1\le-D'K_\Sig-\frac{D'E-1}{2}-1.$$
Thus, $-D'K_\Sig-\frac{D'E-1}{2}-1\le0$, and hence $r = 1$.
Again, we have that
each $\bw\in\Ima(\gamma)$ is a real point in $F^+$, which then belongs to both the complex conjugate curves
$\bn(\widehat C_1)$ and $\bn(\widehat C_2)$. This, however, contradicts the condition that $\bw$ must be a non-singular point
of $\bn(\widehat C_1)\cup\bn(\widehat C_2)$.
\proofend

\subsection{Proof of Theorems \ref{t2} and \ref{t4p}: Moving a real point
}\label{sec1}
Assume that $2m<r$.
In this case, all elements $[\bn:\PP^1\to\Sig]\in{\mathcal V}^{im,\R}_r
(\Sig,E,F,D,\bw)$
for a generic $\bw\in{\mathcal P}_{r,m}(\Sig,F^+)$
verify $\bn(\R\PP^1)\subset\overline F^{\;+}$.

Pick a
generic
($r - 1$)-tuple $\bw'=(w_1,...,w_{r-1})\in{\mathcal P}_{r-1,m}(\Sig, F^+)$
and two generic points $w^{(0)}_r$ and $w^{(1)}_r$ in $F^+$.
Choose
a segment $\sig \subset F^+$ of
a real part of some generic
smooth real algebraic curve
$\C\sig\subset\Sig$
such that $\sig$ starts at
$w^{(0)}_r$,
ends up at
$w^{(1)}_r$, and
avoids all the points of
$\bw'$.
Let $w: I \to F^+$ be a simple parameterization of $\sig$.
Lemmas \ref{lp1-3}, \ref{lp1-4}, and \ref{lp2-1} - \ref{lp2-3}
imply that the path $\gamma: I \to {\mathcal P}_{r, m}(\Sig, F^+)$
defined by $\gamma(t) = (\bw', w(t))$ is qualified.
Let $I_1 \subset I$ be a finite subset certifying that $\gamma$ is qualified.

\begin{lemma}\label{les-net}
No element of $\Gamma(\gamma)$ admits a description given
in Lemma \ref{lp2-2}(ii).
\end{lemma}

{\bf Proof.}
Assuming the contrary, we get that the complex conjugation
on $\widehat C$ must interchange the components $\widehat C_i$,
$i=1,2$. Indeed, otherwise, the curve $C=\bn_*\widehat C$ would have
a pair of real local branches transversally crossing $\R E$, which
is impossible. Then, both the curves $\bn(\widehat C_1)$ and $\bn(\widehat C_2)$ must have a common real point $w$,
which contradicts the condition that $w$ is a non-singular point of $\bn(\widehat C_1)\cup\bn(\widehat C_2)$.
\proofend

Put
\begin{eqnarray}
\overline{\mathcal V}^{\;im}_r(\Sig,E,D, \bw', \C\sig)&=&
\Big\{[\bn:
{\widehat C}
\to \Sig, (\widehat\bw', {\widehat w})] \in \overline{\mathcal M}^*_{0,r}(\Sig,D)\ :\nonumber\\
& &\quad[\bn:
{\widehat C}
\to\Sig]\in\overline{\mathcal V}^{\; l, im}(\Sig,E,D),\nonumber\\
& &\quad \bn(\widehat\bw')= \bw',
\ \bn({\widehat w}) \in \C\sig \Big\}. \nonumber
\end{eqnarray}
We have a
conjugation-invariant evaluation map
$$\Ev:{\overline{\mathcal V}}^{\;im}_r(\Sig,E,D,\bw',\C\sig)\to\C\sig,\quad\Ev([\bn: {\widehat C} \to \Sig,
\widehat\bw'\cup\{\widehat w\})])=\bn(\widehat w)\;.$$

\begin{lemma}\label{les1}
Let $\xi \in \Gamma(\gamma)$
and $\tau(\xi) \not\in I_1$.
Then, the function
$\tau_*\mu_{\Gamma_1(\gamma), \xi}$
is constant.
\end{lemma}

{\bf Proof.}
Let $\xi = [\bn: \PP^1 \to \Sig, \widehat\bw]$
and $\widehat\bw= (\widehat\bw', {\widehat w})$.
Since $\pi_r(\xi) \in \mathcal V^{im}(\Sig,E,D)$, the projection $\rho$ takes the germ
$(\mathcal V^{im}(\Sig,E,D),\pi_r(\xi))$ injectively to $|D|$.
Put $C=\rho(\pi_1(\xi))$.
The image $(\rho(\mathcal V^{im}(\Sig,E,D)),C)\subset|D|$
is determined by the conditions that
it induces an equigeneric deformation of each singular point $z\in C$, and
each local branch centered on $E$
intersects
$E$ in one point
(possibly moving)
with multiplicity $2$.
By Lemmas \ref{l4} and \ref{l5}(1), the germ $(\rho(\mathcal V^{im}(\Sig,E,D))
,C)$, naturally embedded into
\begin{equation}\prod_{z\in\Sing(C)\setminus E}B(C,z)\times\prod_{z\in C\cap E}B_E(C,z,m)
\ ,\label{espace}\end{equation} appears to be the intersection of smooth and linear subvarieties in the above space.
By formulas (\ref{e637}) and (\ref{e2234}), and by Lemma \ref{l637}, the transversality of that intersection, and thereby the smoothness
of the germ $(\rho(\mathcal V^{im}(\Sig,E,D)),C)$, are equivalent to the relation
\begin{equation}h^0(\PP^1,{\mathcal O}_{\PP^1}(\bd
))=r\ ,\label{e637k}\end{equation}
where $\deg\bd=D^2-(D^2+DK_\Sig+2)-l=-DK_\Sig-DE/2-2=r-1$, which holds true by Riemann-Roch in view of
$r-1\ge0>-2$. Furthermore, to specify the germ $\rho\pi_1(\Gamma(\gamma)),C)\subset
(\rho(\mathcal V^{im}(\Sig,E,D)),C)$, we impose the additional condition $\bn(\widehat\bw')=\bw'$.
Then the smoothness of the germ $\rho\pi_1(\Gamma(\gamma),C)$ amounts to the relation
\begin{equation}h^0(\PP^1,{\mathcal O}_{\PP^1}(\bd -\widehat\bw'))=1\ ,\label{e637a}\end{equation}
({\it cf.} (\ref{e637k})), which again holds by Riemann-Roch, since
$\deg(\bd-\widehat\bw')=(r-1)-(r-1)=0>-2$.
At last, a stronger than (\ref{e637a}) relation
\begin{equation}h^0(\PP^1,{\mathcal O}_{\PP^1}(\bd-\widehat\bw))=0\label{e637b}\end{equation} (which once again holds by Riemann-Roch due to $\deg(\bd-\widehat\bw)=-1>-2$) yields that $\tau(\rho\pi_1(\Gamma(\gamma)),C)\to(I,\tau(\xi))$ is
a diffeomorphism.
Thus, the constancy
of the function $\tau_*\mu_{\Gamma_1(\gamma), \xi}$ follows from Lemmas \ref{l4} and \ref{l5}(1).
\proofend

\begin{lemma}\label{les2}
Let $\xi = [\bn: \PP^1 \to \Sig, \widehat\bw] \in \Gamma(\gamma)$, where $\widehat\bw= (\widehat\bw', {\widehat w})$.
Assume that $\tau(\xi)\in I_1$ and
$\xi$ is as in Lemma \ref{lp2-1}(i).
Then, the function
$\tau_*\mu_{\Gamma_1(\gamma), \xi}$
is constant.
\end{lemma}

{\bf Proof.} Since $\bn$ is birational onto its image, the projection
$\rho:(\pi_1\Gamma(\gamma),\pi_1(\xi))\to|D|$ is injective.
Put $C=\rho\pi_1(\xi)$.
The germ $(\rho\pi_1\Gamma(\gamma),C)$, naturally embedded into the space (\ref{espace}),
is an intersection of not necessarily smooth subvarieties obtained as the product
(\ref{espace}), in which one factor $B(C,z)$ is replaced with $B^{eg}(C,z)$, $z\in
\Sing(C)\setminus E$, and smooth subvarieties, similarly obtained by changing one factor in the product
(\ref{espace}). In the considered case, we deduce the constancy of the function
$\tau_*\mu_{\Gamma_1(\gamma), \xi}$ from Lemma \ref{l5}(1) and \cite[Lemma 15]{IKS3}.
To apply \cite[Lemma 15]{IKS3}, we have to show the transversality of the intersection in the space
(\ref{espace}) of the tangent cones to the non-smooth subvarieties mentioned above with the tangent spaces to the
smooth members of the intersection forming the image of the germ $(\rho\pi_1\Gamma(\gamma),C)$.
By formulas (\ref{e637}) and (\ref{e2234}), and by Lemma \ref{l637} the required transversality
turns to be equivalent to relation (\ref{e637b}) that follows from Riemann-Roch in the same manner as in the proof
of Lemma \ref{les1}. \proofend

\begin{lemma}\label{les3}
Let $\xi = [\bn: \PP^1 \to \Sig, \widehat\bw] \in \Gamma(\gamma)$, where $\widehat\bw= (\widehat\bw', {\widehat w})$.
Assume that $\tau(\xi) \in I_1$ and
$\xi$ is as in Lemma \ref{lp2-1}(ii).
Then, the function
$\tau_*\mu_{\Gamma_1(\gamma), \xi}$
is constant.
\end{lemma}

{\bf Proof.} Show, first, that the local branch $P$ of the curve $C=\bn_*\PP^1$, intersecting $E$ with multiplicity $4$ is
smooth. Indeed, the case of a singular branch $P$ does not fit the hypotheses of Theorem \ref{t4p}, thus, we can assume that
$\Sig$ is of degree $2$. Furthermore, by Lemma \ref{lp2-1}(ii), the branch $P$ can be singular only if
$1\le r\le 2$. So, we have
$$DE=4+2s, \ s\ge0,\quad -DK_\Sig=r+s+3.$$
That is,
$-D(K_\Sig+E)=r-s-1\le1\ $, and hence,
the B\'ezout theorem implies
that $C$
is smooth along $E$.

Since $P$ is smooth,
we can
apply Lemma \ref{l5}(1) as in the proof of Lemma \ref{les2}
and reduce the claim on the constancy of the function $\tau_*\mu_{\Gamma_1(\gamma), \xi}$ to the transversality of the
intersection of the tangent cones and tangent spaces to the corresponding subvarieties
in the space (\ref{espace}). In view of the formulas (\ref{e637}), (\ref{e2234}), and (\ref{e2235}), and
in the virtue of Lemma \ref{l637}, we conclude that the required transversality amounts in the relation
(\ref{e637b}), where $\deg\bd=D^2-(D^2+DK_\Sig+2)-(l-2)-2=-DK_\Sig-DE/2-2$, which yields
$\deg(\bd-\widehat\bw)=-1>-2$, and we confirm (\ref{e637b}) by Riemann-Roch. \proofend

\begin{lemma}\label{les4}
Let $\xi = [\bn: \PP^1 \to \Sig, \widehat\bw] \in \Gamma(\gamma)$, where $\widehat\bw= (\widehat\bw', {\widehat w})$.
Assume that $\tau(\xi)\in I_1$ and
$\xi$ is as in Lemma \ref{lp2-1}(iii).
Then, the function
$\tau_*\mu_{\Gamma_1(\gamma), \xi}$
is constant.
\end{lemma}

{\bf Proof.}
The Zariski tangent space to $\overline{\mathcal V}^{\;im}(\Sig, E, D)$
at $[\bn:\PP^1\to\Sig]$, {\it i.e.}, the space of the first order deformations of
the map $\bn:\PP^1\to\Sig$, is contained in $H^0(\PP^1,{\mathcal N}^\bn_{\PP^1})$.
Furthermore, since the map $\bn:\PP^1\to\Sig$ has only two ramification points
mapped to $E$, \cite[Lemma 2.3]{CH} applies (see \cite[Remark 1 in Page 357]{CH}) and yields that
the Zariski tangent space is contained in $H^0(\PP^1,{\mathcal N}^\bn_{\PP^1}/\Tors{\mathcal N}^\bn_{\PP^1})$.
Note that the normal sheaf ${\mathcal N}^\bn_{\PP^1}$ has
two torsion points of order $2$, that is, the line bundle ${\mathcal
N}^\bn_{\PP^1}/\Tors{\mathcal N}^\bn_{\PP^1}$ has degree
$-DK_\Sig-4=2$. Thus, it follows from Riemann-Roch that
$$h^0(\PP^1,{\mathcal N}^\bn_{\PP^1}/\Tors{\mathcal N}^\bn_{\PP^1})=3=r(\Sig,D,l)\ ,$$ which is the dimension of the
germ of $\overline{\mathcal V}^{\;im}(\Sig, E, D)$
at $[\bn:\PP^1\to\Sig]$. Hence, the Zariski tangent space to that germ coincides with
$H^0(\PP^1,{\mathcal N}^\bn_{\PP^1}/\Tors{\mathcal N}^\bn_{\PP^1})$, and the germ itself is smooth.
Moreover, imposing the condition that the marked points $\widehat\bw'$ have fixed evaluation images and using the relation
$$h^0(\PP^1,{\mathcal N}^\bn_{\PP^1}/\Tors{\mathcal N}^\bn_{\PP^1}(-\widehat\bw'))=1$$
(following from Riemann-Roch), we obtain that the germ of $\Gamma(\gamma)$
at $\xi$ is smooth. Moreover, since again by Riemann-Roch
$h^0(\PP^1,{\mathcal
N}^\bn_{\PP^1}/\Tors{\mathcal N}^\bn_{\PP^1}(-\widehat\bw))=0$, we conclude that
the germ of $\Gamma(\gamma)$ at $\xi$ is
diffeomorphically mapped by $\tau$ onto the germ of $I$ at $\tau(\xi)$.

Since
$\gamma(I)\subset F^+$, we have $\widetilde\bn(\PP^1_\R)\subset F^+$. An
easy computation gives
$$\frac{D^2+DK_\Sig}{2}+1=4\left(\frac{C^2+CK_\Sig}{2}+1\right),\quad\text{where}\
\bn_*\PP^1=2C\;,$$ that is, independently of
$w\ne\bn(\widehat w)$ belonging to the germ of $\gamma(I)$ at $\bn(\widehat w)$,
the real nodes of
$\widetilde\bn_*\PP^1$ lie only in neighborhood of real nodes of
$C$:
four non-solitary nodes in a neighborhood of each non-solitary
node of $C$, and two solitary and two complex conjugate nodes in a
neighborhood of each solitary node of $C$, which immediately implies
the constancy of
$\tau_*\mu_{\Gamma_1(\gamma), \xi}$.
\proofend

\begin{lemma}\label{les5}
Let $\xi = [\bn: {\widehat C}_1 \cup {\widehat C}_2 \to \Sig, \widehat\bw] \in \Gamma(\gamma)$,
where $\widehat\bw= (\widehat\bw', {\widehat w})$.
Assume that $\tau(\xi) \in I_1$ and
$\xi$ is as in Lemma \ref{lp2-2}(i-a) or (i-b).
Then, the function
$\tau_*\mu_{\Gamma_1(\gamma), \xi}$
is constant.
\end{lemma}

{\bf Proof.}
If $\xi$ is as in Lemma \ref{lp2-2}(i-a),
then $\widehat\bw'$ splits into
two disjoint parts $\widehat\bw'_i\subset \widehat C_i$ such that
$|\widehat\bw'_i|=r(\Sig,E,D_i)$, $D_i=\bn_*\widehat C_i$, $i=1,2$.
We have \begin{equation}H^1(\widehat
C_i,{\mathcal N}^\bn_{\widehat
C_i}(-\bd_{0,i}-\widehat\bw'_i))=0,\quad H^0(\widehat C_i,{\mathcal
N}^\bn_{\widehat C_i}(-\bd_{0,i}-\widehat\bw'_i))=0\
,\label{e2243}\end{equation} where $2\bd_{0,i}=\bn^*(\bn_*\widehat
C_i\cap E)$, $i=1,2$. Indeed,
$\deg{\mathcal N}^\bn_{\widehat
C_i}(-\bd_{0,i}-\widehat\bw'_i)=-D_iK_\Sig-2-D_iE/2-r(\Sig,E,D_i)=-1>-2$, $i=1,2$, and hence
(\ref{e2243}) by Riemann-Roch.

Consider the normalization $\nu:\PP^1\amalg\PP^1\to\widehat C$ and
explore the exact sequence
\begin{equation}0\to \nu_*{\mathcal
N}^{\bn\circ\nu}_{\PP^1\amalg\PP^1}(-\bd_0-\widehat\bw')\overset{\alpha}{\longrightarrow}{\mathcal
N}^{\bn}_{\widehat
C}(-\bd_0-\widehat\bw')\overset{\beta}{\longrightarrow}{\mathcal
O}_{\widehat z}\to0\;,\label{e2242}\end{equation} where $\widehat
z=\widehat C_1\cap\widehat C_2$,
the map $\alpha$ is an isomorphism outside
$\widehat z$ and acts at $\widehat z$ as follows: identifying
$\left(\sigma_*{\mathcal
N}^{\bn\circ\sigma}_{\PP^1\amalg\PP^1}\right)(-\bd_0-\widehat\bw')_{\widehat
z}$ with $\C\{x\}\oplus\C\{y\}$, we have
$$\alpha_{\widehat z}(f(x),g(y))=xf(x)+yg(y)\in\left({\mathcal
N}^{\bn}_{\widehat C}\right)(-\bd_0-\widehat\bw')_{\widehat
z}\cong\C\{x,y\}/\langle xy\rangle\;.$$ Passing to cohomology in
(\ref{e2242}) and using (\ref{e2243}), we get the isomorphisms
$$H^0(\widehat C,{\mathcal
N}^{\bn}_{\widehat
C}(-\bd_0-\widehat\bw'))\overset{\beta}{\longrightarrow}H^0(\widehat
z,{\mathcal O}_{\widehat z})\overset{\bn}{\longrightarrow}{\mathcal
O}_{\Sig,z}/{\mathfrak m}_z\simeq\C,\quad z=\bn(\widehat z)\;,$$
which yield that the germ of $\overline{\mathcal
V}^{\;im}_{r-1}(\Sig,E,D,\bw')$ at $\xi$ is smooth,
one-dimensional, and its tangent space can be identified with a line
in $|D|$ spanned by $C=\bn_*\widehat C$ and some curve $C'$ which
avoids $z$. In particular, $C\cap C'$ is finite, and by our
assumptions, $\bn(\widehat w)\in C\setminus C'$, that immediately yields that
the germ of $\overline{\mathcal V}^{\;im}_r(\Sig,E,D,\bw')$
at $\xi$ is smooth and diffeomorphically mapped by $\Ev$
onto the germ $(\Sig,\bn(\widehat w))$.
This yields that $\tau$ diffeomorphically maps the germ of $\Gamma(\gamma)$ at $\xi$ onto
the germ of $I$ at $\tau(\xi)$ and that the function
$\tau_*\mu_{\Gamma_1(\gamma), \xi}$ remains constant in the considered case.

If $\xi$ is as in Lemma \ref{lp2-2}(i-b), we notice that the curve $C_2$ is immersed due to condition $T(1)$.
Then, the above argument in the same manner
yields the constancy of
$\tau_*\mu_{\Gamma_1(\gamma), \xi}$.
Indeed, notice that
$\widehat\bw'_2=\emptyset$ and replace in (\ref{e2243}) ${\mathcal
N}^\bn_{\widehat C_2}(-\bd_{0,2}-\widehat\bw'_2)$ with ${\mathcal
N}^\bn_{\widehat C_2}/\Tors{\mathcal N}^\bn_{\widehat C_2}$,
correspondingly, replace in (\ref{e2242}) ${\mathcal
N}^{\bn\circ\sigma}_{\PP^1\amalg\PP^1}(-\bd_0-\widehat\bw')$,
${\mathcal N}^{\bn}_{\widehat C}(-\bd_0-\widehat\bw')$ with
${\mathcal N}^{\bn\circ\sigma}_{\PP^1\amalg\PP^1}/\Tors{\mathcal
N}^{\bn\circ\sigma}_{\PP^1\amalg\PP^1}(-\bd'_{0,1}-\widehat\bw')$,
${\mathcal N}^{\bn}_{\widehat C}/\Tors{\mathcal N}^{\bn}_{\widehat
C}(-\bd'_{0,1}-\widehat\bw')$, where $\bd'_{0,1}=\bd'_0\cap\widehat
C_1$, and then use (\ref{e637b})
together with (\ref{e2243}).
\proofend

\begin{lemma}\label{les6}
Let $\xi = [\bn: {\widehat C}_1 \cup {\widehat C}_2 \to \Sig, \widehat\bw] \in \Gamma(\gamma)$,
where $\widehat\bw= (\widehat\bw', {\widehat w})$.
Assume that $\tau(\xi)\in I_1$ and
$\xi$ is as in Lemma \ref{lp2-3}(ii).
Then, the function
$\tau_*\mu_{\Gamma_1(\gamma), \xi}$
is constant.
\end{lemma}

{\bf Proof.}
By Lemma \ref{lp1-2}(3),
$\bn:\widehat C_1\to
C_1\hookrightarrow\Sig$ is an immersion, the curve $C_1$ is
non-singular along $E$, and is quadratically tangent to $E$ at $l=DE/2$
distinct points. Observe also that there is a well-defined morphism of ${\mathcal V}$, the germ at $\xi$
of the one-dimensional family $\bigcup_{w\in\C\sig}\overline{\mathcal V}^{\;im}_r(\Sig,E,D,\bw'\cup\{w\})$,
onto the germ $(\C\sig,
\bn({\widehat w}))$, which sends $[\bn':\widehat C'\to\Sig]$ to
the (unique) intersection point of $\bn'_*\widehat C'$ with
$(\C\sig,
\bn({\widehat w}))$. Hence, ${\mathcal V}$
can be identified with $V$, the germ at $C=\bn_*(\widehat C)$
of the family of curves $C'=\bn'_*(\widehat C)$ over the elements $[\bn':\widehat C'\to\Sig,\bp']\in{\mathcal V}$.
It follows from \cite[Proposition 2.8(2)]{MS2} that
the germ $V$ consists of $l+1=C_1E/2$
non-singular components, corresponding to smoothing out one of the
intersection points of $C_1$ and $E$. Respectively, real components
correspond to smoothings of real intersection points.
The tangent line to a component of $V$ is spanned by $C$ and
some curve $C'$ not containing $C_1$ (otherwise, we would have
$C'=C$). Hence, $\bn(\widehat w)\not\in C'$, which means that the considered component of $V$ is
diffeomorphically mapped onto $(\C\sig,w_0)$.
Consider a real component $V'$ of $V$, along which a real point $z\in
C_1\cap E$ smoothes out.
This component can be
uniformized in a conjugation-invariant way by a parameter $t\in(\C,0)$ so
that, in local $\Conj$-invariant coordinates $x,y$ in a
neighborhood of $z$ in $\Sig$ with $z=(0,0)$, $E=\{y=0\}$,
$C_1=\{y-2x^2=0\}$, $F^+=\{y>0\}$, $F^-=\{y<0\}$, the curves
$C^{(t)}\in V'$ are given by
$$y^2-2y(x^2+\alpha t)+\alpha^2t^2+yt\cdot O(x,y,t)+xt^2
(1+O(x,t))=0$$ $$\alpha=\const>0,\ t\in(\C,0)\ ,$$ ({\it cf.}
the same deformation in \cite[Formula (56)]{MS2}).

By \cite[Section 2.5.3(4) and Lemma 2.10]{MS2}, the geometry of
curves $C^{(t)}\in V'$ in a neighborhood of $z$ is described by the
deformation patterns $$y^2-2y(x^2+\alpha)+\alpha^2\quad  \text{for}
\ t>0,\quad\text{and}\quad y^2-2y(x^2-\alpha)+\alpha^2\quad \text{for}  \ t<0\;.$$
The former deformation pattern defines a curve with a non-solitary
real node in $F^+$, and the latter one
defines a curve with a solitary
node in $F^-$.
Thus,
the function $\tau_*\mu_{\Gamma_1(\gamma), \xi}$ is constant.
\proofend

\begin{lemma}\label{les8}
Let $\xi = [\bn: {\widehat C}_1 \cup {\widehat C}_2 \cup {\widehat C}_3 \to \Sig, \widehat\bw] \in \Gamma(\gamma)$,
where $\widehat\bw= (\widehat\bw', {\widehat w})$.
Assume that $\tau(\xi)\in I_1$ and
$\xi$ is as in Lemma \ref{lp2-3}(iii-b).
Then, the function
$\tau_*\mu_{\Gamma_1(\gamma), \xi}$
is constant.
\end{lemma}

{\bf Proof.}
Consider a real
irreducible component ${\mathcal V}$ of the germ at $\xi$ of (the one-dimensional family) $\overline{\mathcal
V}_{r-1}^{\;im}(\Sig,E,D,\bw')$. Let $t\in(\C,0)$ be its
conjugation-invariant uniformizing parameter.
Notice that the intersection points $z_1,z_2$ of $E$ and $E_0$ must be real, since one of them,
say, $z_1$ turns into a point of quadratic intersection with $E$, and the other smoothes out in the deformation along ${\mathcal V}$.
Furthermore, $\R E\cup\R E_0\subset F$. Indeed, otherwise
the curve $C_1=\bn(\widehat C_1)$ would intersect $E_0$ only in complex conjugate points, and hence the intersection
point of $\widehat C_1$ and $\widehat C_2$ cannot be real, which is a contradiction.
Observe also that the assumptions of
the lemma
are relevant only for del Pezzo pairs of degree $1$. Hence, the hypotheses of Theorem
\ref{t4p} yield that $C_1$ is disjoint from $E$.

Suppose, first, that $w_0=\bn({\widehat w}) \in E_0$.
Introduce a conjugation-invariant
local coordinate $\kappa:(E,z_1) \to (\C,0)$ of the
germ of $E$ at $z_1$.
There is a natural
morphism $\eta:{\mathcal V}\to(E,z_1)$ sending an element of ${\mathcal V}$ to the
intersection point of its image in $\Sig$ with $(E,z_1)$. This morphism can be expressed as $\kappa = t^n(\alpha+O(t))$,
with some natural $n\ge1$ and $\alpha\ne0$ real. Extend the germ $(E,z_1)$ up to a family of smooth curve germs
centered over an open subset of $E_0$ and transversal to $E_0$. A germ of that family close to $(E,z_1)$ intersects the images in $\Sig$
of the elements of ${\mathcal V}$ in two points, whose coordinates on the given germ have asymptotics of $t^n$. Hence, this holds
for almost all germs of the constructed family. Thus, in view of the general position of $w_0$ in $E_0$,
each element of ${\mathcal V}$ defines a pair of points
in $(\C\sig,w_0)$ which, in a
conjugation-invariant local coordinate $\kappa_0: (\C\sig,w_0) \to  (\C,0)$
of $(\C\sig,w_0)$, can be given by
$\kappa_0=t^n(\alpha_i+O(t))$, $\alpha_i\ne0$, $i = 0, 1$.
Particularly,
these formulas define two morphisms $\eta_i
:{\mathcal V}\to(\C\sig,w_0)$, $i = 0, 1$.
Consider
the families
$$\begin{matrix}w'_0,w'_1 & \hookrightarrow & {\mathcal C}'&
\hookrightarrow & {\mathfrak X}'=\Sig\times(\C,0)\\
\downarrow & {} & \downarrow & {} & \downarrow \\
(\C,0) & = & (\C,0) & = & (\C,0)\end{matrix}$$ where ${\mathcal C}'$
is the family of images of elements of ${\mathcal V}$ so that ${\mathcal C}'_t =
(\bn_t)_*
\PP^1$ as $t\ne0$, ${\mathcal C}'_0
= \bn_*\widehat C$,
and $w'_i(t)=\eta_i(t)$, $i=0,1$, and $t\in(\C,0)$. Blowing up
${\mathfrak X}'$ along $E_0\subset{\mathfrak X}'_0$, $n$ times,
we
get the families \begin{equation}\begin{matrix}w_0,w_1 &
\hookrightarrow & {\mathcal C}& \hookrightarrow & {\mathfrak X}\\
\downarrow & {} &
\downarrow & {} & \downarrow \\
(\C,0) & = & (\C,0) & = &
(\C,0)\end{matrix}\label{e2244}\end{equation} where ${\mathfrak
X}_0=\Sig\cup{\mathcal E}_0\cup...\cup{\mathcal E}_{n-1}$ with the
exceptional surfaces
${\mathcal E}_i\simeq{\mathbf F}_1$ (the plane blown up at one point)
such that
${\mathcal E}_j$ with $j>0$ are disjoint from $\Sig$,
and $\Sig\cap{\mathcal E}_0=E_0$. Here, $E_0$ is a section of ${\mathcal E}_0$ with
self-intersection $1$, disjoint from the $(-1)$-curve
$E'_0\subset{\mathcal E}_0$. Let $\theta_i$ be the fiber of ${\mathcal
E}_0$ through the point $z_i$, $i=1,2$, and $\theta_0$ the fiber
through $w_0\in E_0$. Observe that $w_1(0)\in\theta_2\setminus(E_0\cup
E'_0)$, and $w_0(0)\in\theta_0\setminus(E_0\cup E'_0)$.
The curve ${\mathcal C}_0$ contains
the components $C_1$ and $E$ in $\Sig\subset{\mathfrak X}_0$, whereas in ${\mathcal E}_0\simeq\F_1$ it
has several fibers and a conic $C_2$ crossing $E_0$ at $z_2$ and $z^*$ (the image of the node $\widehat C_1
\cap\widehat C_2$), passing through the point $p_0(0)$ (the limit position of $w_0$) belonging to some fiber $\theta_0$.
Since $\bn:\widehat C_2\to E_0$ is the double covering ramified at $E\cap E_0$, the conic
$C_2$ is tangent to the fiber $\theta_1$ at some point $p_1(0)$ and
tangent to the fiber $\theta_2$ at $z_2$. Note that these conditions define two
conics,
and they both must be real in the
situation considered.

We claim that the choice of a curve $C_1\cup C_2\subset{\mathfrak X}_0$ determines a unique smooth germ
$V$ regularly parameterized by the germ of $\C\sig$ at $w_0$. Indeed,
the considered situation fits the hypotheses
of \cite[Lemma 2.19]{MS2} ({\it cf.} also \cite[Lemma 2.15]{MS2}). Namely, in our situation,
the parameter $k$ in \cite[Lemma 2.19]{MS2} equals $2$, and the polynomial $h^{(1)}_z(x,y)$
describes the chosen conic $C_2$.
At last, the sufficient conditions for applying \cite[Lemma 2.19]{MS2} amount to the equality
\begin{equation}H^1(\widehat C_1,{\mathcal N}^{\bn}_{\widehat C_1}(-\bd_{0,1}-\widehat\bw'))=0\label{e637c}\end{equation}
that can be established in the same way as (\ref{e2243}),
and to the fact that the conditions imposed on the conic $C_2$ are transversal in the space of conics.
It follows that $n=1$, that the germ $V$ is smooth and is diffeomorphically mapped onto $(\C\sig,w_0)$, which finally yields the constancy of
$\tau_*\mu_{\Gamma_1(\gamma), \xi}$, since the solitary nodes of ${\mathcal C}_t$, $t\in(\R,0)$, $t\ne0$,
come only from those of the component $C_1$.

The case of $w_0\in C_1$ can be reduced to the preceding case. Namely, consider the germ of a smooth real curve transversally crossing $E_0$
at some generic real point,
and let $\sig'$ be the real part of this germ. Since $C_1$ intersects $E_0$ at some real point in $E_0\cap
F^+$, and $\widehat C_2$ doubly covers $E_0$ with ramification at $E_0\cap E$, we derive that $\sig'$ intersects with
each of the curves in $|D|$ induced by the real part of the germ ${\mathcal V}$ (in fact, we have two real intersection points, and we
choose one of them). The previous consideration yields that ${\mathcal V}$  is smooth. Since $w_0\in C_1$ is generic,
we obtain that ${\mathcal V}$ diffeomorphically maps onto $\C\sig$, and hence the constancy of $\tau_*\mu_{\Gamma_1(\gamma), \xi}$ follows.
\proofend

\subsection{Proof of Theorems \ref{t2} and \ref{t4p}: Moving a pair of complex conjugate points
}\label{sec1a}
Assume
that $r \ge 2m$ and $m\ge1$.
Pick a
generic
($r - 2$)-tuple $\bw'=(w_1,...,w_{r - 2})\in{\mathcal P}_{r - 2, m}(\Sig, F^+)$
and two generic points $w^{(0)}$ and $w^{(1)}$ in $\Sig \setminus \R\Sig$.
Choose
a segment $\sig \subset \Sig \setminus \R\Sig$
on some generic
smooth real algebraic curve
such that $\sig$ starts at
$w^{(0)}$,
ends up at
$w^{(1)}$, and
avoids all the points of
$\bw'$.
Let $w: I \to \Sig$ be a regular parametrization of $\sig$.
Lemmas \ref{lp1-3}, \ref{lp1-4}, and \ref{lp2-1}--\ref{lp2-3}
imply that the path $\gamma: I \to {\mathcal P}_{r, m}(\Sig, F^+)$
defined by $\gamma(t) = (w(t), \Conj(w(t)), \bw')$ is qualified.
Let $I_1 \subset I$ be a subset certifying that $\gamma$ is qualified.

Due to the generic choice of $\bw'$ and $\sig$ (and dimension arguments), we can suppose
that
$\Gamma(\gamma)$ avoids
the real elements $\xi\in\overline{{\mathcal V}}^{\;im}_{r-2}(\Sig,E,D,\bw')$
with $\pi_{r - 2}(\xi) \in {\overline{\mathcal V}}^{\;im}(\Sig,E,D)$
belonging to the equisingular strata of dimension $\le r-2$,
and that
$\sig$
avoids the
isolated singularities of the real elements
$\xi\in\overline{{\mathcal V}}^{\;im}_{r-2}(\Sig,E,D,\bw')$.
The
images $\bn(\hat C)$ for the real elements
$\xi=[\bn:\widehat C\to\Sig, \widehat\bw'] \in\overline{\mathcal V}^{\;im}_{r-2}(\Sig,E,D,\bw')$
such that $\pi_{r - 2}(\xi) \in {\overline{\mathcal V}}^{\;im}(\Sig,E,D)$
belong to the $(r-1)$-dimensional
strata, sweep a three-dimensional real-analytic variety
in
$\Sig$, and we suppose that $\sig$ crosses it transversally and only at the
points which are generic on the corresponding curves $\bn_*\widehat C$.

\begin{lemma}\label{l-ima}
Let $\xi \in \Gamma(\gamma)$, and let one of the following conditions
hold:
\begin{itemize}
\item either $\tau(\xi) \in I\setminus I_1$,
\item or $\tau(\xi) \in I_1$
and $\xi$ satisfies the hypotheses of one of the Lemmas \ref{lp2-1},
\ref{lp2-2}(i-a,i-b),
and \ref{lp2-3}(ii,iii-b).
\end{itemize}
Then, the function
$\tau_*\mu_{\Gamma_1(\gamma), \xi}$
is constant.
\end{lemma}

{\bf Proof.}
The proof literally coincides with that of Lemmas \ref{les1}--\ref{les8}.
Indeed, the key point of the argument consists in checking appropriate transversality conditions, which depend only on the complex data
(in particular, those based on the cohomology computations (\ref{e637a}),
(\ref{e637b}),
(\ref{e2243}), (\ref{e637c})
or those based on \cite[Proposition 2.8(2), Lemma 2.15, Lemma 2.19]{MS2}),
and they are the same both when moving either real or a pair of complex conjugate points of the point constraint.
\proofend

In view of Lemma \ref{les-net1},
it remains to examine the wall-crossing described in the following statement.

\begin{lemma}\label{les9}
Let $\xi = [\bn: {\widehat C}_1 \cup {\widehat C}_2 \to \Sig, \widehat\bw] \in \Gamma(\gamma)$.
Assume that $\xi$ is as in Lemma \ref{lp2-2}(ii).
Then, the function
$\tau_*\mu_{\Gamma_1(\gamma), \xi}$
is constant.
\end{lemma}

{\bf Proof.}
Clearly, there are no real fixed point, {\it i.e.}, $r=2m$.
Then, $\widehat C_i$, $i=1,2$ are complex conjugate as well as $C_i=\bn_*(\widehat C_i)$, $i=1,2$.
Furthermore, $D=D'+\Conj_*D'$ with $D'E$ odd
({\it cf.} the assertion of Theorem \ref{cor1}), the configuration $\bw'$ splits into
disjoint complex conjugate subsets
subsets $\bw'_1\subset C_1$ and $\bw'_2\subset C_2$, and
respectively $\widehat \bw=\widehat\bw'\cup\{\widehat w,\conj(\widehat w)\}$, where
$\bn(\widehat\bw')=\bw'$, $\bn(\widehat w)\in\sig$, splits so that
$\widehat\bw_1\cup\{\widehat w\}\subset\widehat C_1$, $\widehat\bw_2\cup\{\conj(\widehat
w)\}\subset\widehat C_2$, $\bn(\widehat\bw'_i)=\bw'_i$, $i=1,2$.
Since $r=2m\ge2$
and the pair $(\Sig,E)$ possesses property $T(1)$, Lemma
\ref{lp1-3}(iii)
yields that
the curve $C_1$ is immersed outside $E$, the divisor $(\bn\big|_{\widehat C_1})^*(E)$ consists of
$l_1$ double points and one simple point;
the same holds for $C_2$, where $l_2=l_1$.
The simple points of $(\bn\big|_{\widehat C_1})^*(E)$ and $(\bn\big|_{\widehat C_2})^*(E)$
glue up into the node of $\widehat C=\widehat C_1\cup\widehat C_2$, which is mapped to a real
point $z\in E$. In the deformation along the family $\overline{\mathcal V}^{\;im}_{r-1}(\Sig,E,D,\bw')$,
the local branches of
$C_1,C_2$ crossing at $z$, glue up into a smooth branch tangent to
$E$ in a nearby point. In particular, we obtain that $F\cap\R E\ne\emptyset$.

Choose local real coordinates $x,y$ so that $z=(0,0)$, $E=\{y=0\}$, \mbox{$C=C_1\cup C_2=\{x^2+y^2=0\}$}, and consider the versal deformation
$x^2+y^2+\eta_1x+\eta_2y+\eta_3$, $\eta_1,\eta_2,\eta_3\in(\C,0)$, of the singular point $z$ of $C$. Its base $B=
\{(\eta_1,\eta_2,\eta_3)\in(\C^3,0)\}$ contains a smooth
two-dimensional locus $B'=\{\eta_3=\eta_1^2/4\}$ of curves tangent to $E$. We claim that the (three-dimensional) germ of
$\overline{\mathcal V}^{\;l-1,\;im}_{r-2}(\Sig,E,D,\bw')$ at $(\xi,\widehat\bw')$ is smooth and
is isomorphically projected onto $B$: both statements follow from
\begin{equation}H^0(\widehat C_i,{\mathcal N}^\bn_{\widehat C_1}(-\bd_i-\widehat z_i-\widehat\bw'_i))=0,\quad i=1,2\;,
\label{eh0}\end{equation}
where $(\bn\big|_{\widehat C_i})^*(C_i\cap E)=2\bd_i+\widehat z_i$, $\bn(\widehat z_i)=z$
for $i=1,2$. Relation (\ref{eh0}) comes from the Riemann-Roch theorem and the $h^1$-vanishing
\begin{equation}H^1(\widehat C_i,{\mathcal N}^\bn_{\widehat C_1}(-\bd_i-\widehat z_i-\widehat\bw'_i))=0,\quad i=1,2\;,
\label{eh1}\end{equation} that in turn is a consequence of the inequalities
$$\deg{\mathcal N}^\bn_{\widehat C_1}(-\bd_i-\widehat z_i-\widehat\bw'_i)=-D_iK_\Sig-2-(D_iE-1)/2-1-\#\widehat\bw_i$$
$$=-D_iK_\Sig-D_iE/2-5/2-(-D_iK_\Sig-D_iE/2-3/2)=-1>-2,\quad i=1,2\ .$$

Hence, the germ of
$\overline{\mathcal V}^{\;im}_{r-2}(\Sig,E,D,\bw')$ at $(\xi,\widehat\bw')$, isomorphic to $B'$, is smooth. This yields the smoothness of the germ of
$\overline{\mathcal V}^{\;im}_r(\Sig,E,D,\bw')$
at $(\xi,\widehat\bw)$,
which, in view of (\ref{eh0}) (where we substitute $\bn(\widehat w)$ and $\conj(\bn(\widehat w))$ for
$\widehat z_1$ and $\widehat z_2$, respectively),
is isomorphically mapped by
 the evaluation map onto $(\Sig\times\Sig,(\bn(\widehat w),\conj(\bn(\widehat w))))$.
Hence, $\Gamma(\gamma)$
is smooth and is isomorphically mapped by $\Ev$ onto $(\sig,w_0)$.
The real one-dimensional global branch of a current curve is a
circle tangent to $E$; it collapses to the point $z$ and then
appears again, but may be on the other side of $E$. The constancy of the
function
$\tau_*\mu_{\Gamma_1(\gamma), \xi}$ follows.
\proofend

\subsection{Proof of Theorems \ref{t2} and \ref{t4p}: final arguments}\label{sec1b}
Suppose that the non-tangential uninodal DP-pair $(\Sig,E)$ lacks property $T(1)$.
Consider (a germ of) a real elementary deformation $f: (\mathfrak X, \mathfrak E) \to \Delta_a$
of the pair $(\Sig,E)$, in which all other uninodal DP-pairs
$(\Sig^t, E^t)$,
$t\ne0$, possess property $T(1)$.
We can also suppose that $D^2>0$ ({\it cf.} (\ref{ed2g0})).

Pick $r$ disjoint smooth analytic sections $w_i: \Delta_a \to \mathfrak X$,
$i=1,...,r$, such that, for any $t\in(-a,a)$,
the configurations $\bw(t)=\{w_i(t)\}_{i=1,...,r}\in{\mathcal P}_{r,m}(\Sig^t, F^{t, +})$ are
generic, that is, the sets ${\mathcal V}_r(\Sig^t, E^t, D, \bw(t))$,
are finite (see Lemma \ref{lem-fin})
and their elements satisfy conditions of Lemmas \ref{lp1-3}
and \ref{lp1-4}.
The following lemma shows
that the
numbers $RW_m(\Sig^t, D, F^{t, +}, \varphi^t,\bw(t))$
form a constant function, $t\in(-a,a)$, and thus reduces the statement to the case of
uninodal DP-pairs possessing property $T(1)$ settled in Sections \ref{sec1} and \ref{sec1a}.

Consider the one-dimensional family ${\mathcal V}\to\Delta_a\setminus\{0\}$ formed by the sets
${\mathcal V}_r(\Sig_t, E_t, D, \bw(t))$, $t\ne0$, and its closure $\overline{\mathcal V}\to\Delta_a$.

\begin{lemma}\label{l22b}
The set $\overline{\mathcal V}\setminus{\mathcal V}$ does not
contain elements with a reducible source curve.
Let $\xi_0=[\bn:\PP^1\to\Sig, \widehat\bw]\in{\mathcal V}_r(\Sig, E, D,\bw(0))$.

(1) The element $\xi_0$ is a center of a unique local branch $V$ of
$\overline{\mathcal V}$.

(2) If $\bn:\PP^1\to\Sig$ is a double covering onto the image, then each element of $V$ is a double covering onto its image.

(3) If $\bn:\PP^1\to\Sig$ is birational onto its image, then $V$ is smooth, and, if in addition
$\xi_0$ is real, the function
$\xi \mapsto \mu(F^+, \varphi, \xi)$
is constant on $V$.
\end{lemma}

{\bf Proof.}
Since $r>0$ and since we can assume that $DE>0$,
the lack of the property $T(1)$ does not affect irreducible curves in $\overline{\mathcal V}$.
Thus, it follows from Lemmas \ref{lp2-2} and \ref{lp2-3} that the reducible elements in
${\overline V}(\Sig,E,D)$ form strata of positive codimension, and hence cannot
appear in $\overline{\mathcal V}\setminus{\mathcal V}$ due to the general position of $\bw(0)$.

Assume that $\bn:\PP^1\to\Sig$
be birational onto its image.
From the hypotheses of Lemma \ref{lp1-3}, we derive relations (\ref{e637a}) and (\ref{e637b}).
This means,
first, that $[\bn:\PP^1\to\Sig, \widehat\bw]$ is not isolated in $\overline{\mathcal V}
\cup{\mathcal V}_r(\Sig, E, D,\bw(0))$,
and hence belongs to $\overline{\mathcal V}$.
Furthermore, by Lemma \ref{l5}(1)
the unique branch of $\overline{\mathcal V}$ centered at $[\bn:\PP^1\to\Sig, \widehat\bw]$ is smooth and
diffeomorphically projects onto $\Delta_a$. The constancy of
the
function $\xi \mapsto \mu(F^+, \varphi, \xi)$
follows from Lemmas \ref{l4} and \ref{l5}(1).

Assume that $\bn:\PP^1\to\Sig$ is a double covering of a curve $C\in|-K_\Sig+E_0|$ with a node
in $\Sig\setminus E$ and ramification at $C\cap E$ (see Lemma \ref{lp1-4}(ii)).
We claim that it cannot be deformed into
a map, birational onto its image.
Indeed, a map birational onto its image, obtained in a deformation of $\bn:\PP^1\to\Sig$,
would have an image with (at least) four nodes born from a node of
$C$, which contradicts the fact $p_a(D)=3$, $D=2(-K_\Sig+E_0)$.
\proofend

\subsection{Proof of Theorem \ref{cor1}}\label{sec4.5}
If $r>2m$, then (in the notation
of Theorem \ref{cor1}) we necessarily have
$\bn(\R P^1)\subset\overline{F^+}$, and hence
$RW^+_{m}(\Sig,E,F^+,\varphi,D,\bw)=RW_{m}(\Sig,E,F^+,\varphi,D,\bw)$. Thus, Theorem \ref{cor1}
follows from Theorems \ref{t2} and \ref{t4p}.

If $r=2m$,  split the set
${\mathcal V}^{im,\R}_r(\Sig,E,F,D,\bw)$ into subsets
${\mathcal
V}^{im,\R}_{r, +}(\Sig,E,F,D,\bw)$ and ${\mathcal
V}^{im,\R}_{r, -}(\Sig,E,F,D,\bw)$ specified by the conditions
\mbox{$\bn(\R P^1)\subset\overline{F^+}$} and \mbox{$\bn(\R P^1)\subset\overline{F^-}$}, respectively, as
$[\bn:\PP^1\to\Sig, \widehat\bw]\in{\mathcal V}^{im,\R}(\Sig,E,F,D,\bw)$.
Then, we follow the proof of Theorem \ref{t4p} in Sections \ref{sec1} and \ref{sec1a} and
consider all bifurcations relevant in the case $r=2m$ when the configuration $\bw$ varies along the path
introduced in Section \ref{sec1a}. Notice, first, that the bifurcation described in Lemma
\ref{les9} is excluded by the hypotheses of Theorem \ref{cor1} and, second, that in all other bifurcations
the sets ${\mathcal V}^{im,\R}_{r, +}(\Sig,E,F,D,\bw)$ and
${\mathcal V}^{im,\R}_{r, -}(\Sig,E,F,D,\bw)$ never mix with each other. Hence, in each case we get the constancy of the numbers
$RW^+_{m}(\Sig,E,F^+,\varphi,D,\bw)$ and $RW^-_{m}(\Sig,E,F^+,\varphi,D,\bw)$.

\section{Proof of Theorem \ref{t1}}\label{sec-t1}

Under
the hypotheses of Theorem \ref{t1}, by Proposition \ref{prop1}
we have to study only (germs of) real elementary deformations $f: (\mathfrak X, \mathfrak E) \to \Delta_a$
of nDP-pairs $(\Sig,E)$ of degree $1$, which are either non-tangential uninodal DP-pairs lacking property $T(1)$,
or binodal, cuspidal, or tangential DP-pairs which are non-ridged, other members
of the deformation being
non-tangential uninodal DP-pairs possessing property $T(1)$.

\begin{lemma}\label{l5.1}
Theorem \ref{t1} holds when either $DE=0$, or $D^2\le0$, or $r(\Sig,D,l)>0$.
\end{lemma}

{\bf Proof.} The case of $DE=0$ is covered by Theorems \ref{t3n} and \ref{t3}, and \cite[Theorem 6]{IKS3}.

In view of (\ref{ed2g0}) and $DE>0$, the case
$D^2\le0$ reduces to the situation of $D^2=0$, $DE=2$,
$r(\Sig,D,l)=0$, $\dim|D|=1$, and $p_a(D)=0$. Thus, we count smooth real rational curves in $|D|$
tangent to $E$. If $\R E\cap F=\emptyset$, the invariant vanishes for all deformation equivalent tuples
$T$. If $\R E\cap F\ne\emptyset$, we have an equivariant ramified double covering $\rho:E\to|D|\simeq\PP^1$, and
the value of $RW_0(\Sig,E,F^+,\varphi,D)$ depend only on the following topological data that remains invariant in the deformation class: whether the ramification points of $\rho$ are real, or complex conjugate, and the sign $(-1)^{C_{1/2}\circ\varphi}$ ({\it cf.} (\ref{esign})).
\proofend

Thus, further on, we suppose that $D^2>0$ and $DE>0$.

Pick $r$ disjoint smooth analytic sections $w_i: \Delta_a \to \mathfrak X$,
$i=1,...,r$, such that, for any $t\in(-a,a)$,
the configuration $\bw(t)=\{w_i(t)\}_{i=1,...,r}\in{\mathcal P}_{r,m}(\Sig^t, F^{t, +})$ is
generic.

As in Section \ref{sec1b}, consider the sets
$${\mathcal V}^l(\Sig^t, E^t, D, \bw(t)) = \{[\bn^t:\PP^1 \to \Sig^t, \widehat\bw^t] \in
{\mathcal V}^l(\Sig^t, E^t, D)\ :\
\bn^t(\widehat\bw^t) = \bw(t)\},
$$
where $t \in \Delta_a \setminus \{0\}$ and $l = DE_t/2$.
These sets are finite (see Lemma \ref{lem-fin})
and
form a family ${\mathcal V} \to \Delta_a \setminus \{0\}$..
Its closure $\overline{\mathcal V} \to \Delta_a$ is the union of irreducible one-dimensional components.
Denote by ${\mathcal V}^\R$
the real part of
${\mathcal V}$, and denote by $\overline{{\mathcal V}^\R}$ the closure
of ${\mathcal V}^\R \subset \overline{\mathcal V}$.

\begin{lemma}\label{l22c}
Assume that $(\Sig, E)$ is a non-tangential uninodal DP-pair which lacks property $T(1)$
and that $D^2>0$ and $DE>0$.
Then the following holds:

(1) The elements of $\overline{\mathcal V}^\R\setminus{\mathcal V}^\R$ are
represented
by maps $\bn:\PP^1\to\Sig$, which are birational onto their image.

(2) If either $r(\Sig,D,l)>0$ or $D\ne-2K_\Sig-E$, the projection $\overline{\mathcal V^\R} \to (-a, a)$
is a trivial covering and
the function $\xi^t \mapsto \mu(F^{t,+}, \varphi^t, \xi^t)$ is constant
along each sheet of that covering.

(3) If $D=-2K_\Sig-E$, and $\xi=[\bn:\PP^1\to\Sig]\in{\mathcal V}^1(\Sig,E,D)$ represents a real
curve with a cusp in $\Sig\setminus E$, then $\xi$ is a center of a singular branch $V$
of $\overline{\mathcal V}$, and one has
\begin{equation}\sum_{\xi^t\in V\cap f^{-1}(t)}\mu(F^{t,+}, \varphi^t, \xi^t)=0
\quad\text{for all}\quad t\in(-a,a),\ t\ne0\ .\label{elast}\end{equation}
\end{lemma}

{\bf Proof.}
Let us show that $[\bn:\widehat C\to\Sig]\in\overline{\mathcal V}^\R\setminus{\mathcal V}^\R$
with a reducible $\widehat C$ or with a multiple covering of the image cannot exist. If $r(\Sig,D,l)>0$, the absence of reducible elements follows from
Lemmas \ref{lp2-2} and \ref{lp2-3}, which yield that the reducible elements of
$\overline{\mathcal V}^l(\Sig,E,D)$ form substrata of positive codimension; hence they do not hit
the configuration $\bw(0)$ in general position.
Respectively, the absence of multiple covers for $r(\Sig,D,l)>0$ follows from
Lemma \ref{l22b}. If $r(\Sig,D,l)=0$ a possible reducible element
must be such that $\widehat C=\widehat C_1\cup\widehat C_2$, where $\widehat C_1\simeq\PP^1$ and
$\bn:\widehat C_1\to\Sig$
satisfies either conditions of Lemma \ref{lp1-2}(2) or conditions of Lemma \ref{lp1-2}(3), and $\bn_*\widehat C_2=sE_0$, $s\ge1$.
The case of Lemma \ref{lp1-2}(2iii) is excluded by the condition $D^2>0$, the cases (2ii) and (3i) of Lemma
\ref{lp1-2} are not possible, since $\bn_*\widehat C_1$ and $E_0$ must intersect. In the remaining cases (2i) and (3ii) of Lemma \ref{lp1-2},
we have $(\bn_*\widehat C_1)E_0=1$, and hence $s=1$ in view of $0\le DE_0=1-s$, but then $DE=DE_0=1$ contrary to condition (\ref{e1}).
Thus, the elements of $\overline{\mathcal V}^\R\setminus{\mathcal V}^\R$ are
represented by maps $\bn:\PP^1\to\Sig$ satisfying either conditions of Lemma \ref{lp1-2}(2) or conditions of Lemma \ref{lp1-2}(3).
The case of Lemma \ref{lp1-2}(2iii) is again excluded by the condition $D^2>0$, so the above maps $\bn:\PP^1\to\Sig$
are birational onto their image.

Under the hypotheses of statement (2), the elements $\xi_0=[\bn:\PP^1\to\Sig,\widehat w]
\in\overline{\mathcal V}\setminus{\mathcal V}$ are immersions along $\bn^{-1}(\Sig\setminus E)$
(see Lemmas \ref{lp1-2} and \ref{lp1-3}), and then the statement follows from Lemma \ref{l22b}
(in fact, Lemma \ref{l22b} concerns the case $r(\Sig,D,l)>0$, but literally the same argument
applies to the the case $r(\Sig,D,l)=0$).

Finally, if $D=-2K_\Sig-E$ and $\xi=[\bn:\PP^1\to\Sig]
\in\overline{\mathcal V}\setminus{\mathcal V}$ represents a curve with a cusp in $\Sig\setminus E$,
then we have the $h^0$-vanishing relation (\ref{e637b}),
which brings the considered case to the framework of \cite[Lemma 14]{IKS3}: the branch
of $\overline{\mathcal V}$ at $\xi$ is isomorphic to the discriminantal semicubical parabola in the versal deformation base of the ordinary cusp, and hence the equality
(\ref{elast}) ({\it cf.} \cite[Lemma 2.6(2)]{W1}).
\proofend

\smallskip

We now pass to elementary deformations of binodal, cuspidal, and tangential DP-pairs,
and derive the statement of the theorem from
the following lemmas.

\begin{lemma}\label{l22}
(1) Let $(\Sig, E)$ be a binodal DP-pair. Then, the elements of $\overline{\mathcal V}^\R\setminus{\mathcal V}^\R$ are
represented
\begin{enumerate}\item[(1i)] either by maps $\bn:\PP^1\to\Sig$,
\item[(1ii)] or
by maps $\bn:\widehat C_1\cup\widehat C_2\to\Sig$ such
that
\begin{itemize}\item $\widehat C_1\cup\widehat C_2$ is a nodal tree of
$\PP^1$'s; \item $\widehat C_2$ is the
disjoint union of $s\ge1$ copies of $\PP^1$ isomorphically mapped
onto $E'$;
\item $\widehat C_1\simeq\PP^1$
intersects each component of $\widehat C_2$ at one point,
$C_1=\bn_*\widehat C_1\in|D-sE'|$, and $\bn:\widehat C_1\to\Sig$ satisfies conditions of
Lemma \ref{p3}(1-a).
\end{itemize}\end{enumerate}

(2) Let $(\Sig, E)$ be a cuspidal DP-pair. Then, the elements of $\overline{\mathcal V}^\R\setminus{\mathcal V}^\R$ are
represented by maps $\bn:\PP^1\to\Sig$ satisfying conditions of Lemma
\ref{p3}(1-a).

(3) Let $(\Sig, E)$ be a tangential DP-pair. Then, the elements of $\overline{\mathcal V}^\R\setminus{\mathcal V}^\R$ are
represented
\begin{enumerate}
\item[(3i)] either by maps $\bn:\PP^1\to\Sig$ satisfying conditions of Lemma
\ref{p3}(1-a),
\item[(3ii)] or by maps $\bn:\widehat C_1\cup\widehat C_2\to\Sig$, where $\widehat C_1\simeq\PP^1$ is birationally taken onto a curve
$C\ne E_0$, $\widehat C_2$ is the disjoint union of $s\ge1$ copies of $\PP^1$ isomorphically taken onto $E_0$, and $\widehat C_1$
intersects each component of $\widehat C_2$ at one point;
furthermore,
the divisor $\bn^*(E_0)$ consists of $CE_0 \geq s$ simple points.
\end{enumerate}\end{lemma}

{\bf Proof.}
{\bf (1)}
Let $(\Sig, E)$ be a binodal DP-pair.
By the semistable reduction theorem, the limit maps are
$\bn:\widehat C_1\cup\widehat C_2\cup \widehat C_3\to\Sig$, where
$\widehat C_1\cup\widehat C_2\cup\widehat C_3$ is a nodal tree of
$\PP^1$'s, the components of $\widehat C_2$ are mapped onto $E'$,
the images of the components of $\widehat C_1$ are curves different
from $E'$, and the components of $\widehat C_3$ are contracted to
points.

Let us show that $\widehat C_1\simeq\PP^1$. Indeed,
$\bn_*\widehat C_1\in|D-sE'|$ for some $s\ge0$. Observe that
$r(\Sig,D-sE',l)=r(\Sig,D,l)=r$, and that $[\bn:\widehat
C_1\to\Sig_0]\in\overline{{\mathcal V}_l(\Sig, E, D-sE')}$. Then,
the irreducibility of $\widehat C_1$ can be proved as Lemmas \ref{lp2-2} and \ref{lp2-3},
where all strata formed by maps of reducible curves are shown to be of
intersection dimension $<r$:
indeed, the referred argument is entirely based on the dimension count of
Lemmas \ref{lp1-2} and \ref{lp1-3}, which is identical to that in Lemma
\ref{p3}, and hence holds in our situation.

If $\bn:\widehat C_1\to\Sig$ satisfies conditions of Lemma \ref{p3}(1-a), then $\bn_*\widehat C_1$
intersects $E'$ transversally.
To prove the statement (1) of Lemma, it remains to show that
$\bn:\widehat C_1\to\Sig$ cannot satisfy conditions of Lemma \ref{p3}(1-b) or Lemma \ref{p3}(1-c).
The case of Lemma \ref{p3}(1-b) is
excluded by the assumption $D^2>0$.
Suppose that $\bn:\widehat C_1\to\Sig$ satisfies conditions of Lemma \ref{p3}(1-c).

First, $\bn:\widehat C_1\to\Sig$ cannot be a double cover of a uninodal rational curve $C$. Indeed,
in this case $C\in|-K_\Sig+E_0|$, and hence
$C$ intersects with any $(-1)$-curve, disjoint from $E$, with multiplicity $2$. Thus, $CE'=0$, which yields
$\widehat C_2=\emptyset$. Furthermore, the arithmetic genus $p_a(D)$ of the divisor $D=2C$ equals $3$, whereas
any deformation of $\bn:\widehat C_1\to\Sig$ would exhibit at least four nodes in a neighborhood of the node of $C$, a
contradiction.

Thus, $\bn:\widehat C_1\to\Sig$ is a double covering of a smooth rational curve $C$.
Since $CE=2$, $CE_0=0$, we have $0\le CE'\le2$. If $CE'=0$, then $\widehat C_2=\emptyset$, which, however, is excluded by the assumption
$D^2>0$. If $CE'=1$, then the relation $D^2>0$ yields that $D=2C+E'$.
Blow up the family $\mathfrak X \to \Delta_a$
along the curve $C \subset \Sig$. Then, the central fiber of the obtained family is the union of $\Sig$ and $\F_0 \simeq \PP^1 \times \PP^1$,
intersecting along $C$. Notice that $r(\Sig,D,DE/2)=1$ and $DC=E'C=1$. Hence,
the limit of the family
${\mathcal V}$ (over $\Delta_a\setminus\{0\}$) in the central fiber has the following curve in
the component $\F_0$: the fibers $\Phi$ passing through the intersection point of $C$ and $E'$, different from $w_1(0)$, and a rational curve $T\in|\Phi+2C|$, intersecting $C$ at $w_1(0)$ and tangent to the two fibers passing through the points of $C\cap E$. However, then $T$ and $\Phi$ intersect in two distinct points, and only one of these points
smoothes out in the deformation along the family ${\mathcal V}$ leaving a node of a rational
curve $(\pi_t)_*\PP^1\in|D|$ contrary to $p_a(D)=0$.
If $CE'=2$, then the same relation $D^2>0$ yields that
$D=2C+sE'$, $1\le s\le 3$. Again we blow up the family $\mathfrak X \to \Delta_a$ along the curve
$C \subset \Sig$ and obtain in the component $\F_0$ of the central fiber the limit curve,
consisting of two $s$-multiple fibers $\Phi_1,\Phi_2$, passing through the points of $C\cap E'$,
and a curve $T\in|\Phi_1+2C|$, intersecting $C$ at $w_1(0)$ and tangent to the two fibers passing through the points of $C\cap E$. Notice that $T$ meets $\Phi_1\cup\Phi_2$ at four distinct points. Hence,
\begin{itemize}
\item
if $s=1$, three of the points $T\cap(\Phi_1\cup\Phi_2)$ persist in the deformation against
$p_2(D)=2$,
\item
if $s=2$, then (at least) two of the points $T\cap(\Phi_1\cup\Phi_2)$ persist in the deformation bearing $4$ nodes against $p_a(D)=3$,
\item
and if $s=3$, then (at least) one of the points
$T\cap(\Phi_1\cup\Phi_2)$ persist in the deformation bearing $3$ nodes against $p_a(D)=2$.
\end{itemize}

\smallskip

{\bf (2)} Let $(\Sig, E)$ be a cuspidal DP-pair.
If an element of $\overline{\mathcal V}^\R\setminus{\mathcal V}^\R$ is represented by
a map $\bn:\widehat C\to\Sig$ such that $\bn_*\widehat C\not\supset E'$, then as in
the proof of (1)
we get $\widehat C\simeq\PP^1$, and $\bn:\PP^1\to\Sig$ must satisfy conditions
of Lemma
\ref{p3}(1-a).

So, suppose that there exists an element $[\bn:\widehat C_1\cup\widehat C_2\to\Sig]\in\overline{\mathcal V}^\R
\setminus{\mathcal V}^\R$ such that $\bn_*\widehat C_2=sE'$, $s\ge1$, and $C=\bn_*\widehat C_1\not\supset E'$.
Since $E$ and $E'$ are real and intersect in one point,
one has $\R E'\ne\emptyset$ and $\R E'\cap F=\emptyset$.
Thus, the intersection $C\cap E'$ lifts to several pairs of complex conjugate points on $\widehat C_1\simeq\PP^1$, $\widehat C_2$ is the disjoint union of $2k$ copies of $\PP^1$, which form $k$ complex conjugate pairs
respectively attached to complex conjugate pairs in $\bn^*(C\cap E')$. Furthermore, each component
of $\widehat C_2$ covers $E'$ with even degree and even ramifications at $E\cap E'$, thus,
for example $s\ge4k$. Since
$$-D(K_\Sig+E)=(C+sE')(E'+E_{-1})=CE_{-1}+(CE'-s)\ ,$$ the $2k$ local branches of
$\widehat C_1$ at $\widehat C_1\cap
\widehat C_2$ are mapped by $\bn$ to $2k$ local branches of $C$ centered on $E'$
and
totally intersecting $E'$ with multiplicity at least $s$. Since $CE=DE-s=2l-s$ and $r(\Sig,C,l-s/2)=
r(\Sig,D,l)+s/2$, Lemma \ref{gus}, applied to the family ${\mathcal V}_{l-s/2}(\Sig,E,C)$, allows $2k$
branches of $C$ centered on $E'$ to have at most $2k+s/2$ as the total intersection multiplicity
with $E'$. Hence, in particular, $s=4k$, the total intersection multiplicity with $E'$
of those local branches of $C$ centered on $E'$, which do not glue up with $E'$,
equals $CE'-s=CE'-4k$, and each component of $\widehat C_2$ doubly covers $E'$
with ramification at $E\cap E'$ (and some other point). However, then
each curve in ${\mathcal V}$ must have at least
$$p_a(C)+4k(CE'-4k)+2k>p_a(C)+4k(CE'-4k)=p_a(D)$$ nodes, where the summand $2k$ in the
left-hand side correspond to additional nodes appearing in local deformations of $2k$ tangency points
of $C$ and $E'$. Thus, $k=0$ and $\widehat C_2=\emptyset$, proving the second statement of
Lemma.

\smallskip

{\bf (3)} Let $(\Sig, E)$ be a tangential DP-pair. If an element of $\overline{\mathcal V}^\R\setminus{\mathcal V}^\R$ is represented by
a map $\bn:\widehat C\to\Sig$ such that $\bn_*\widehat C\not\supset E_0$, then as in
the proof of (1)
we get $\widehat C\simeq\PP^1$, and $\bn:\PP^1\to\Sig$ must satisfy conditions
of Lemma
\ref{p3}(1-a).

So, suppose that there exists an element $[\bn:\widehat C_1\cup\widehat C_2\to\Sig]\in\overline{\mathcal V}^\R
\setminus{\mathcal V}^\R$ such that $\bn_*\widehat C_2=sE_0$, $s\ge1$, and $C=\bn_*\widehat C_1\not\supset E_0$. Observe that $$r(\Sig,D,DE/2)=-DK_\Sig-1-\frac{DE}{2}=-CK_\Sig-1-\frac{CE}{2}=r(\Sig,C,CE/2)\ .$$
Hence, $\widehat C_1\simeq\PP^1$, and $\bn:\widehat C_1\to\Sig$ satisfies conditions of Lemma
\ref{p3}(1-a) (the cases (1-b) and (1-c) are excluded by the assumptions $D^2>0$ and $CE_0>0$). Since
$CE_0=DE_0+s$, at least $s$ local branches of $C$ centered on $E_0$ must glue up with $\bn_*\widehat C_2$, the rest
of statement (3) follows.\proofend

\begin{lemma}\label{l22e}
Let $(\Sig, E)$ be a binodal, cuspidal or tangential DP-pair.
Then, the projection $\overline{\mathcal V^\R} \to (-a, a)$
is a trivial covering and
the function $\xi^t \mapsto \mu(F^{t,+}, \varphi^t, \xi^t)$ is constant
along each sheet of that covering.
\end{lemma}

{\bf Proof.} First, we complexify our family of surfaces into
$\mathfrak X \to \Delta_a=\{z\in\C\ :\ |z|<a\}$, forget the real structure and reduce this family to a trivial one as follows.
If $(\Sig, E)$ is binodal (respectively, cuspidal or tangential), we take a model of $\Sig$
as in Lemma \ref{negtype},
and denote by $z_8$ one of the blown up points in $C'_2\setminus C_2$ (respectively, in $C_1\setminus C_2$).
Identify
$\Delta_a$ with the germ of a smooth plane curve $\Lam$, passing through $z_8$ and transversal to $C'_2$
(respectively, to $C_1$),
fix seven blown up points and vary $z_8$ along $\Lam$, and then perform the blowing down
$\Pi: {\mathfrak X} \to {\mathfrak X}' \simeq \Sig' \times\Delta_a$ of the family of exceptional divisors $E_8$ to the family
of points $z_8(t)\in\Lam$. We also take $\bw(t)\equiv\bw$ to be the constant family of configurations
in $\Sig'$. The family of maps $V$ turns into the family $V'$ of maps $\nu':\PP^1\to\Sig'$ subject to conditions:
\begin{itemize}\item $\bn'_*\PP^1\in|D'|_{\Sig'}$, where $D'=\Pi_*(D)$, \item $\bn'$ is an immersion, the curve
$C'=\bn'_*\PP^1$ passes through $\bw$, is simply tangent to $E$ at $l=DE/2=D'E/2$ distinct points,
and has an ordinary singular point of multiplicity $d_8=DE_8$ at some point $z_{8,\bn'}\in\Lam$.
\end{itemize}

Given a branch $B$ of $\overline V'$,
we show that the
map $[\bn':\PP^1\to\Sig']\in B\cap V'\mapsto z_{8,\bn'}\in\lam\simeq\Delta_a$ is
1-to-1 onto $\Delta_a\setminus\{0\}$. We argue on the contrary, assuming that $B$ multiply covers $\Delta_a$.
Consider two preimages $[\bn':\PP^1\to\Sig']$, $[\bn'':\PP^1\to\Sig']$ of the same point $z_8\in\Lam$.
Denote $C'=\bn'_*\PP^1$, $C''=\bn''_*\PP^1$. We will see that
\begin{equation}C'C''>(D')^2\ ,\label{e-new1}\end{equation} providing the required contradiction.

Suppose that $B$ is centered at $[\bn'_0:\PP^1\to\Sig']$. Then $C',C''$ appear in the equinormalizable deformation of
the curve $C'_0=\bn'_{0,*}\PP^1$, and hence in view of
\begin{equation}-K_{\Sig'}D'=-K_\Sig D+d_8,\quad l=DE/2=D'E/2\ ,\label{e-new2}\end{equation}
and by Lemma \ref{gus},
$$C'C''\ge2\delta(C'_0)+d_8+\frac{D'E}{2}+r(\Sig,D,l)\qquad\qquad\qquad\qquad\qquad\qquad$$ $$=(D')^2+K_{\Sig'}D'+2+d_8+\frac{DE}{2}+
\left(-DK_\Sig-\frac{DE}{2}-1\right)$$ $$=(D')^2+1\ .\qquad\qquad\qquad\qquad\qquad\qquad\qquad\qquad\quad\;\quad$$

Suppose that $B$ is centered at $[\bn'_0:\widehat C\to\Sig']$ with a reducible source curve $\widehat C$. By Lemma \ref{l22},
this is possible only if $(\Sig, E)$ is either binodal, or tangential.

Suppose that $(\Sig, E)$ is binodal, and $B$ is centered at $[\bn'_0:\widehat C_1\cup\widehat C_2\to\Sig']$ with $\widehat C_2$ consisting of
$s\ge1$ copies of $\PP^1$ mapped onto $E'_*
=\Pi(E')$. Then the curve $\bn'_{0,*}\widehat C_1=C'_1\in|D'-sE'_*|_{\Sig'}$ has an ordinary $(d_8-s)$-multiple singularity
at $z_8\in E'_*$ and $C'_1E'_*-(d_8-s)=D'E'_*+2s-d_8$ smooth branches centered on $E'_*\setminus\{z_8\}$, and, in the deformation along
$B$, $s$ of these branches glue up with $E'_*$ and the other $D'E'_*+s-d_8$ branches persist. In particular, in a neighborhood of each of the
latter set of branches of $C'_1$, the curves $C',C''$ have at least $2s$ intersection points. Hence, in view of (\ref{e-new2}),
relations $(E'_*)^2=-1$ and $-K_{\Sig'}E'_*=1$, and Lemma \ref{gus},
$$C'C''\ge2\delta(C'_1)-(d_8-s)(d_8-s-1)+d_8^2+\frac{D'E}{2}+r(\Sig,D,l)+2s(D'E'_*+s-d_8)$$
$$=(C'_1)^2+K_{\Sig'}C'_1+2-d_8^2+2sd_8-s^2+d_8-s+d_8^2+\frac{D'E}{2}\qquad\qquad\quad$$ $$+
(-K_\Sig D-\frac{DE}{2}-1)+2sD'E'_*+2s^2-2sd_8\qquad\qquad\qquad\qquad$$ $$=((D')^2-2sD'E'_*-s^2)+(-K_{\Sig'}D'+s)+1+s^2+d_8-s\qquad\qquad\quad$$
$$-K_\Sig D+2sD'E'_*\qquad\qquad\qquad\qquad\qquad\qquad\qquad\qquad\qquad\qquad$$
$$=(D')^2+1\ .\qquad\qquad\qquad\qquad\qquad\qquad\qquad\qquad\qquad\qquad\qquad\qquad\quad$$

Suppose that $(\Sig, E)$ is tangential, and $B$ is centered at $[\bn'_0:\widehat C_1\cup\widehat C_2\to\Sig']$,
where $\widehat C_2$ consists of $s\ge1$ copies of $\PP^1$ mapped onto $E'_0=\Pi(E_0)$. Observe that the curve
$\bn'_{0,*}\widehat C_1=C'_1\in|D'-sE'_0|_{\Sig'}$ has an ordinary $(d_8-s)$-multiple singularity
at $z_8\in E'_0$ and $C'_1E'_0-(d_8-s)=D'E'_0+s-d_8$ smooth branches centered on $E'_0\setminus\{z_8\}$, and, in the deformation along
$B$, $s$ of these branches glue up with $E'_0$ and the other $D'E'_0-d_8$ branches persist. Notice also that
the local branch of each of $s$ copies of $E'_0$ at the point $z_0=E'_0\cap E$ deforms into
a smooth branch quadratically tangent to $E$ in a neighborhood of $z_0$, and hence the curves
$C',C''$ have at least $s^2$ intersection points in a neighborhood of $z_0$. Thus, in view of
(\ref{e-new2}), relations $(E'_0)^2=0$ and $-K_{\Sig'}E'_0=2$, and Lemma \ref{gus} we obtain
$$C'C''\ge2\delta(C'_1)-(d_8-s)(d_8-s-1)+d_8^2+\frac{C'_1E}{2}+r(\Sig,D,l)+2s(D'E'_0-d_8)+s^2$$
$$=((C'_1)^2+K_{\Sig'}C'_1+2)-d_8^2+2sd_8-s^2+d_8-s+d_8^2+\left(\frac{D'E}{2}-s\right)\quad$$
$$+(-K_\Sig D-\frac{DE}{2}-1)+2sD'E'_0-2sd_8+s^2\qquad\qquad\qquad\qquad\quad$$
$$=((D')^2-2sD'E'_0+K_{\Sig'}D'+2s+2)+d_8-2s+\frac{D'E}{2}\qquad\qquad\qquad\quad$$
$$-K_\Sig D-\frac{DE}{2}-1+2sD'E'_0\qquad\qquad\qquad\qquad\qquad\qquad\qquad\qquad$$
$$=(D')^2+1\ ,\qquad\qquad\qquad\qquad\qquad\qquad\qquad\qquad\qquad\qquad\qquad\qquad\quad
$$ which is a contradiction.

Thus, it follows that the projection $\overline{\mathcal V^\R} \to (-a, a)$
is a trivial covering as required in the lemma.
Furthermore, we observe that the Welschinger signs of generic curves in each sheet of that covering
for positive and negative values of the parameter are the same.
Indeed,
\begin{itemize}
\item in the case of a bifurcation through an irreducible curve,
the distribution of solitary and non-solitary nodes persists by Lemmas \ref{l4} and \ref{l5};
\item in the case of a binodal bifurcation through a curve $C_1\cup sE'$, $s\ge1$,
the real nodes of generic curves in
$B^\R$ come

- either from the singularities of $C_1$, and their contribution to the Welschinger sign
is constant by Lemmas \ref{l4} and \ref{l5},

- or from intersections of $C_1$ and $E'$, which
produce only non-solitary nodes in the component of
$\R\Sig_t\setminus E$ that merges to the component of $\R\Sig\setminus E$ containing $\R E'$, and
the parity of their number is constant,

- or from self-intersections of components $E'$, in which case the real nodes again appear only
in the component of $\R\Sig\setminus E$ containing $\R E'$, and the parity of the numbers of solitary and non-solitary nodes depends only on the numbers of components $E'$ which are images of
real or complex conjugate components of $\widehat C_2$;
\item in the case of a tangential bifurcation through a curve $C_1\cup sE_0$, $s\ge1$, the only new
set of nodes to be considered are those, which pop up in a neighborhood of the point $E\cap E_0$:
here, solitary nodes appear as intersections of complex conjugate local branches tangent to $E$; namely,
in suitable conjugation-invariant coordinates in a neighborhood of the intersection point $E\cap E_0$,
we can take $E_0=\{y=0\}$, $E=\{y=x^2\}$, then a pair of local branches tangent to $E$ can be written as $y=2\xi x+\xi^2$, $y=2\overline\xi x+\overline\xi^2$, $\xi\in(\C,0)\setminus\R$, and hence their intersection point is $x=-(\xi+\overline\xi)/2$, $y=-|\xi|^2$, which means that the solitary nodes
under consideration always belong to the component of $\R\Sig_t\setminus E$ that merges to the same component of $\R\Sig\setminus E$.
\end{itemize}

Our final remark is that, for a real branch of $\overline V$, the nodes of current curves in a neighborhood of $E'$ come from
intersections of $C_1$ with $E'$, thus, are not solitary. \proofend

The proof of Theorem \ref{t1} is completed.

\section{Concluding remarks}\label{sec-cr}

\subsection{Examples}\label{sec-se}

Here, we exhibit a few
examples
where the relative invariants of uninodal DP-pairs can be found in a rather simple way.

Let $C_2\subset\PP^2$ be a smooth real conic with $\R C_2\simeq S^1$, and $\pi:\PP^2_{(0,3)}\to\PP^2$ the blow up
at three pairs of complex conjugate points on $C_2$. Then $(\PP^2_{(0,3)},E)$ is a real uninodal DP-pair of degree $3$, where
$E$ is the strict transform of $C_2$. Observe that $\R\PP^2_{(0,3)}\setminus\R E\simeq\R P^2\setminus\R C_2$ consists
of an open disc $F^{o}$ (orientable component) and an open M\"obius band $F^{no}$ (non-orientable component). Thus, we have a series of invariants
$RW_m(\PP^2_{(0,3)},E,F^+,\varphi,D)$, for $F^+=F^{o}$ or $F^{no}$,
all $\Conj$-invariant classes
$\varphi\in H_2(\PP^2_{(0,3)}\setminus \R\PP^2_{(0,3)})$, real effective divisor classes $D\in\Pic_+(\PP^2_{(0,3)},E)$ matching condition (\ref{e1}),
and integers $0\le m\le r/2$, $r=-DK-DE/2-1$.

In a particular case $D=dL$, $d\ge1$,
where $L$ is
the pull-back of a generic line, the irreducible curves in the linear system $|D|$ do not intersect
the exceptional
divisors of the blow up $\pi:\PP^2_{(0,3)}\to\PP^2$. Thus,
by blowing down the exceptional
divisors,
we can consider the pair $(\PP^2,C_2)$,
relative invariants $RW_m(\PP^2,C_2,F^+,\varphi,D)$
by letting
$$RW_m(\PP^2,C_2,F^+,\varphi,D)= RW_m(\PP^2_{(0,3)},
E, F^+,\varphi,D), $$
and
observe that it counts those
real plane rational curves of degree $d$ passing through
a generic tuple $\bw\in{\mathcal P}_{r,m}(\PP^2,F^+)$,
where $r=2d-1$, $0\le m<d$,
that  are simply tangent to $C_2$ at $d$ distinct points
({\it cf.} Lemma \ref{lp1-3}).

Since the class of $C_2$ is divisible by $2$, we have a double
covering $\rho:Q\to\PP^2$ ramified along $C_2$. The two ruling linear systems $|L_1|$ and $|L_2|$ of $Q$ are interchanged by the
covering automorphism. The real structure on $\PP^2$ lifts to
two real structures $Q^+$ and $Q^-$ on the quadric $Q$ so that $\R Q^+$, homeomorphic to a $2$-sphere,
doubly covers the disc $\overline F^{\;o}$, and $\R Q^-$, homeomorphic to $(S^1)^2$, doubly covers the
M\"obius band $\overline F^{\;no}$. Each counted plane rational curve $C$ of degree $d$, quadratically tangent to $C$ at
$d$ points lifts to a pair of rational curves $C_1\in|iL_1+(d-i)L_2|$,
$C_2\in|(d-i)L_1+iL_2|$, $1\le i\le d-1$, interchanged by the covering automorphism.
They are not constrained by relative conditions, and $C_1$, resp. $C_2$ passes through one point in each
pair $\rho^{-1}(w)$, $w\in\bw$.
Thus, we get
 $$RW_m(\PP^2,C_2,F^o,0,dL)=\begin{cases}2^{2d-2-m}W_m(Q^+,\frac{d}{2}(L_1+L_2),\R Q^+,0),
 \ &d\equiv0\mod2,\\
 0,\ &d\equiv1\mod2,\end{cases}$$
$$RW_m(\PP^2,C_2,F^{no},0,dL)=2^{2d-2-m}\sum_{i=0}^dW_m(Q^-,iL_1+(d-i)L_2,\R Q^-,0)\ $$
where $W_m(*)$ are Welschinger invariants. A few particular values obtained in such a way are
presented
in
the table
below.

\begin{table}[ht]
\begin{center}
{\begin{tabular}{|c||c|c|c|c|c|c|}
\hline $d$  & $1$ & $2$ & $3$ & $4$ & $5$ &
$6$
  \\
\hline\hline $RW_0(\PP^2,C_2,F^o,0,dL)$ & $0$ & $4$ & $0$ & $384$ & $0$ & $589824$
  \\
\hline $RW_0(\PP^2,C_2,F^{no},0,dL)$  & $2$ & $4$ & $32$ & $640$ & $43008$ & $1523712$
  \\
\hline
\end{tabular}}\end{center}
\end{table}

Note
that, as it
follows from the properties of Welschinger invariants of quadrics
(see \cite[Theorem 2.2]{IKS2}),
one has
$RW_0(\PP^2,C_2,F^o,0,2kL)>0$ and $RW_0(\PP^2,C_2,F^{no},0,dL)>0$ for any positive integers $k$ and $d$;
furthermore,
$$
\displaylines{
\log RW_0(\PP^2,C_2,F^o,0,2kL) = 4k \log k + O(k), \cr
\log RW_0(\PP^2,C_2,F^{no},0,dL) = 2d \log d + O(d).
}
$$

We also exhibit a series of relative invariants that are not directly linked to
Welschinger invariants.
They
coincide with the so-called sided $w$-numbers
\cite[Section 3.8]{IKS2}, and they can be computed via the recursive formula in
\cite[Theorem 3.2(3)]{IKS2} and the initial values in \cite[Proposition 3.4]{IKS2}.

First,
we
consider the triple
$(Q,E,F^+)$,
where $Q\subset\PP^3$ is a real quadric surface with
$\R Q$ being a $2$-dimensional sphere,
$E$ is a smooth plane section with $\R E\simeq S^1$, and $F^+$ is one of the components of
$\R Q\setminus\R E$.
Here, we get the following values:
\begin{table}[ht]
\begin{center}
{\begin{tabular}{|c||c|c|c|c|}
\hline $D$
& $E$ & $2E$ &
$3E$ & $4E$
  \\
\hline\hline $RW_0(Q,E,F^+,0,D)$ & $2$
& $16$ & $1024$ &
$259584$
  \\
\hline
\end{tabular}}\end{center}
\end{table}

Next,
we compute
relative
invariants for the tuples $(\PP^2_{(0,k)},E,F^{no},-K_\Sigma)$, $k=1,2,3$,
where $\PP^2_{(0,k)}$ is the plane blown up at $k$ pairs of distinct complex conjugate
points on a real conic $C_2$ with
non-empty real part, $E$
is
the strict transform of $C_2$, and $F^{no}$
is
the non-orientable component of
$\R\PP^2_{(0,k)} \setminus \R E$
(for the orientable component $F^{o}$ of
$\R\PP^2_{(0,k)} \setminus \R E$
the invariants
vanish, since the real part of a real plane cubic cannot lie in the disc bounded
by a real conic):

\begin{table}[ht]
\begin{center}
{\begin{tabular}{|c||c|c|c|}
\hline $\Sig$
& $\PP^2_{(0,1)}$ & $\PP^2_{(0,2)}$ &
$\PP^2_{(0,3)}$
  \\
\hline\hline $RW_0(\Sig,E,F^{no},0,-K_\Sig)$ & $16$
& $8$ & $4$
  \\
\hline
\end{tabular}}\end{center}
\end{table}

\FloatBarrier

The following examples illustrate dependence
of relative invariants
on the
choice of the component $F^+$ of $F\setminus E$:
$${}$$

\begin{table}[ht]
\begin{center}
{\begin{tabular}{|c||c|c|}
\hline $F^+$
& $F^{o}$ & $F^{no}$
  \\
\hline\hline $RW_0(\PP^2_{(0,3)},E,F^+,0,-2K_\Sig-E)$ & $48$
& $80$
  \\
\hline
\end{tabular}}\end{center}
\end{table}

\medskip

The
last series of examples
concern nodal del Pezzo surfaces,
$\Sigma_+$ and $\Sigma_-$, of degree $2$
that are
obtained from $\PP^2_{(0,3)}$ by blowing up a real point lying in the
orientable or non-orientable component of
$\R\PP^2_{(0,3)} \setminus E$,
respectively. Notice that $\R\Sigma_+\setminus E$ consists of two components homeomorphic to the open
M\"obius band, and we denote one of them by $\widehat F$. In turn, $\R\Sigma_-\setminus E$ consists of
an open disc $\widehat F^o$ and a non-orientable component $\widehat F^{no}$.

\begin{table}[ht]
\begin{center}
{\begin{tabular}{|c||c|c|c|}
\hline $(\Sig,F^+)$
& $(\Sigma_+,\widehat F)$ & $(\Sigma_-,\widehat F^{o})$ &
$(\Sigma_-,\widehat F^{no})$
  \\
\hline\hline $RW_0(\Sig,E,F^+,0,-2K_\Sig-E)$ & $16$
& $8$ & $24$
  \\
\hline
\end{tabular}}\end{center}
\end{table}

\subsection{Lack of invariance for uninodal DP-pairs of degree $1$}\label{non-invariance}

The following examples demonstrate that the restrictions in Theorem \ref{t4p} (as compared to Theorem \ref{t2})
cannot be removed.

\smallskip

\subsubsection{A tangential DP-pair, $DE=2$}\label{ce-t-nodal}

\begin{figure}
\setlength{\unitlength}{1cm}
\begin{picture}(14.5,5)(0,-0.7)
\epsfxsize 145mm \epsfbox{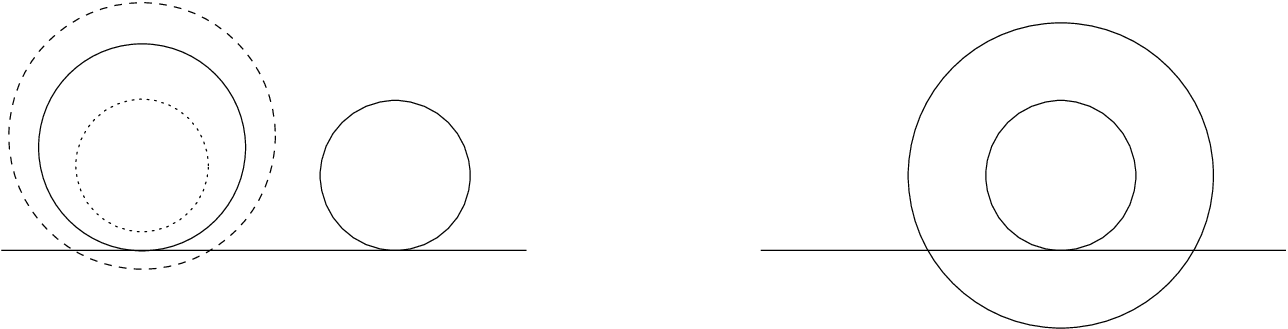}
\put(-11.5,-0.7){(a)}\put(-3,-0.7){(b)}
\put(-10,2.7){$C_2$}\put(-8.3,0.7){$C_1$}\put(-4.5,2.8){$C_2$}
\put(-3,2.8){$C'_2$}\put(-2.8,0.5){$z_0$}\put(-11.3,0.77){$\bullet$}
\put(-9.3,0.77){$\bullet$}\put(-5,0.77){$\bullet$}\put(-0.5,0.77){$\bullet$}
\put(-9.3,0.4){$z_8$}\put(-5,0.4){$z_1$}\put(-0.5,0.4){$z_2$}\put(-11.3,0.4){$z_7$}
\put(-12.8,2.4){$\bullet$}\put(-12.6,3){$\bullet$}\put(-12.4,3.4){$\bullet$}
\put(-13,2){$w^-$}\put(-12.9,2.8){$w$}\put(-12.4,3.7){$w^+$}
\end{picture}
\caption{Loss of invariance: Examples} \label{fig-example}
\end{figure}
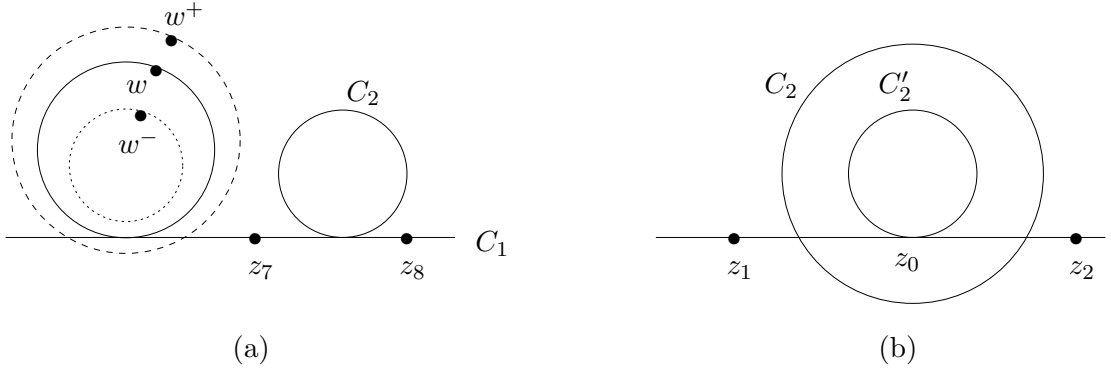

Let $C_2,C'_2\subset\PP^2$ be two real conics with $\R C_2\simeq\R C'_2\simeq S^1$ such that
they intersect in four non-real points $z_1,z_2,z_3,z_4$, and $\R C'_2$ lies in the non-orientable component
$F^+$ of $\R\PP^2\setminus\R C_2$.
Take a real line $C_1$ tangent both to $C_2$ and $C'_2$. In addition, we pick two generic complex conjugate points $z_5,z_6\in C_2\setminus C'_2$,
two generic points
$z_7,z_8\in \R C_1\setminus(C_2\cup C'_2)$, and a generic point $w\in\R C'_2$. There exist two close to $w$
points $w^+,w^-\in\R\PP^2$ such that the conic $C^{\prime,+}_2$ passing through $z_1,...,z_4$ and $w^+$ intersects $C_1$ in two
real points $q_1,q_2$, and the conic $C^{\prime,-}_2$ passing through $z_1,...,z_4$ and $w^-$ intersects $C_1$ in
a two complex conjugate points (see Figure \ref{fig-example}(a)), where
the conics $C'_2$, $C^{\prime,+}_2$, and $C^{\prime,-}_2$ are shown by the solid, dashed, and
dotted lines, respectively).

The blowing up $\pi:\Sig\to\PP^2$ at $z_1,...,z_8$ produces a
real tangential DP-pair $(\Sig,E)$ of degree $1$, with $E$ being the strict transform of $C_2$. For the divisor
$D=3L-E_1-E_2-E_3-E_4-E_7-E_8$, we have $DE=2$ and $r(\Sig,D,1)=1$.

\begin{lemma}\label{l22f}
One has the following equality:
$$RW_0(\Sig,E,F^+,0,D,w^-)-RW_0(\Sig,E,F^+,0,w^+)=2.$$
\end{lemma}

{\bf Proof.} Notice, first, that the generic choice of the points $z_7,z_8$ ensures that
there are no cuspidal curves in the linear system $|D|_\Sig$ passing through the point $w$.

Construct an elementary deformation of the pair $(\Sig,E)$ by
varying the point $z_8$ in the above construction along a germ of a real curve transversal to $C_1$.
Any member $(\widetilde\Sig,E)\ne(\Sig,E)$ of this deformation is not tangential, hence
satisfies the conditions of
Theorem
\ref{t4p}. There is a one-to-one correspondence between the set
${\mathcal V}_1(\widetilde\Sig,E,F^+,D,w^{\eps})$ and the set
\begin{equation}{\mathcal V}_1(\Sig,E,F^+,D,w^{\eps})\cup\left\{[\bn:\widehat C_1\cup\widehat C_2\to\PP^2]\ \big|\
\bn:\widehat C_1\overset{\sim}{\to}C_1,\ \bn:\widehat C_2\overset{\sim}{\to}C^{\prime\;,\eps}_2\right\}\ ,\label{e22g}
\end{equation}
for $\eps=\pm$. An injection of the former set to the latter one comes from Lemmas \ref{l22}(3) and \ref{l22e}.
To obtain an inverse map, we observe that any element of the set (\ref{e22g}) is a center of a unique
smooth branch of the family $\overline{\mathcal V}$ (in the notation of Lemma \ref{l22}). Indeed,
this follows from the facts that
the conditions to pass through the point $w^\varepsilon$
and to be tangent to $E$ are transversal for immersed rational curves in the linear system $|D|$,
and the conditions to pass through the point $z_7$
and to be tangent to $C_2$ are transversal for lines in $\PP^2$.

Hence,
the invariant
$$RW_0(\widetilde\Sig,E,F^+,0,D)=RW_0(\widetilde\Sig,E,F^+,0,D,w^+)=
RW_0(\widetilde\Sig,E,F^+,0,D,w^-)$$
on one side equals $RW_0(\Sig,E,F^+,0,D,w^-)$, and on the other side equals
$RW_0(\Sig,E,F^+,0,D,w^+)+2$, where the summand $2$ corresponds to the maps
$\bn:\widehat C_1\cup\widehat C_2\to\PP^2$ taking $\widehat C_1\simeq\PP^1$ isomorphically onto $C_1$, taking
$\widehat C_2$ isomorphically onto $C^{\prime\;,+}$, and projecting the intersection
point $\widehat C_1\cap\widehat C_2$ either to $q_1$, or to $q_2$.
\proofend

\smallskip

\subsubsection{A non-tangential uninodal DP-pair, $DE>2$}\label{ce-non-t-nodal}

Let $C_2\subset\PP^2$ be a real conic with $\R C_2\simeq S^1$.
Let $C_1\subset\PP^2$ be a real line intersecting $C_2$ at two real points, and let $C'_2$ be a real conic
tangent to $C_1$ at some point $z_0\in C_1\cap F^+$
and having the real part $\R C'_2\simeq S^1$ inside the open disc $F^+\subset\R\PP^2$ bounded by
$\R C_2$.
We assume that
$C_2\cap C'_2$ consists of four distinct imaginary points
(see Figure \ref{fig-example}(b)). Pick two generic real points $z_1,z_2\in C_1\setminus F^+$
and consider the family $V$ of plane rational quartics passing through $C_2\cap C'_2$,
having double points at $z_1,z_2$,
and
tangent to $C_2$ at some two points. One can easily extract from Lemma \ref{lp1-3} that $V$ is
a one-dimensional variety (if nonempty).

\begin{lemma}\label{lex1}
(1) The closure $\overline V$ contains the curve $C_4=C'_2+2C_1$, and the germ $(\overline V,C_4)$
contains a unique real branch $B$.

(2) The branch $B$ is smooth. For one component $\mathcal B$ of $\R B \setminus \{C_4\}$,
the double point of
any curve $C \in {\mathcal B}$ in $\R P^2\setminus\{z_1,z_2\}$
is solitary; for the other component of $\R B \setminus \{C_4\}$,
such a double point is non-solitary;
in the both cases, the double point is near $z_0$
{\rm (}{\it cf}. Figure \ref{fig-bif}(c){\rm )}.
\end{lemma}

{\bf Proof.} We construct a real branch $B\subset(\overline V,C_4)$ using the patchworking procedure described in
\cite[Section 5.3]{Sh0}. It starts with a flat family ${\mathfrak X}\to(\C,0)$ of surfaces
such that ${\mathfrak X}_t$ is smooth connected as $t\ne0$, and ${\mathfrak X}_0$ is reducible equipped with
a curve $C^{(0)}\subset{\mathfrak X}_0$. The patchworking extends $C^{(0)}$ to a flat family of curves $C^{(t)}
\subset{\mathfrak X}_t$,
$t\in(\C,0)$, possessing preassigned properties.

\begin{figure}
\setlength{\unitlength}{1cm}
\begin{picture}(10,7)(-2,0)
\thinlines

\put(1,1){\vector(1,0){3}}\put(1,1){\vector(0,1){5}}
\put(6,1){\vector(1,0){3}}\put(6,1){\vector(0,1){5}}

\thicklines

\put(3,1){\line(0,1){2}}\put(3,1){\line(-1,1){2}}
\put(1,3){\line(0,1){2}}\put(1,3){\line(1,0){2}}\put(3,3){\line(-1,1){2}}
\put(6,2){\line(0,1){3}}\put(6,2){\line(2,1){2}}
\put(6,4){\line(2,-1){2}}\put(6,5){\line(1,-1){2}}
\put(8,1){\line(0,1){2}}
\put(2.9,0.6){$2$}\put(0.6,2.8){$2$}
\put(0.6,4.8){$4$}\put(5.6,1.8){$1$}
\put(7.9,0.6){$2$}
\put(5.6,3.8){$3$}\put(5.6,4.8){$4$}
\put(1.5,3.5){$T_1$}\put(2.2,2.1){$T_2$}
\put(6.2,4){$T'_1$}\put(6.7,1.6){$T'_2$}
\put(6.5,2.8){$T_0$}

\put(2.2,0){(a)}\put(7.2,0){(b)}

\end{picture}
\caption{Newton polygons and their subdivisions in the proof of Lemma \ref{lex1}}\label{f1}
\end{figure}
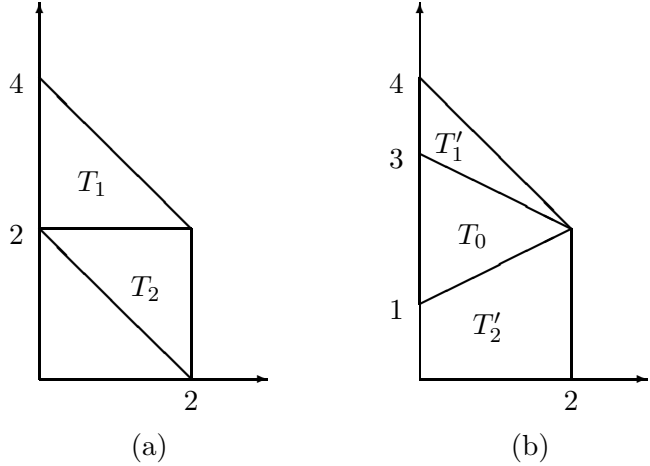

Introduce conjugation-invariant
homogeneous coordinates \mbox{$(x_0:x_1:x_2)$} in $\PP^2$ so that
\begin{eqnarray}&C_1=\{x_2=0\},\quad z_1=(1:0:0),\ z_2=(0:1:0),\ z_0=(1:q_0:0), \ q_0>0\nonumber\\
& C_1\cap C_2=\{(1:q_1:0),(1:q_2:0)\},\quad 0<q_1<q_0<q_2\ .\nonumber\end{eqnarray} In the affine coordinates $x=x_1/x_0$, $y=x_2/x_0$, the quartics $C\in V$
can be regarded as curves on the toric surface $\Tor(\Delta)$ with $\Delta=\conv\{(2,0),(0,2),(0,4),(2,2)\}$
(see Figure \ref{f1}(a)).
Consider the subdivision of $\Delta$ into two triangles $T_1$ and $T_2$
by the segment $\conv\{(0, 2), (2, 2)\}$ (see Figure \ref{f1}(a)).
The convex piecewise-linear function
$\nu:\Delta\to\R$, $\nu\big|_{T_1}(i,j)=0$, $\nu\big|_{T_2}(i,j)=4-2j$, defines
a family of surfaces ${\mathfrak X}\to(\C,0)$,
${\mathfrak X}_t\simeq\Tor(\Delta)$ as $t\ne0$, ${\mathfrak X}_0\simeq\Tor(T_1)\cup\Tor(T_2)$.
Note that $C_1$ and $C_2$, naturally mapped into ${\mathfrak X}_t$, $t\ne0$, respectively degenerate in ${\mathfrak X}_0$ into the
line $\Tor(T_1)\cap\Tor(T_2)$, and into the union of $C_2\subset\Tor(T_1)\simeq\PP^2$ with the lines $\{x=q_1\}$,
$\{x=q_2\}\subset\Tor(T_2)$.

Now we set $C^{(0)}\cap\Tor(T_1)$ to be defined by $y^2P_1(x,y)=0$, an equation of the curve $C'_2+2C_1$ in the plane,
whose truncation to the edge $T_1\cap T_2$ is $y^2(x-q_0)^2$, and we set $C''_2=C^{(0)}\cap\Tor(T_2)$ to be defined by
a polynomial $P_2(x,y)$, having Newton triangle $T_2$, truncation $y^2(x-q_0)^2$ on the edge
$T_1\cap T_2$, coefficient $1$ at $x^2$, and defining a reduced, irreducible curve, which
is simply tangent to the lines $x=q_1$ and $x=q_2$.
Notice that the vanishing of the discriminants of the
quadratic polynomials $P_2(q_1,y)$ and $P_2(q_2,y)$ gives four solutions for the pair of the
remaining coefficients of $P_2$ at $xy$ and $x_2y$. It is easy to check that all four solutions are real, two
of them correspond to the case of $P_2$ being an exact square,
and only one of the two remaining solutions satisfies the condition
\begin{equation}\lambda_1\lambda_2>0,\quad\text{where}\quad \lambda_1=\frac{\partial P_1}{\partial y}(q_0,0),\ \lambda_2=\frac{\partial P_2}{\partial y}(q_0,0)>0\ .
\label{eex1}\end{equation} Thus, we choose the latter solution for $P_2(x,y)$.

We want to deform $C^{(0)}$ in a conjugation invariant family of curves $C^{(t)}\in V$. In particular,
the point of simple tangency on the divisor $\Tor(T_1)\cap\Tor(T_2)$ deforms into a node. Such deformations are described
by deformation patterns in the sense of \cite[Section 3.5]{Sh0}: namely, the convex hull of the Newton polygons $T'_1,T'_2$ of the polynomials
$y^2P_1(x'+q_0,y)$, $P_2(x'+q_0,y)$, respectively, includes an additional triangle
$T_0=\conv\{(0,1),(0,3),(2,2)\}$ (see Figure \ref{f1}(b)). A deformation pattern is a nodal curve in
$\Tor(T_0)$ defined by a polynomial, whose coefficients of $y^3,y,(x')^2y^2$ induced by $y^2P_1$ and $P_2$ are $\lambda_1,\lambda_2,1$, respectively,
and the only extra coefficient is that of $y^2$; the nodality condition yields its value $\pm2\sqrt{\lambda_1\lambda_2}$. Note that the function $\nu$ induces a
convex piecewise-linear function $\nu':T'_1\cup T'_2\cup T_0\to\R$ such that
$\nu'\big|_{t'_1}=0$, $\nu'\big|_{T'_2}(i,j)=4-2j$, and $\nu'(0,0)=1$.

The existence and uniqueness of the desired deformation of $C^{(0)}$ follow from
\cite[Theorem 5]{Sh0}. The hypotheses of this patchworking theorem reduce to the following transversality
claims: \begin{itemize}\item in the germ at $C'_2$ of the space of conics, the conditions to pass through
the four points of $C_2\cap C'_2$ and be tangent to $C_1$ are smooth and transversal, and determine
the unique curve $C'_2$, \item in the germ at $P_2$ of the space of polynomials with Newton triangle $T_2$,
the conditions to have truncation $y^2(x-q_0)^2$ on the edge $T_1\cap T_2$ and coefficient $1$ at $x^2$, and to define
a curve tangent to the lines $\{x=q_1\}$ and $\{x=q_2\}$, are smooth and transversal, and determine the unique polynomial $P_2$,
\item the required transversality condition for the deformation patterns follows from
\cite[Lemma 5.5(i)]{Sh0}.
\end{itemize}
Furthermore, in the coordinates $(x',y)$ introduced above, by \cite[Formula (5.3.22)]{Sh0} the deformation
can be expressed as
$$C^{(t)}=\Big\{\sum_{(i,j)\in T'_1}(x')^iy^j(a_{ij}+O(t))+\sum_{(i,j)\in T'_2\setminus T'_1}t^{4-2j}(x')^iy^j(a_{ij}+O(t))$$ $$+
ty^2(2\sqrt{\lambda_1\lambda_2}+O(t))+x'y^2\cdot O(t)=0\Big\},\quad t\in(\C,0)\ .$$
It follows that
\begin{itemize}\item the Welschinger sign of the
node of $C^{(t)}$, $t\in(\R,0)\setminus\{0\}$, in a neighborhood of
the point $(q_0,0)$, is that of the node of the deformation pattern
$(x')^2y^2+\lambda_1y^3+\lambda_2y+\sign(t)\cdot2\sqrt{\lambda_1\lambda_2}y^2$, and hence it changes
when $t$ changes its sign;
\item the constructed local branch of $V$ at $C_4=C'_2+2E_0$ is smooth; thus, if its
(projective) tangent is spanned by $C_4$ and some curve $C'_4$, then $(V,C_4)$ is
conjugation-invariant diffeomorphic to the germ of a real line $L$ transversally crossing $C'_2$ at
a generic point $w_0\in\R C'_2\setminus C'_4$ via the map $C^{(t)}\in V\mapsto
C^{(t)}\cap (L,w_0)$.
\end{itemize}
Observe that the germ $(V,C_4)$ does not contain any other real branch. Indeed, let $B'$ be a real branch
of $(V,C_4)$,
regularly parameterized by $t\in(\C,0)$. Mapping curves $C\in B'\setminus\{C_4\}$ to their intersection
point with the germ of $C_2$ at
$z\in C_2\cap C_1$, we obtain a covering of an even degree, since, for real $t$, these intersection points are real, and they belong to the same
local real half-branch of $C_2$ at $z$, which lies in the component of
$F^+\setminus C_1$ that contains $\R C'_2$ (follows from the fact that the intersection of any curve
$C\in V\setminus\{C_4\}$ with $C_1$ is concentrated at the double points $z_1,z_2$).
Hence, in the above affine coordinates $(x,y)$, the spoken intersection points are
$(q_i+O(t),t^{2k}(\xi_i+O(t))$, $\xi_i\in\R^*$, $i=1,2$. We then inscribe the branch
$B'\to(\C,0)$ into the family of surfaces ${\mathfrak X}\to(\C,0)$, parameterized by
$\tau=t^{2k}$ and obtained from the trivial family $\PP^2\times(\C,0)$ by
blowing up the line $C_1\subset\PP^2\times\{0\}$. In the central fibre ${\mathfrak X}_0=\PP^2\cup
\F_1$, the constant family of curves $C_2$ degenerates into the union of $C_2\subset\PP^2$ and two
fibres $\{x=q_1\}$, $\{x=q_2\}$. Respectively, the central curve $C_0$ of $B'$ turns into
the union of $C'_2\subset\PP^2$ and a rational curve in $\F_1$ having double points at
$z_1.z_2$, tangent to the fibre $\{x=q_i\}$ at the point $(q_i,\xi_i)$, $i=1,2$, and tangent to
the $(-1)$-line at the point $x=q_0$. Thus, we get to the initial data of the construction of the
above real branch $B$. Note that a local deformation of the tacnode of $C_0$ of the line
$\PP^2\cap\F_1$ is described by gluing up a deformation pattern, which we mentioned above
(see \cite[Section 3.5, Lemma 3.10]{Sh0}), we conclude that $B'=B$.
\proofend

Now we blow up the points $z_1,z_2$, the four points of $C_2\cap C'_2$ and two more generic complex conjugate points of $C_2$, and obtain a real nodal del Pezzo pair $(\Sig,E)$ with $E$ being the strict transform of $C_2$. The family of quartics $\overline V$ turns into the family
$\overline{\mathcal V}_2(\Sig,E,D)$ with $D=4L-E_1-...-E_4-2E_7-2E_8$
(in the notation of (\ref{e5n})), and, finally, the statement of Lemma \ref{lex1} yields that
the number $RW_0(\Sig,E,F^+,\varphi,D,w)$ jumps by $\pm2$ as $w$ moves along
a smooth
real-analytic curve germ, transversally to the strict transform of $C'_2$ at a generic real point $w_0$.
Lemma \ref{lex1} can be easily generalized in order to produce the following non-invariance statement.

\begin{theorem}\label{t4}
Let $(\Sig, E)$ be a real non-tangential DP-pair of degree $1$ possessing property $T(1)$.
Let $F \subset \Sig$ be an admissible component,
$\varphi\in H_2(\Sig\setminus F,\Z/2)$ a $\conj$-invariant class,
and $D\in\Pic_+(\Sig,E)$ a real
divisor class matching conditions (\ref{e1})
and satisfying $r = -DK_{\Sig} - DE/2 - 1 > 0$.
Suppose that
$E_0\cap F\ne\emptyset$,
$DE=2l\ge4$, $\R E_0 \subset F$, and that
there exists a real rational curve $C\in|D-2E_0|$ such that
\begin{itemize}
\item $C$ has
$l-2$ local
branches centered on $E$, each branch intersecting $E$ with
multiplicity $2$, \item $C$ has a one-dimensional real branch in $F^+$,
\item $C$ intersects $E_0$ at $CE_0-1$ distinct points, at $CE_0-2$ of them transversally
and at one of them $C$ has a smooth branch
simply tangent to $E_0$.
\end{itemize}
Then, for any integer $0\le m \le r/2$,
the number $RW_m(\Sig,E,F^+,\varphi,D,\bw)$ does depend on
the choice of a generic tuple $\bw\in{\mathcal P}_{r,m}(\Sig,F^+)$.
\end{theorem}

Details of the proof
are left to the reader.

\medskip

{\bf Acknowledgements}.
A considerable part of this work was done during our visits to the {\it Max-Planck-Institut f\"ur Mathematik}, Bonn
and authors' Research in Pairs stay at the {\it Mathematisches Forschungsinstitut Oberwolfach}.
We thank MPIM and MFO for the hospitality and excellent working conditions.
The
third author has been supported in part from the Israeli Science
Foundation grants no. 448/09 and 176/15
, from the German-Israeli Foundation
grant no. 1174-197.6/2011, and from the Hermann-Minkowski-Minerva
Center for Geometry at the Tel Aviv University.

{\ncsc Institut de Math\'ematiques de Jussieu - Paris Rive Gauche\\[-21pt]

Sorbonne Universit\'e\\[-21pt]

4 place Jussieu,
75252 Paris Cedex 5,
France} \\[-21pt]

{\it and} {\ncsc D\'epartement de math\'ematiques et applications\\[-21pt]

Ecole Normale Sup\'erieure\\[-21pt]

45 rue d'Ulm, 75230 Paris Cedex 5, France} \\[-21pt]

{\it E-mail address}: {\ntt     ilia.itenberg@imj-prg.fr}

\vskip10pt

{\ncsc Universit\'e de Strasbourg et IRMA \\[-21pt]

7, rue Ren\'e Descartes, 67084 Strasbourg Cedex, France} \\[-21pt]

{\it E-mail address}: {\ntt kharlam@math.unistra.fr}

\vskip10pt

{\ncsc School of Mathematical Sciences \\[-21pt]

Raymond and Beverly Sackler Faculty of Exact Sciences\\[-21pt]

Tel Aviv University \\[-21pt]

Ramat Aviv, 69978 Tel Aviv, Israel} \\[-21pt]

{\it E-mail address}: {\ntt shustin@post.tau.ac.il}

\end{document}